\documentclass{amsart}
\usepackage{mathtools,bm,comment,amssymb,enumitem,hyperref,cite,xcolor,array,float,subcaption,float,pgfplots,tikz,filecontents}
\usepackage[capitalise,nameinlink]{cleveref}
\usepackage[english]{babel}
\usepackage[margin=1in]{geometry}
\usepackage[foot]{amsaddr}
\usepackage[linesnumbered,ruled]{algorithm2e}
\usepackage[noend]{algpseudocode}
\usepackage[export]{adjustbox}

\numberwithin{equation}{section}
\crefname{subsection}{Subsection}{Subsections}
\crefformat{equation}{(#2#1#3)}
\crefformat{enumi}{(#2#1#3)}
\crefname{figure}{Figure}{Figures}

\newcommand{\weak}{\rightharpoonup}
\newcommand{\longweak}{\relbar\joinrel\!\!\weak}
\newcommand{\eps}{\varepsilon}

\newcommand{\wdr}{\widetilde{R}}
\newcommand{\wdk}{\widetilde{K}_v}
\newcommand{\ECM}{{E\hspace{-.05em}C\hspace{-.15em}M}}
\newcommand{\MDE}{{M\hspace{-.15em}D\hspace{-.1em}E}}
\newcommand{\TAF}{{T\hspace{-.15em}A\hspace{-.1em}F}}
\newcommand{\ecm}{\phi_{\ECM}}
\newcommand{\taf}{\phi_{\TAF}}
\newcommand{\mde}{\phi_{\MDE}}
\newcommand{\con}{\hookrightarrow}
\newcommand{\com}{%
    \mathrel{\mathrlap{{\mspace{4mu}\lhook}}{\hookrightarrow}}}
\newcommand{\dd}{\textup{d}}
\newcommand{\dt}{\,\textup{d}t}
\newcommand{\ds}{\,\textup{d}s}
\newcommand{\dS}{\,\textup{d}S}
\newcommand{\ddt}{\frac{\dd}{\dd t}}

\newcommand{\p}{\partial}
\newcommand{\pt}{\p_t}
\newcommand{\ov}{\overline}

\renewcommand{\div}{\textup{div}}
\newcommand{\tr}{\textup{tr}}
\newcommand{\A}{\mathcal{A}}
\newcommand{\C}{\mathcal{C}}
\newcommand{\D}{\mathcal{D}}

\newcommand{\nablala}{\nabla_{\!\Lambda}}

\renewcommand{\L}{\mathcal{L}}
\renewcommand{\H}{\mathcal{H}}
\newcommand{\R}{\mathbb{R}}

\newcommand{\N}{\mathbb{N}}
\newcommand{\RHS}{\textup{RHS}}

\newcommand{\blue}{\textcolor{black}}

\newcommand{\CH}{{\mathcal{C}\mathcal{H}}}
\newcommand{\RD}{{\mathcal{R}\mathcal{D}}}
\newcommand{\phib}{\bm{\phi}}

\newcommand{\Lip}{\textup{Lip}}

\newtheorem{assumption}{Assumption}
\newtheorem{lemma}{Lemma}
\newtheorem{theorem}{Theorem}
\newtheorem{definition}{Definition}

\allowdisplaybreaks

\makeatletter
\@namedef{subjclassname@2010}{\textup{2010} Mathematics Subject Classification}
\makeatother

\makeatletter
\renewcommand{\email}[2][]{%
	\ifx\emails\@empty\relax\else{\g@addto@macro\emails{,\space}}\fi%
	\@ifnotempty{#1}{\g@addto@macro\emails{\textrm{(#1)}\space}}%
	\g@addto@macro\emails{#2}%
}
\makeatother

\allowdisplaybreaks
\title[Analysis of a new multispecies tumor growth model]{Analysis of a new multispecies tumor growth model coupling \\ 3D phase-fields with a 1D vascular network}
\author[M. Fritz, P. K. Jha, T. K\"oppl, J. T. Oden, and B. Wohlmuth]{Marvin Fritz$^{1,*}$, Prashant K. Jha$^2$, Tobias K\"oppl$^1$,  J. Tinsley Oden$^2$, and Barbara Wohlmuth$^1$}
\subjclass[2010]{35K35, 35A01, 35D30, 35Q92, 65M60.}
\keywords{tumor growth,  3D-1D coupled blood flow models, ECM degradation, existence of weak solutions, energy inequality, Galerkin method}
\thanks{${}^*$Corresponding author}
\email{\{fritzm, koepplto, wohlmuth\}@ma.tum.de}
\email{pjha@utexas.edu}
\email{oden@oden.utexas.edu}

\begin{document}
	\maketitle
	\vspace*{-2mm}
	\begin{center} \footnotesize
		$^1$Department of Mathematics, Technical University of Munich, Germany \\
		$^2$Oden Institute for Computational Engineering and Sciences, The University of Texas at Austin, USA
	\end{center}
	
	\vspace{1mm}
	\begin{abstract}
		In this work, we present and analyze a mathematical model for tumor growth incorporating ECM erosion, interstitial flow, and the effect of vascular flow and nutrient transport. The model is of phase-field or diffused-interface type in which multiple phases of cell species and other constituents are separated by smooth evolving interfaces. The model involves a mesoscale version of Darcy's law to capture the flow mechanism in the tissue matrix. Modeling flow and transport processes in the vasculature supplying the healthy and cancerous tissue, one-dimensional (1D) equations are considered. Since the models governing the transport and flow processes are defined together with cell species models on a three-dimensional (3D) domain, we obtain a 3D-1D coupled model. 
	\end{abstract}

\section{Introduction}
We develop and analyze a mathematical model of vascular tumor growth designed to simulate abstractions of many of the key phenomena known to be involved in the growth-decline of tumors and therapeutic treatment in living tissue. The complex vascular structure of tissue and the network of blood vessels supplying nutrients to a solid tumor mass embedded in the tissue are modeled as a network of one-dimensional capillaries within a three-dimensional tissue domain, while the growth of the tumor is represented by a phase-field model involving multiple cell species and other constituents. Our tumor models may be regarded as mesoscale depictions of physical and biological events employing continuum mixture theory to construct general forms of the Ginzburg--Landau--Helmholtz free energy of biological materials in terms of volume fractions or mass concentrations of the cell phenotypes and principal mechanical and chemical fields. The equations governing the tumor growth are derived from the balance laws of continuum mixture theory as in e.g. \cite{byrne2003modelling, cristini2009nonlinear, oden2016toward, oden2010general, lima2014hybrid}, and representations of the principal mechanisms governing the development and evolution of cancer \cite{lima2014hybrid, hanahan2011hallmarks}. In the tissue containing the tumor cells, the microvascular network is represented by a graph structure with 1D filaments through which nutrient-containing blood may flow. The exchange of nutrients between the network and tissue is depicted by a Kedem--Katalchsky type law \cite{ginzburg1963frictional}.  We briefly describe the construction of approximations of these models, see also \cite{garcke2016global, garcke2017well, cristini2010multiscale, lima2017selection, fritz2019local, fritz2018unsteady}.

There is a significant and growing volume of published work on various aspects of this subject. Continuum mixture theory as a framework for developing meaningful models of materials with many interacting constituents is proposed in \cite{oden2016toward, oden2010general, lima2014hybrid, hawkins2012numerical, garcke2018multiphase, cristini2010multiscale}.
Of particular interest are the comprehensive developments of diffuse-interface multispecies models described in \cite{wise2008three, frieboes2010three}, the four- and ten species models presented in \cite{hawkins2012numerical, lima2014hybrid}, and the multispecies nonlocal models of adhesion and tumor invasion described in \cite{fritz2019local}. The book compiled by Lowengrub and Cristini \cite{cristini2010multiscale} contains over 700 references to relevant cancer cell biology and mathematical models of cancer growth. 
The complex processes underlying angiogenesis which are key to vascular tumor growth present formidable challenges to the goal of predictive computer modeling. Angiogenesis models embedded in models of hypoxic and cell growth or decline were presented in \cite{lima2014hybrid, xu2016mathematical, xu2017full, xu2020phase, wise2008three, santagiuliana2016simulation, santagiuliana2019coupling}. More recent developments have included models of the vascular network interwoven in tissue containing solid tumors, and the sprouting of capillaries in response to concentrations of various tumor angiogenesis factors so as to supply nutrients to hypoxic tumor cells. Such network-tissue models are discussed in \cite{xu2016mathematical,  kremheller2019approach, koppl2018mathematical}. These models generalized the lattice-probabilistic network models of \cite{anderson1998continuous}.

This article is organized as follows: In \cref{Sec:Derivation}, we introduce various components of the complete model, such as the tissue domain, the 1D network domain, the species in the multi-species phase-field model. Further, we present the governing partial differential equations. The resulting model is a highly non-linear coupled system of partial differential equations. We give some analytical preliminaries in \cref{Sec:Preliminaries}, e.g. Sobolev embeddings and interpolation inequalities in Bochner spaces, which will be used in the following sections. In \cref{Sec:Theorem}, we state a theorem for the existence of weak solutions of the coupled non-linear 3D-1D model under certain given assumptions. In \cref{Sec:Proof}, we give the proof of the theorem via the Faedo--Galerkin approximation and compactness methods. 

\section{Derivation of the model} 
\label{Sec:Derivation}

\subsection{Setup and notation}

We consider a region of vascularized tissue in a living subject, e.g., within an organ, which is host to a colony of tumor cells and other constituents that make up the so-called microenvironment of a solid tumor. The tumor is contained in an open bounded domain $\Omega \subset \R^3$ and is supported by a network of macromolecules within $\Omega$ consisting of collagen, enzymes, and various proteins, that constitute the extracellular matrix (ECM). We focus on developing phenomenological characterizations of the evolutions of the tumor cell colony that attempt to capture mesoscale and macroscale events. 

The primary feature of our model of tumor growth is that it employs the framework of continuum mixture theory in which multiple mechanical and chemical species can exist at a point $x \in \Omega$ at time $t>0$. Thus, for a \blue{medium} with $N$ interacting constituents, the volume fraction of each species $\phi_\alpha$, $1\leq \alpha \leq N$, is represented by a field $\phi_\alpha$ with value $\phi_\alpha(t,x)$ at $x\in \Omega$, and time $t\geq 0$, and $\sum_{\alpha} \phi_\alpha(t,x) = 1$. Setting $\alpha=1=T$, the volume fraction of tumor cells $\phi_T(t,x)$ is understood to represent an averaged cell concentration, a homogenized depiction over many thousands of cells, since in volumes as small as a voxel in modern tumor imaging techniques, $4-5 \times 10^4$ cells can exist.

 We could also develop equivalent models in terms of mass concentration, $c_\alpha = \rho_\alpha \phi_\alpha$, $\rho_\alpha$ being the mass density of species $\alpha$. Moreover, we assume that $\rho_\alpha = \rho_0$ = constant, $1\leq \alpha\leq N$, and thus, $C_\alpha$ and $\phi_\alpha$ are up to a fixed scaling equivalent. This simplification is regarded as a reasonable assumption in many investigations since the mass densities of species are generally close to that of water at room temperature. 

As another important feature of our model, we depict the evolving interfaces in which a smooth boundary layer exists and which is defined intrinsically as a feature of the solution of the forward problem. This feature is a property of phase-field or diffuse-interface models and avoids complex interface tracking while producing characterizations of interfaces between cell species which are in good agreement with actual observations. 

\begin{figure}[h]
    \centering
    \begin{subfigure}[t]{0.4\textwidth}
        \centering
        \includegraphics[page=1,width=\textwidth]{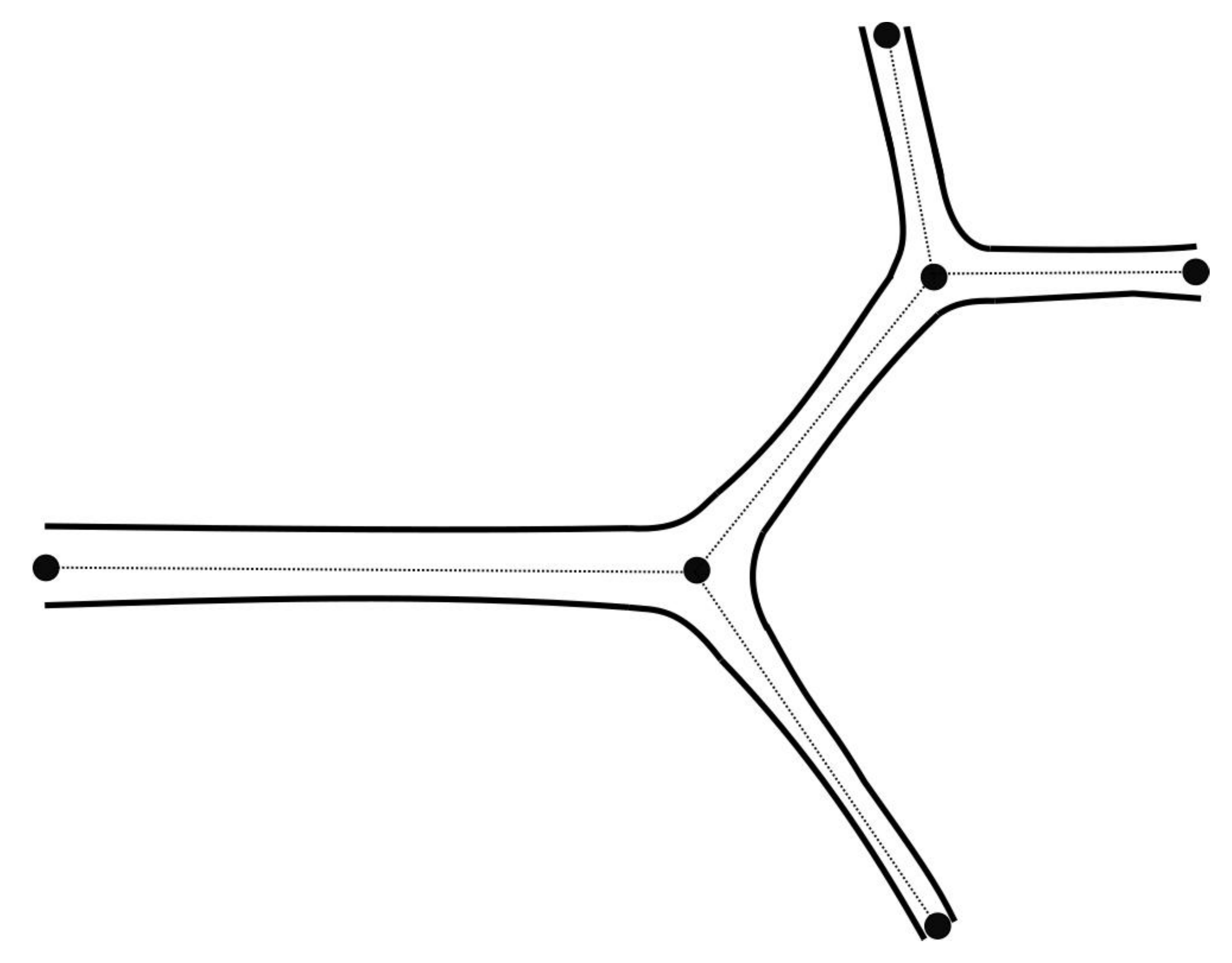}
        \caption{Outline of a blood vessel network}
	\end{subfigure}
	\quad
	\begin{subfigure}[t]{0.4\textwidth}
        \centering
        \includegraphics[page=2,width=\textwidth]{derivation3.pdf}
        \caption{Approximation of the blood vessels by cylinders with a constant radius}
	\end{subfigure} \\
	\begin{subfigure}[t]{0.4\textwidth}
        \centering
        \includegraphics[page=3,width=\textwidth]{derivation3.pdf}
        \caption{\label{fig:1DNetworkc}Reduction to a 1D graph-like structure}
	\end{subfigure}
	\quad
	\begin{subfigure}[t]{0.4\textwidth}
        \centering
        \includegraphics[page=4,width=\textwidth]{derivation3.pdf}
        \caption{\label{fig:1DNetworkd}Notation related to a single vessel}
	\end{subfigure} 
    \caption{\label{fig:1DNetwork}Modeling a blood vessel network \textsc{(a)} by means of a 1D graph-like structure \textsc{(c)}. At first the surface of the blood vessels is approximated by cylinders with constant radius whose surfaces are denoted by $\Gamma_i$, see \textsc{(b)}. Then, the blood vessels are lumped to the center lines $\Lambda_i$ of the cylinders.}
\end{figure}

Moreover, we consider a one-dimensional graph-like structure $\Lambda$ inside of $\Omega$ forming a microvascular network. The single edges of $\Lambda$ are denoted by $\Lambda_i$ such that $\Lambda$ is given by
$
\Lambda = \bigcup_{i=1}^N \Lambda_i.
$
The edge $\Lambda_i$ is parameterized by a curve parameter $s_i$, such that $\Lambda_i$ is given by:
$$
\Lambda_i = \left\{ x \in \Omega \left|\; x = \Lambda_i ( s_i ) = x_{i,1} + s_i \cdot ( x_{i,2}-x_{i,1} ),\;s_i \in (0,1 )  \right. \right\}.
$$
Thereby, $x_{i,1} \in \Omega$ and $x_{i,2} \in \Omega$ mark the boundary nodes of $\Lambda_i$, see \cref{fig:1DNetwork}. For the total 1D network $\Lambda$, we introduce a global curve parameter $s$, which has to be interpreted in the following way: $s=s_i$, if $x = \Lambda ( s ) = \Lambda_i (s_i )$. At each value of the curve parameter $s$, we study 1D constituents, which couple to their respective 3D counter-part in $\Omega$. In order to formulate the coupling between 3D and 1D constituents in \cref{sec:3Dmodel} and \cref{sec:1Dmodel}, we need to introduce the surface $\Gamma$ of the microvascular network. For simplicity, it is assumed that the surface for a single vessel is approximated by a cylinder with a constant radius, see \cref{fig:1DNetwork}. The radius of a vessel that is associated with $\Lambda_i$, is given by $R_i$ and the corresponding surface is denoted by $\Gamma_i$. \blue{In fact, $\Gamma_i$ is the surface of the cylinder whose center line is given by $\Lambda_i$, i.e.,}
$$
\blue{\Gamma_i = \left\{ x \in \Omega \left|\; \text{dist}( x, \Lambda_i ( s_i ) ) = R_i, \; s_i \in ( 0,1 ) \right. \right\}.}
$$
According to the definition of $\Lambda$, the total surface $\Gamma$ is given by the union of the single vessel surfaces,  i.e.,  
$
\Gamma = \bigcup_{i=1}^N \Gamma_i.
$

\begin{figure}[h]
	\centering
    \begin{subfigure}[t]{0.5\textwidth}
        \centering
        \includegraphics[width=\textwidth]{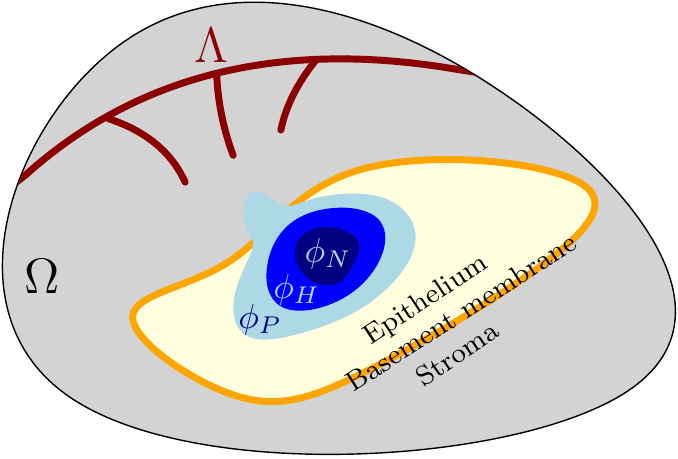}
        \caption{Sketch of the domain $\Omega$ containing a 1D network $\Lambda$.}
	\end{subfigure}
	\quad
	\begin{subfigure}[t]{0.38\textwidth}
        \centering
        \includegraphics[width=\textwidth]{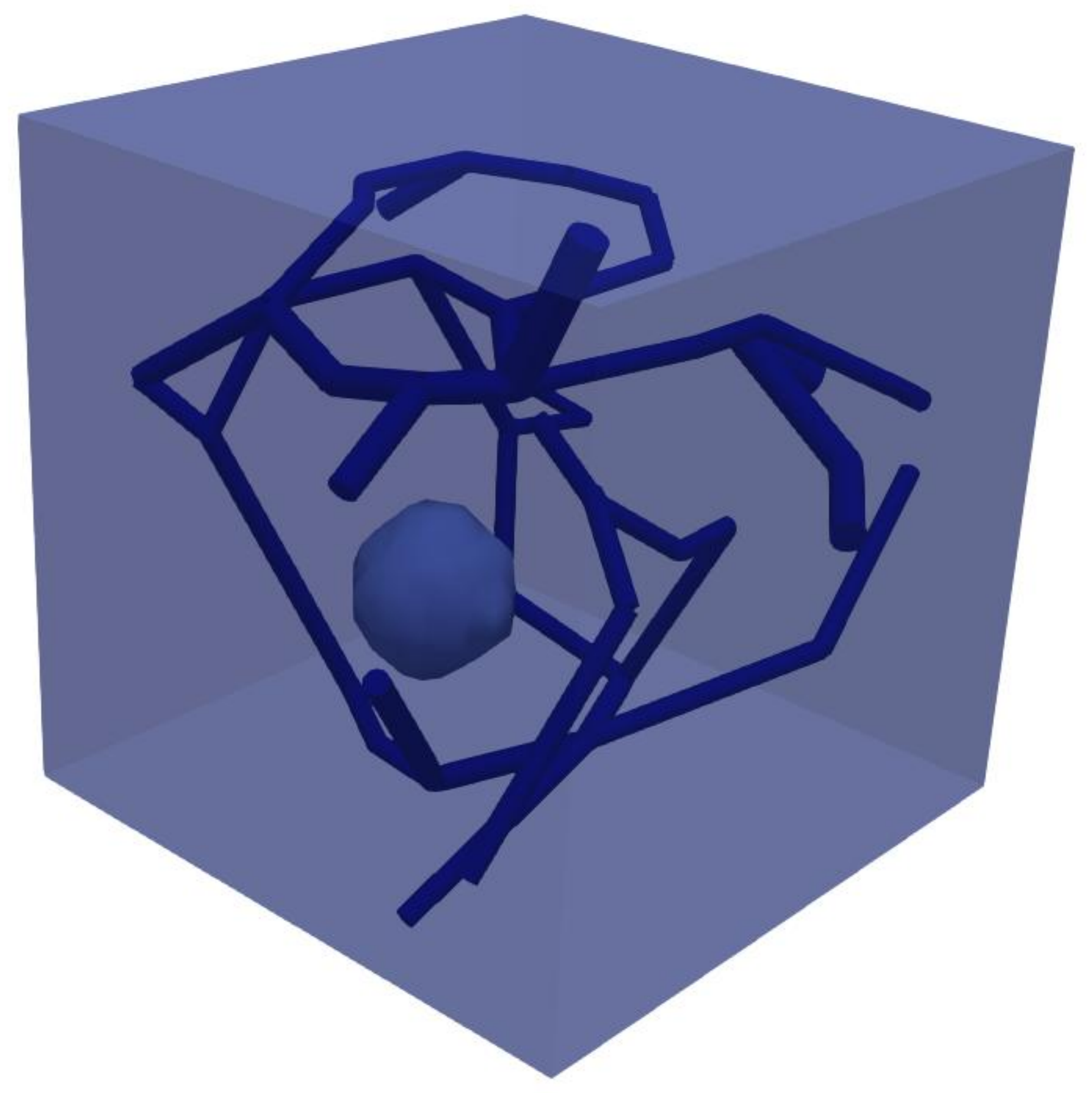}
        \caption{Tumor core surrounded by a network}
	\end{subfigure} \\
	\caption{Setup of the domain $\Omega$ with the microvascular network $\Lambda=\cup \Lambda_i$ and a tumor mass, which is composed in its proliferative ($\phi_P$), hypoxic ($\phi_H$) and necrotic ($\phi_N$) phase (left). Three dimensional presentation of a given tumor core surrounded by a capillary network (right).}
\end{figure}

\subsection{Constituents}
After introducing the domains on which the 1D and 3D models are defined, we describe in a next step all the dependent variables occurring in our model.

The tumor cell's field, $\phi_T=\phi_T(t,x)$, can be represented as the sum of three components, $\phi_T = \phi_P + \phi_H + \phi_N$, where $\phi_P=\phi_P(t,x)$ is the volume fraction of proliferative cells, $\phi_H=\phi_H(t,x)$ that of hypoxic cells, and $\phi_N=\phi_N(t,x)$ is the volume fraction of necrotic cells. Proliferative cells are those which have a high probability of mitosis, division into twin cells, and to produce growth of tumor. Hypoxic cells are those tumor cells deprived of sufficient nutrient (e.g., oxygen) to become or remain proliferative and necrotic cells have died due to the lack of nutrients. The local nutrient concentration is represented by a field $\phi_\sigma = \phi_\sigma(t,x)$. The tumor cells response to hypoxia (e.g., low oxygen), i.e., $\phi_\sigma$ is below a certain threshold, by the production of an enzyme (hypoxia-inducible factor) that accumulates and increases cell mobility and activates the secretion of angiogenesis promoting factors characterized by another field, $\taf=\taf(t,x)$, tumor angiogenesis factor. Of several such factors, that most frequently addressed, is VEGF, Vascular Endothelial Growth Factor, which induces sprouting of endothelial cells forming the tubular structure of blood vessels, the lumins, which grow into new vessels that supply nutrient to the hypoxic cells. In this article, we treat a stationary network of endothelial cells and neglect the sprouting.

Moreover, at lower oxygen levels the hypoxic cells release matrix-degenerative enzymes such urokinase-plasminogen and matrix metalloproteinases, labelled MDEs, with volume fraction denoted by $\mde=\mde(t,x)$, that can erode the extracellular matrix, whose density is denoted by $\ecm=\ecm(t,x)$, and make room for invasion of tumor cells, increasing $\phi_T$ in the ECM domain and increasing the likelihood of metastasis. Below a certain level of nutrient, or sustained periods of hypoxia, cells may die and enter the necrotic phase represented by the field $\phi_N$. In many forms of cancer, necrotic cells undergo calcification and become inert and can be removed as waste from the organism. 


On the one-dimensional network $\Lambda$, we consider the constituents $\phi_v=\phi_v(t,s)$ and $v_v=v_v(t,s)$, which represent the one-dimensional counter-part of the local nutrient concentration $\phi_\sigma$ and the volume-averaged velocity $v$. In addition, we consider both in the vascular system and the tissue domain pressure variables that are denoted by $p_v$ and $p$, respectively. The different constituents are coupled by the source terms of the different partial differential equations governing the behavior of the constituents. 

For convenience, we collect the constituents within the following 7-tuple:
$$
\phib  = ( \phi_P, \phi_H, \phi_N, \phi_\sigma, \mde, \taf, \ecm  ) = (\phi_\alpha)_{\alpha \in \A},
$$
where $\A=\{P,H,N,\sigma,\ECM,\MDE,\TAF\}$, and further, we distinguish between the tumor phase-field indices $\CH=\{P,H,N\}$, the reaction-diffusion indices $\RD=\{\sigma,\MDE,\TAF\}$ and the evolution index $\{\ECM\}$, which corresponds to an abstract ordinary differential equation.

\subsection{Three-dimensional model}
\label{sec:3Dmodel}

The constituents $\phi_\alpha$, $\alpha \in \A$, are \blue{governed by the following mass balance law}, see e.g., \cite{lima2014hybrid,lima2015analysis},
\begin{equation} \label{Eq:MassBalance}
	\partial_t \phi_\alpha+\text{div}(\phi_\alpha v_\alpha)=- \text{div} J_\alpha(\phib) +S_\alpha(\phib),
\end{equation}
for all $\alpha \in \A$, where $v_\alpha$ is the cell velocity of the $\alpha$-th constituent, and $S_\alpha$ describes a mass source term depending on all species $\phib$. Moreover, $J_\alpha$ denotes the flux of the $\alpha$-th constituent, which is given by
\begin{equation} \label{Eq:Flux}
	J_\alpha(\phib) = - m_\alpha(\phib) \nabla \mu_\alpha.
\end{equation}
Here, $\mu_\alpha$ denotes the chemical potential of the $\alpha$-th species and $m_\alpha$ the mobility function of it. In our applications, we consider the mobilities
$$\begin{aligned} m_\alpha(\phib) &=M_\alpha \phi_\alpha^2 (1-\phi_\alpha)^2 I_d, && \alpha \in \CH, \\ m_\beta(\phib) &= M_\beta I_d, && \beta \in \RD, \\
	m_{\ECM}(\phib) &=0, &&
\end{aligned}$$ where $M_\alpha$ are mobility constants and $I_d$ is the $(d\times d)$-dimensional identity matrix.  Especially, we choose $m_{\ECM}=0$ in accordance to the non-diffusivity of the ECM, see \cite{nargis2016effects}.  Following \cite{hawkins2012numerical,lima2015analysis,lima2014hybrid,wise2008three}, we define the chemical potential as
	$$\mu_\alpha=\frac{\delta \mathcal{E}(\phib)}{\delta \phi_\alpha},$$
where $\delta \mathcal{E}/\delta \phi_\alpha$ denotes the first variation (G\^{a}teaux derivative) of the Ginzburg--Landau--Helmholtz free energy functional,
\begin{equation}
	\mathcal{E}(\phib)= \int_\Omega \Big\{ \Psi(\phi_P, \phi_H, \phi_N) + \sum_{\alpha \in \CH}  \frac{\eps_\alpha^2}{2} |\nabla \phi_\alpha|^2  + \sum_{\beta \in \RD} \frac{D_\beta}{2} \phi_\beta^2 -  (\chi_c \phi_\sigma+\chi_h \ecm) \!\!\sum_{\alpha \in \{P,H\}}\!\! \phi_\alpha \Big\} \text{ d}x.
	\label{Eq:GinzburgLandau}
\end{equation}
Here, $\chi_c$ is the chemotaxis parameter, see \cite{hillen2013convergence},  $\chi_h$ represents the haptotaxis parameter, see \cite{fritz2019local,tao2011global}, and $\varepsilon_\alpha$, $\alpha \in \CH$, is a parameter associated with the interface thickness separating the different cell species. Lastly, $\Psi$ represents a double-well potential, e.g., \blue{it can be of Landau type, where we mention the three possibilities} $$\begin{aligned} \Psi(\phi_P,\phi_H,\phi_N)&=C_{\Psi_T} \phi_T^2(1-\phi_T)^2, \\
\Psi(\phi_P,\phi_H,\phi_N) &= C_{\Psi_P} \phi_P^2(1-\phi_P)^2 + C_{\Psi_H} \phi_H^2(1-\phi_H)^2 + C_{\Psi_N} \phi_N^2(1-\phi_N)^2 + C_{\Psi_T} \phi_T^2(1-\phi_T)^2, \\
\Psi(\phi_P,\phi_H,\phi_N) &= C_{\Psi_P} \phi_P^2(1-\phi_T)^2 + C_{\Psi_H} \phi_H^2(1-\phi_T)^2 + C_{\Psi_N} \phi_N^2(1-\phi_T)^2, \\
\end{aligned}$$ 
where $C_{\Psi_\alpha}$ are appropriate prefactors. \blue{Alternatively}, one can also select a logarithmic potential of Flory--Huggins type, e.g., see \cite{frigeri2018on, cherfils2011cahn},
$$\begin{aligned}\Psi(\phi_P,\phi_H,\phi_N) &=C_{\Psi_P} \phi_P\log\phi_P + C_{\Psi_H}\phi_H \log\phi_H +C_{\Psi_N}\phi_N \log \phi_N + C_{\Psi_T} (1-\phi_T) \log (1-\phi_T) \\
&\qquad + \frac{1}{2} \left( C_{\Psi_P}\phi_P(1-\phi_P) + C_{\Psi_H}\phi_H(1-\phi_H) + C_{\Psi_N}\phi_N(1-\phi_N) +C_{\Psi_T} \phi_T(1-\phi_T) \right).\end{aligned}$$
Lastly, we also mention potentials, which are used for abstract multiphase models, see \cite{garcke2018multiphase},  $$\Psi(\phi_P,\phi_H,\phi_N)= C_{\Psi_P} \phi_P^2 \phi_H^2  + C_{\Psi_H} \phi_H^2 \phi_N^2 + C_{\Psi_N} \phi_N^2 \phi_P^2 ,$$

The chemical potentials read
\begin{equation} \label{Eq:ChemicalPotential}
\begin{aligned}
	\mu_\alpha &= \partial_{\phi_\alpha} \Psi(\phi_P,\phi_H,\phi_N) - \eps^2_\alpha \Delta \phi_\alpha - \chi_c \phi_\sigma-\chi_h \ecm, && \alpha \in \CH \backslash \{N\}, \\
	\mu_\beta &= D_\beta \phi_\beta, && \beta \in \RD \backslash \{\sigma\}, \\
	\mu_N &= \partial_{\phi_N} \Psi(\phi_P,\phi_H,\phi_N) - \eps^2_\alpha \Delta \phi_N, &&  \\
	\mu_{\sigma} &= D_\sigma \phi_\sigma - \chi_c (\phi_P + \phi_H ),   && \\
	\mu_{\ECM} &= - \chi_h (\phi_P + \phi_H ).   && \\
	\end{aligned}
\end{equation}

The necrotic cells \blue{are non-moving} and only gain mass from the nutrient-lacking hypoxic cells. \blue{Therefore, the mobility of the necrotic cells is set to zero. Consequently, we have $m_N=v_N=0$.} Consequently, inserting \cref{Eq:Flux} and \cref{Eq:ChemicalPotential} into the mass balance equation \cref{Eq:MassBalance},  we arrive at the equations for $(\phi_\alpha)_{\alpha \in \CH}$ 
\begin{equation} \label{Eq:DerivationCH} \begin{aligned}
		\pt \phi_P+ \div(\phi_P v) & = \div (m_P( \phib) \nabla \mu_P) + S_P( \phib),    \\
		\mu_P                      & =   \partial_{\phi_P} \Psi(\phi_P,\phi_H,\phi_N) - \eps^2_P \Delta \phi_P - \chi_c \phi_\sigma-\chi_h \ecm, \\
		\pt \phi_H+ \div(\phi_H v) & = \div (m_H( \phib) \nabla \mu_H)+S_H( \phib),      \\
		\mu_H                      & =   \partial_{\phi_H} \Psi(\phi_P,\phi_H,\phi_N) - \eps^2_H \Delta \phi_H - \chi_c \phi_\sigma-\chi_h \ecm, \\
		\pt \phi_N                 & = S_N( \phib).
	\end{aligned} \end{equation}
Further, we propose the source functions
\begin{equation} \label{Eq:DerivationSourceCH}
	\begin{aligned}
S_P(\phib) &= \lambda_P \phi_\sigma \phi_P(1- \phi_T) - \lambda_A \phi_P - \lambda_{P\!H}  \H(\sigma_{P\!H} - \phi_\sigma)\phi_P + \lambda_{H\!P} \H(\phi_\sigma - \sigma_{H\!P}) \phi_H,  \\
S_H(\phib) &=  \lambda_{P_h} \phi_\sigma \phi_H (1-\phi_T)-\lambda_{A_h} \phi_H + \lambda_{P\!H}  \H(\sigma_{P\!H} - \phi_\sigma) \phi_P - \lambda_{H\!P} \H(\phi_\sigma - \sigma_{H\!P}) \phi_H \\&\qquad - \lambda_{H\!N} \H(\sigma_{H\!N} - \phi_\sigma)\phi_H,  \\
S_N(\phib) &= \lambda_{H\!N} \H(\sigma_{H\!N} - \phi_\sigma) \phi_H.
	\end{aligned}
\end{equation}
In \cref{Eq:DerivationCH}, $v=v_\alpha$ is a volume-averaged velocity for the fields $\phi_P$ and $\phi_H$. In \cref{Eq:DerivationSourceCH}, $\lambda_P$ is the rate of cellular mitosis of tumor cells, $\lambda_A$ and $\lambda_{A_h}$ are the apoptosis rates of the proliferative and hypoxic cells, respectively, $\lambda_{P_h}$ is the proliferation rate of hypoxic cells, $\lambda_{P\!H}$ the transition rate from the proliferative to the hypoxic phase below the nutrient level $\sigma_{P\!H}$, $\lambda_{H\!P}$ the transition rate from the hypoxic to the proliferative phase above the nutrient level $\sigma_{H\!P}$, and $\lambda_{H\!N}$ the transition rate from the hypoxic to the necrotic phase below the nutrient level $\sigma_{H\!N}$. Finally, $\H$ denotes the Heaviside step function.

Related models of extracellular matrix (ECM) degradation due to matrix-degenerative enzymes (MDEs) released by hypoxic cell concentrations and subsequent tumor invasion and metastasis are discussed in \cite{chaplain2005mathematical, chaplain2011mathematical, tao2011chemotaxis, tao2014energy, engwer2017structured,sfakianakis2020hybrid}. Following these references, we introduce the equation for the ECM evolution, 
	\begin{equation}
	\label{Eq:DerivationECM} 	
	\begin{aligned}
	\pt \ecm &= S_{\ECM}( \phib) \\
	&=-\lambda_{\ECM_{\!D}}\ecm\mde + \lambda_{\ECM_{\!P}} \phi_\sigma (1- \ecm) \H(\ecm - \phi_{\ECM_{\!P}}),
	\end{aligned}   
	\end{equation}
where $\lambda_{\ECM_{\!D}}$ is the degradation rate of ECM fibers due to the matrix degrading enzymes, and $\lambda_{\ECM_{\!P}}$ is the production rate of ECM fibers above the threshold level $\phi_{\ECM_{\!P}}$ for the ECM density.

Further, for $(\phi_\beta)_{\beta \in \RD}$ we arrive at the following system of equations
\begin{equation} \label{Eq:DerivationRD} \begin{aligned}
	\pt \phi_\sigma + \div(\phi_\sigma v) & = \div (m_\sigma( \phib) (D_\sigma \nabla \phi_\sigma - \chi_c  \nabla (\phi_P +\phi_H ))+S_\sigma( \phib) + S_{\sigma v}(\phi_\sigma,p,\phi_v,p_v), \\
	\pt \mde                        & = \div (m_{\MDE}( \phib) D_{\MDE} \nabla \mde)+S_{\MDE}( \phib),
	\\
	\pt \taf                        & = \div (m_{\TAF}( \phib) D_{\TAF} \nabla \taf)+S_{\TAF}( \phib),                        
\end{aligned} \end{equation}
with source functions
\begin{equation} \label{Eq:DerivationSourceRD}
	\begin{aligned}
S_{\sigma}(\phib) &= -\lambda_P \phi_\sigma\phi_P(1-\phi_T) -\lambda_{P_h} \phi_\sigma\phi_H(1-\phi_T) + \lambda_A \phi_P + \lambda_{A_h} \phi_H + \lambda_{\ECM_{\!D}}\ecm\mde \\
&\quad - \lambda_{\ECM_{\!P}} \phi_\sigma (1- \ecm) \H(\ecm - \phi_{\ECM_{\!P}}) ,  \\
S_{\MDE}(\phib) &= -\lambda_{\MDE_{\!D}}\mde + \lambda_{\MDE_{\!P}}(\phi_P + \phi_H) \ecm\frac{\sigma_{H\!P}}{\sigma_{H\!P} + \phi_\sigma} (1-\mde) - \lambda_{\ECM_{\!D}} \ecm \mde,  \\
S_{\TAF}(\phib) &= \lambda_{\TAF_{\!P}} (1- \taf) \phi_H \H(\phi_H-\phi_{H_{\!P}}) - \lambda_{\TAF_{\!D}} \taf.  
\end{aligned}
\end{equation}
Here, $\lambda_{\MDE_{\!D}}$ and $\lambda_{\TAF_{\!D}}$ denote the decay rates of the MDEs and TAFs, respectively, $\lambda_{\MDE_{\!P}}$ the production rate of MDEs, and $\lambda_{\TAF_{\!P}}$ is the production rate of the $\taf$ due to the release by hypoxic cells above a threshold value of $\phi_{H_{\!P}}$. \blue{We note that the cell species $\phi_\alpha$, $\alpha \in \{P,H,N,\sigma,\ECM\}$, form a mass conserving subsystem in the sense that their source terms add to zero. The fields $\phi_\MDE$ and $\phi_\TAF$ do not belong to this mass exchanging closed subsystem system since they show natural degradation factors.}

Additionally, we have introduced a source term $S_{\sigma v}$ in  \cref{Eq:DerivationRD} for the nutrient volume fraction $\phi_\sigma$, which depends on the 1D constituents $\phi_v$ and $p_v$, and therefore, this source term is responsible for the coupling between the constituents in $\Omega$ and $\Lambda$. In particular, it governs the exchange of nutrients between the vascular network and the tissue. In order to quantify the flux of nutrients across the vessel surface, we use the Kedem--Katchalsky law, see e.g., \cite{ginzburg1963frictional},
\begin{equation}
\label{eq:KedemKatalchsky}
J_{\sigma v} (\ov \phi_\sigma, \ov p, \phi_v, p_v ) = ( 1-r_{\sigma} ) J_{pv} ( \ov p,p_v )  \phi_{\sigma}^v + L_\sigma ( \phi_v - \ov \phi_\sigma ),
\end{equation}
where $J_{\sigma v}$ represents the flux of nutrients between the vascular network and the tissue. The Kedem--Katchalsky law \cref{eq:KedemKatalchsky} consists of two parts: The first part quantifies the nutrient flux caused by the flux of blood plasma $J_{pv}$ from the vessels into the tissue or vice versa. It is determined by Starling's law, which is given by the pressure difference between $p_v$ and $p$ weighted by a parameter $L_p$ for the permeability of the vessel wall,
\begin{equation}
\label{eq:Starling}
J_{pv}( \ov p, p_v ) = L_p ( p_v - \overline{p} ).
\end{equation}
Here, $\overline{p}$ denotes an averaged pressure over the circumference of cylinder cross-sections.. For each parameter $s_i$, we consider a point on the curve $\Lambda_i ( s_i )$. Around this point a circle $\partial B_{R_i}( s_i )$ of radius $R_i$ and perpendicular to $\Lambda_i$ is constructed and the tissue pressure $p$ is averaged with respect to $\partial B_{R_i}( s_i )$,
$$
\overline{p} ( s_i ) = \frac{1}{2 \pi R_i} \int_{\partial B_{R_i}( s_i )} p_{|\Gamma}( x )\,\dd S.
$$
From a physical point of view, the averaging reflects the fact that the 3D-1D coupling is a reduced model, whereas in a fully coupled 3D-3D model, the exchange occurs through the surface. 

In order to account for the permeability of the vessel wall with respect to the nutrients, $J_{pv} \phi_\sigma^v$ is weighted by a factor $1-r_{\sigma}$, where $r_{\sigma}$ is considered as a reflection parameter. The value of $\phi_\sigma^v$ is either set to $\ov \phi_{\sigma}$ or $\phi_v$ depending on the sign of $J_{pv}$, 
$$
\phi_\sigma^v = 
\begin{cases}
\phi_v, & p_v \geq \overline{p}, \\
\ov \phi_{\sigma}, & p_v < \overline{p}.
\end{cases}
$$
The second part of the law \cref{eq:KedemKatalchsky} is a Fickian type law, accounting for the tendency of the nutrients to balance out their concentration levels. Again, the 3D quantity $\phi_{\sigma}$ has to be averaged such that it can be related to the 1D quantity $\phi_v$,
$$
\overline{\phi}_\sigma ( s_i ) = 
\frac{1}{2 \pi R_i} \int_{\partial B_{R_i}( s_i )} \phi_{\sigma}|_\Gamma( x )\,\dd S.
$$
The permeability of the vessel wall is represented by another parameter $L_\sigma$.

Since the exchange processes between the vascular network and the tissue occur at the vessel surface $\Gamma$, we concentrate the flux $J_{\sigma v}$ by means of the Dirac measure $\delta_\Gamma$, i.e., with the distributional space $\D' = (C_c^\infty(\Omega))'$ we define 
$$\langle \delta_\Gamma, \varphi \rangle_{\D' \times \D} = \int_\Gamma \varphi|_{\Gamma}(x) \, \dd S \quad \text{ for all } \varphi \in \D.$$ This yields the following source term in \cref{Eq:DerivationRD},
$$
S_{\sigma v}(\phi_\sigma,p,\phi_v,p_v) = J_{\sigma v} ( \phi_\sigma, p, \Pi_\Gamma \phi_v, \Pi_\Gamma p_v) \delta_\Gamma,
$$
where $\Pi_\Gamma \in \L(L^2(\Lambda);L^2(\Gamma))$ is the projection of the 1D quantities onto the cylindrical surface $\Gamma$ via extending the function value $\Pi_\Gamma \phi_v(s) = \phi_v(s_i)$ for all $s \in \partial B_{R_i}(s_i)$. In particular, we have
$$\int_{\p B_{R_i}(s_i)}\Pi_\Gamma \phi_v(x) \,\dd S = 2 \pi R_i \phi_v(s_i). $$

We assume a volume-averaged velocity $v$ for the prolilferative cells, hypoxic cells, and the nutrients. This assumption of a volume-averaged velocity is reasonable since the cells are tightly packed. Therefore, we assume $v$ to obey the compressible Darcy law 
\begin{equation}
\label{Eq:DerivationDarcy}
\begin{aligned}
v  & = - K(\nabla p -S_p(\phib,\mu_P,\mu_H)) , \\
-\div ( K \nabla p )  & =  J_{pv}( p,\Pi_\Gamma p_v ) \delta_\Gamma -  \div (KS_p(\phib,\mu_P,\mu_H)) , 
\end{aligned}
\end{equation} 
where $K>0$ is the permeability and $J_{pv}( p, \Pi_\Gamma p_v ) \delta_\Gamma$ models the flux between the vascular system and the tissue. Moreover, the source $S_p$ is assumed to represent a form of the elastic Korteweg force, e.g., see \cite{frigeri2018on}, and we correct the chemical potential by the haptotaxis and chemotaxis adhesion terms as done in \cite{garcke2018multiphase}, giving
\begin{equation}
	\label{Eq:Korteweg}
S_p(\phib,\mu_P,\mu_H)=-(\nabla \mu_P+\chi_c \nabla \phi_\sigma + \chi_h \nabla \ecm) \phi_P -(\nabla \mu_H+\chi_c \nabla \phi_\sigma+ \chi_h \nabla \ecm) \phi_H.
\end{equation}

Collecting \cref{Eq:DerivationCH}--\cref{Eq:DerivationDarcy}, we arrive at a model governed by the system,
\begin{equation} \label{Eq:Model3D} \begin{aligned}
		\pt \phi_P+ \div(\phi_P v)            & = \div (m_P( \phib) \nabla \mu_P)+ S_P( \phib),                                                               \\
		\mu_P                                 & =   \p_{\phi_P} \Psi(\phi_P,\phi_H,\phi_N) - \eps^2_P \Delta \phi_P - \chi_c \phi_\sigma-\chi_h \ecm ,                                                           \\
		\pt \phi_H+ \div(\phi_H v)            & = \div (m_H( \phib) \nabla \mu_H)+S_H( \phib)   ,                                                              \\
		\mu_H                                 & =   \p_{\phi_H}\Psi(\phi_P,\phi_H,\phi_N) - \eps^2_H \Delta \phi_H - \chi_c \phi_\sigma-\chi_h \ecm ,                                                           \\
		\pt \phi_N                            & = S_N( \phib)   ,                                                                                                                                   \\
		\pt \phi_\sigma + \div(\phi_\sigma v) & = \div (m_\sigma( \phib) (D_\sigma \nabla \phi_\sigma - \chi_c  \nabla (\phi_P +\phi_H  ))+S_\sigma( \phib) + J_{\sigma v} (\phi_\sigma, p, \Pi_\Gamma \phi_v,  \Pi_\Gamma p_v ) \delta_\Gamma ,\\
		\pt \mde    & = \div (m_{\MDE}( \phib) D_{\MDE} \nabla \mde)+S_{\MDE}( \phib), \\
		\pt \taf  & = \div (m_{\TAF}( \phib) D_{\TAF} \nabla \taf)+S_{\TAF}( \phib), \\
		\pt \ecm                        & = S_{\ECM}(\phib),  \\
		v  & = - K(\nabla p -S_p(\phib,\mu_P,\mu_H)) , \\
        -\div ( K \nabla p ) &=    J_{pv}( p, \Pi_\Gamma p_v ) \delta_\Gamma - \div (KS_p( \phib,\mu_P,\mu_H)),
	\end{aligned} \end{equation}
in the time-space domain $(0,T)\times \Omega$ with source functions $S_P, S_H, S_N, S_\sigma, S_{\MDE}, S_{\TAF}, S_{\ECM}, S_p$, recall \cref{Eq:DerivationSourceCH}, \cref{Eq:DerivationSourceRD} and \cref{Eq:Korteweg}, with properties laid down in \cref{Ass:Assumption} of \cref{Sec:Theorem}. We supplement the system with the following boundary and initial conditions,
\begin{equation} 
\label{Eq:Initial3D} 
\begin{aligned}
	 m_\alpha(\phib)  \p_n \mu_\alpha - \phi_\alpha v \cdot n = m_\beta (\phib)  \p_n \phi_\beta =  \p_n \phi_\gamma    & =0  &  & \text{on } (0,T) \times \p\Omega,  \\
	 p &= p_\infty & & \text{on} (0,T) \times \p \Omega_D, \\
	 \p_n p &= 0 & & \text{on} (0,T) \times \p \Omega \backslash \p \Omega_D, \\
		\phi_\delta(0)  & =\phi_{\delta,0} &  & \text{in } \Omega,                                 
\end{aligned}
\end{equation}
for $\alpha \in \{P,H\}$, $\beta \in \RD$, $\gamma \in \CH \cup \{\ECM\}$, and $\delta \in \A$. Here, $\phi_{\delta,0}$ are given functions with regularity as in \cref{Ass:Assumption} of \cref{Sec:Theorem}, $\partial_n f = \nabla f \cdot n$ denotes the normal derivative of a function $f$ at the boundary $\partial \Omega$ with the outer unit normal $n$ and $\p \Omega_D$ is a part of the boundary with positive measure representing an inlet where the pressure is set to the time-dependent function $p_\infty:(0,T) \times \Omega \to \R$.

\subsection{One-dimensional model for flow and nutrient transport in the vascular network}
\label{sec:1Dmodel}
Since the vascular network typically forms a system of small inclusions, we average all the physical units across the cross-sections of the single blood vessels and set them to a constant with respect to the angular and radial component. This means that the 1D variables $\phi_v$ and $p_v$ on a 1D vessel $\Lambda_i$ depend only on $s_i$. For further details related to the derivation of 1D pipe flow and transport models, we refer to \cite{koppl20203d}. Accordingly, the 1D model equations for flow and transport on $\Lambda_i$ read as follows,
\begin{equation}  \label{Eq:Model1D}
		\begin{aligned}
			\pt \phi_v + \partial_{s_i} (v_v \phi_v) & =  \partial_{s_i} (m_v(\phi_v)D_v \partial_{s_i} \phi_v) -2\pi R_i J_{\sigma v} (\ov \phi_\sigma, \ov p, \phi_v,  p_v),\\
			-  \;\partial_{s_i} ( R_i^2 \pi K_{v,i} \; \partial_{s_i} p_v )   & = -2 \pi R_i J_{pv}( \overline{p}, p_v ).
		\end{aligned} 
\end{equation}
As in \cref{Eq:Model3D}, the fluxes $J_{\sigma v}$ and $J_{pv}$ account for the exchange processes between the blood vessels and the tissue. The permeability is given by the relation
$
K_{v,i} = \frac{R_i^2}{8 \mu_{bl}}, 
$
where $\mu_{bl}$ represents the viscosity of blood. For convenience, we fix it to a constant value, i.e., the non-Newtonian behavior of blood is not considered in this work. The diffusivity parameter $D_v$ is the same as the one of the nutrients in the blood. The blood velocity $v_v$ is calculated as follows via a Darcy-tpye model,
$$
v_v = - R_i^2 \pi K_{v,i} \partial_{s_i} p_v. 
$$

In order to interconnect the different solutions on $\Lambda_i$ at inner networks nodes on intersections
$x \in \partial \Lambda_i \setminus \partial \Lambda,$ we require the continuity of pressure and concentration as well as the conservation of mass to obtain a physically relevant solution. To formulate these coupling conditions in a mathematical way, we define for each bifurcation point $x$ an index set $N(x) \subset \left\{ 1,\ldots,N \right\}$:
$$
N(x ) = \left\{\left. i \; \right| \; x \in \partial \Lambda_i,\;i \in \left\{ 1,\ldots,N \right\} \right\}.
$$
Using this notation, we have for $p_v$ and $\phi_v$ four different coupling conditions at an inner node $x \in \partial \Lambda_i$:
\ \\
	\begin{enumerate} \itemsep.8em
	\item Continuity of $p_v$:
	$$
	p_v \big|_{\Lambda_i}(x) = p_v \big|_{\Lambda_j}(x) \quad \text{ for all } j \in N(x) \setminus \left\{i \right\}.
	$$
	\item Mass conservation with respect to $p_v$:
	$$
	\sum_{j \in N(x)}  -\frac{  R_j^4 \pi}{8 \mu_{\mathrm{bl}} } \frac{\partial p_v}{\partial s_j} \bigg|_{\Lambda_j}(x)  = 0.
	$$	
	\item Continuity of $\phi_v$:
	$$
	 \phi_v \big|_{\Lambda_i}(x) = \phi_v \big|_{\Lambda_j}(x) \quad \text{ for all }~ j \in N(x) \setminus \left\{i \right\}.
	$$
	\item Mass conservation with respect to $\phi_v$:
	$$
	\sum_{j \in N(x)}  \left( v_v \phi_v - m_v(\phi_v)D_v \frac{\partial \phi_v}{\partial s_j} \right) \bigg|_{\Lambda_j}(x)  = 0.
	$$
\end{enumerate}
Further, we decompose the boundary of $\Lambda$ into a Dirichlet boundary $\partial \Lambda_D$ and a Neumann boundary $\partial \Lambda_N$ such that
$
\partial \Lambda = \partial\Lambda_D \,\dot\cup \, \partial\Lambda_N.
$
We introduce the inlet functions $\phi_{v,\infty}, p_{v,\infty}:(0,T) \to \R$ on $\p \Lambda_D$  and prescribe the following boundary data  for $\phi_v$ and $p_v$,
\begin{equation} \begin{aligned}
	 \phi_v - \phi_{v,\infty}=p_v - p_{v,\infty}  &= 0 &&\text{ on } (0,T) \times \partial\Lambda_D, \\ \p_{n_\Lambda} \phi_v =\partial_{n_\Lambda} p_v  &= 0  &&\text{ on } (0,T) \times \partial\Lambda_N.
\end{aligned}
\label{Eq:Initial1D}
\end{equation}

\section{Analytical Preliminaries} \label{Sec:Preliminaries}

Notationally, we equip the function spaces $L^p(\Omega)$, $L^p(\Lambda)$, $W^{m,p}(\Omega)$, $W^{m,p}(\Lambda)$ with the norms $|\cdot|_{L^p(\Omega)}$, $|\cdot|_{L^p(\Lambda)}$, $|\cdot|_{W^{m,p}(\Omega)}$, $|\cdot|_{W^{m,p}(\Lambda)}$. In the case of $d$-dimensional vector functions, we write $L^p(\Omega;\R^d)$ and in the same way for the other Banach spaces, but we do not make this distinction in the notation of norms, scalar products and applications with its dual.

Throughout this paper, $C<\infty$ stands for a generic constant, which may change from line to line. For brevity, we write $x \lesssim y$ for $x \leq Cy$. We recall the Poincar\'e--Wirtinger and Sobolev inequalities, see \cite{brezis2010functional,roubicek,evans2010partial},
\begin{equation}\begin{aligned}
	|f-f_\Omega|_{L^p(\Omega)} & \lesssim |\nabla f|_{L^p(\Omega)} && \text{for all } f \in W^{1,p}(\Omega),   \\
	|f|_{L^p(\Omega)}              & \lesssim  |\nabla f|_{L^p(\Omega)} && \text{for all } f \in W^{1,p}_0(\Omega), \\[-.18cm]
	 |f|_{W^{m,q}(\Omega)} & \lesssim |f|_{W^{k,p}(\Omega)} && \text{for all } f \in W^{k,p}(\Omega), \quad k-\frac{d}{p}  \geq m-\frac{d}{q} , \quad k\geq m, 
\end{aligned} \label{Eq:SobolevInequality} \end{equation}
where $p,q \in [1,\infty)$ and $f_\Omega=\frac{1}{|\Omega|}\int_\Omega f(x) \,\textup{d} x$ denotes the mean of $f$ with respect to $\Omega$. Also, the last inequality yields the continuous embedding $W^{k,p}(\Omega) \hookrightarrow W^{m,q}(\Omega)$.

For a given Banach space $X$, we define the Bochner space, see e.g.,\cite{diestel1977vector},
$$L^p(0,T;X)=\{ u :(0,T) \to X: u \text{ is strongly measurable, } \int_0^T |u(t)|_X^p \,\text{d$t$} < \infty \},$$
where $1 \leq p < \infty$, with the norm $\|u\|_{L^p X}^p = \int_0^T |u(t)|_X^p \,\textup{d} t$. For $p=\infty$, we equip $L^\infty(0,T;X)$ with the norm $\|u\|_{L^\infty X} = \text{ess\,sup}_{t\in (0,T)} |u(t)|_X$. Moreover, we introduce the Sobolev--Bochner space,
$$W^{1,p}(0,T;X)=\{ u \in L^p(0,T;X) : \pt u \in L^p(0,T;X) \}.$$

Let $X$, $Y$, $Z$ be Banach spaces such that $X$ is compactly embedded in $Y$, and $Y$ is continuously embedded in $Z$, i.e.,  $X\com Y \con Z$. In the proof of the existence theorem below, we make use of the Aubin--Lions--Simon compactness lemma, see \cite[Corollary 4]{simon1986compact},
\begin{equation}\begin{alignedat}{2} L^p(0,T;X) \cap W^{1,1}(0,T;Z) &\com L^p(0,T;Y), &&\quad 1\leq p<\infty, \\
		L^\infty(0,T;X) \cap W^{1,r}(0,T;Z) &\com C^0([0,T];Y), &&\quad r >1,
	\end{alignedat}
	\label{Eq:EmbeddingComp}
\end{equation}
where we equip an intersection space $X \cap Y$ with the norm $\|\cdot\|_{X\cap Y}=\max\{\|\cdot\|_X,\|\cdot\|_Y\}$. Further, we make use of the following continuous embeddings, see \cite[Theorem 3.1, Chapter 1]{lions2012non},
\begin{equation}\begin{alignedat}{2} L^2(0,T;Y)  & \cap H^1(0,T;Z)   &  & \con C^0([0,T];[Y,Z]_{1/2}),  \\
	L^\infty(0,T;Y) & \cap C_w([0,T];Z) &  & \con C_w([0,T];Y),
\end{alignedat}
\label{Eq:EmbeddingCont}
\end{equation}
where $[Y,Z]_{1/2}$ denotes the interpolation space between $Y$ and $Z$, see \cite[Definition 2.1, Chapter 1]{lions2012non} for more details. Also, $C_w([0,T];Y)$ denotes the space of the weakly continuous functions on the interval $[0,T]$ with values in $Y$.

We note the following special case of the Gagliardo--Nirenberg inequality, see \cite[Lemma II.2.33]{boyer2012mathematical},
$$|f|_{L^p(\Omega)} \lesssim |f|_{H^1(\Omega)}^\alpha |f|_{L^2(\Omega)}^{1-\alpha} \quad  \text{for all } f \in H^1(\Omega), \quad \frac{1}{p} = \frac12 - \frac{\alpha}{3}, \quad \alpha \in [0,1],$$
which gives in a time-dependent setting, choosing $\alpha=2/q$ with $q \geq 2$,
\begin{equation} \label{Eq:GagliardoTime}
\begin{aligned}\|u\|_{L^q(0,T;L^p(\Omega))}^q = \int_0^T |u(t)|_{L^p(\Omega)}^q \dt &\lesssim \int_0^T |u(t)|_{H^1(\Omega)}^{q \alpha} |u(t)|_{L^2(\Omega)}^{q(1-\alpha)} \dt  \\ &= \int_0^T | u(t)|_{H^1(\Omega)}^2 |u(t)|_{L^2(\Omega)}^{q-2} \dt  \\ &\leq \| u\|_{L^2(0,T; H^1(\Omega))}^2 \|u\|_{L^\infty(0,T;L^2(\Omega))}^{q-2}
\\ &\leq (\max \{ \|u\|_{L^\infty(0,T;L^2(\Omega))}, \| u\|_{L^2(0,T; H^1(\Omega))} \})^q
.\end{aligned}
\end{equation}
In particular, it yields the continuous embedding
$$L^\infty(0,T;L^2(\Omega)) \cap L^2(0,T;H^1(\Omega)) \con L^q(0,T;L^p(\Omega)), \quad \frac{1}{p} + \frac{2}{3q} = \frac12. $$

We also make use of the classical Gr\"onwall--Bellman lemma in the energy estimates to absorb solution-dependent terms on the right hand side of the energy inequalities.
\begin{lemma}[Gr\"onwall--Bellman, cf. {\cite[Lemma II.4.10]{boyer2012mathematical}}] Let $u \in L^\infty(0,T)$, $g \in L^1(0,T;\R_{\geq 0})$ and $u_0 \in \mathbb{R}$. If we have
	$$u(t) \leq u_0 +  \int_0^t g(s) u(s) \, \dd s \quad \text{ for a.e. } t \in (0,T),$$
	then it holds $u(t) \leq u_0 \exp(\int_0^t g(s) \, \dd s)$ for almost every $t \in (0,T)$.
	\label{Lem_Gronwall}
\end{lemma}

\section{Existence of Solutions} \label{Sec:Theorem}

In this section, we lay down some general assumptions on the model that are in force throughout this paper. Under these assumptions, we state the definition of a weak solution, and we then state a theorem, which provides the existence of a weak solution.

For simplicity, we write
$$\begin{gathered}
			S_\alpha = S_\alpha( \phib), \quad
			m_\beta = m_\beta( \phib), \quad
			\Psi = \Psi(\phi_P,\phi_H,\phi_N), \\  J_{p v} = J_{p v}( \ov p, p_v), \quad J_{p v,\Gamma} = J_{p v}(  p,\Pi_\Gamma p_v), \quad J_{\sigma v} = J_{\sigma v}(\ov \phi_\sigma, \ov p, \phi_v,p_v), \quad J_{\sigma v,\Gamma} = J_{\sigma v}(\phi_\sigma, p, \Pi_\Gamma \phi_v,\Pi_\Gamma p_v),
\end{gathered} $$
where $\alpha \in \A$ and $\beta \in \A \backslash \{N,\ECM\}$. \blue{We introduce the scaled parameters $\wdr=2\pi R_i$ and $\wdk=R_i^2 \pi K_{v,i}$ in order to express the 1D model \cref{Eq:Model1D} in a shorter way. Moreover, we define the cut-off operator 
	\begin{equation} \mathcal{C}(x)=\max\{0,\min\{1,x\}\}.\label{Eq:Cutoff}
	\end{equation} }
	Moreover, we introduce the following abbreviations for frequently appearing function spaces, 
$$\begin{alignedat}{4}
		&V = H^1(\Omega) &&\con \quad H = L^2(\Omega) \quad&&\con\quad  V'=(H^1(\Omega))' ,\\
			&V_0 = H^1_D(\Omega) &&\con \quad H = L^2(\Omega) &&\con\quad  V_0'=(H^1_D(\Omega))' ,\\
		&W = W^{1,3/2}(\Omega) ~&&\con \quad H=L^2(\Omega) &&\con\quad  W'=(W^{1,3/2} (\Omega))', \\
		&X= H^1(\Lambda) &&\con \quad Y = L^2(\Lambda) &&\con\quad  X'=(H^1(\Lambda))', \\
		&X_0 = H^1_D(\Lambda) &&\con \quad Y = L^2(\Lambda) &&\con\quad  X_0'=(H^1_D(\Lambda))',
	\end{alignedat}$$
where we have denoted the Sobolev space of vanishing trace on $\p \Omega_D \subset \p \Omega$ by $H_D^1(\Omega) = \{u \in H^1(\Omega) : u|_{\partial \Omega_D} =0 \}$  and in the same way $H_D^1(\Lambda) = \{u \in H^1(\Lambda) : u|_{\partial \Lambda_D} =0 \}$. We equip these spaces of vanishing trace with the norms $|\cdot|_{V_0}=|\nabla \cdot|_H$ and $|\cdot|_{X_0} = |\nablala \cdot|_Y$, respectively. \blue{Here, we use the notation $\nablala$ for the space derivative of the 1D fields.}

The space $W$ with the Lebesgue order $3/2$ becomes useful in the application of the H\"older inequality. Indeed, we have the relation $\frac23 = \frac16 + \frac12$, and therefore, we obtain
	$$|u \varphi|_{L^{3/2}(\Omega)} \leq |u|_{L^6(\Omega)} |\varphi|_H \lesssim |u|_V |\varphi|_H \quad \text{for all } u \in V, \varphi \in H,$$
	where we also applied the Sobolev embedding theorem $V\con L^6(\Omega)$ in the three-dimensional domain $\Omega$. Hence, we have for all $u , \varphi \in V$,
	\begin{equation} \label{Eq:Holder32}
		|u\varphi|_W = \Big(|u\varphi|_{L^{3/2}(\Omega)}^{3/2} + |\nabla(u\varphi)|_{L^{3/2}(\Omega)}^{3/2} \Big)^{2/3} \leq |u\varphi|_{L^{3/2}(\Omega)} + |\nabla(u\varphi)|_{L^{3/2}(\Omega)} \lesssim |u|_V |\varphi|_V, 
	\end{equation}
	where we used the Bernoulli inequality to obtain $(a+b)^r \leq a^r + b^r$ with $a,b \geq 0$, $r \in [0,1]$. 
	

\begin{assumption} \, ~
	\begin{enumerate}[label=(A\arabic*), ref=A\arabic*] \itemsep.25em
		\item $\Omega\subset \R^3$ is a bounded domain with $C^{1,1}$-boundary, \blue{$\Lambda$ is a 1D structure as depicted in \cref{fig:1DNetworkc}, $\Gamma$ is the 2D associated cylindrical surface, see \cref{fig:1DNetworkd},} and $T>0$ denotes a finite time horizon, 
		\item $\phi_{\alpha,0} \in V$ for all $\alpha \in \CH \cup \{\ECM\}$, $\phi_{\beta,0} \in H$ for all $\beta \in \RD$, $\phi_{v,0} \in Y$, $\phi_{v,\infty}, p_{v,\infty} \in H^1(0,T) \subset C([0,T])$ and $p_\infty \in H^1(0,T;H) \cap L^2(0,T;V) \subset C([0;T];H)$, \label{Ass:Initial}
		\item $\chi_c,\chi_h \geq 0$ and $\eps_\alpha, D_\beta, \blue{C_\sigma,} \wdk, \wdr >0$ for $\alpha \in \{P,H\}$, $\beta \in \RD$, \label{Ass:Parameters}
		\item $S_\alpha$ are of the form
		      $$\begin{aligned}
				      S_\alpha(\phib) & =  \sum_{\gamma \in \A} \phi_\gamma f_{\alpha,\gamma} (\phib) ,  &&\alpha \in  \A \backslash \{ N, \ECM\}, \\
				      S_\beta(\phib)  & = f_\beta( \phib), &&\beta \in \{N,\ECM\}, \\
				      \hspace{-5cm}\!\! S_p(\phib,\mu_P,\mu_H) &= \rlap{$\displaystyle \begin{multlined}[t] -\mathcal{C}(\phi_P) (\nabla \mu_P+\chi_c \nabla \phi_\sigma+\chi_h \nabla \ecm) \\ \qquad -\mathcal{C}(\phi_H) (\nabla \mu_H + \chi_c \nabla \phi_\sigma+\chi_h \nabla \ecm), \end{multlined}$}  
			      \end{aligned}$$
			  where $f_{\alpha,\gamma} \in C_b(\R^{|\A|})$, $f_\beta \in \Lip(\R^{|\A|}) \cap PC^1(\R^{| \A|})$, such that $|f_{\alpha,\gamma}|, |f_\beta|, |\p_{\phi_\gamma} f_\beta| \leq f_\infty$ for all $\alpha \in \A\backslash\{N,\ECM\}$, $\beta \in \{N,\ECM\}$, $\gamma \in \A$, \label{Ass:Source}
				\item $J_{pv}$ and $J_{\sigma v}$ are of the form $$\begin{aligned} J_{pv}(y_1,y_2)&=L_p (y_2-y_1), \\ J_{\sigma v}(x_1,y_1,x_2,y_2) &=f_{\sigma,v} (x_1,x_2) J_{pv}(y_1,y_2) + L_\sigma (x_2-x_1),\end{aligned}$$ where $f_{\sigma,v} \in C_b(\mathbb{R}^2)$ such that $|f_{\sigma,v}(x)|\leq f_\infty$ \blue{for all $x \in \R^2$} and $L_p, L_\sigma, K \geq 0$ are sufficiently small in the sense that the prefactors in \cref{Eq:Collect} are positive,
		\label{Ass:Jv}
		\item $m_\alpha \in C_b(\R^{| \A|})$ such that  $0<m_0 \leq m_\alpha(x) \leq m_\infty$ \blue{for all $x \in \R^{|\A|}$} for all $\alpha \in \A \backslash \{N,\ECM\},$ \label{Ass:Mobility}
		\item \label{Ass:Psi} $\Psi \in C^{1}(\R^3)$ non-negative such that $\Psi(0,0,0)=\Psi'(0,0,0)=0$, and there are constants $C_{\Psi_j}$, $j \in \{1,\dots,3\}$, such that \blue{for all $(x,y,z) \in \R^3$ it holds} $$\begin{aligned} \Psi(x,y,z) &\geq C_{\Psi_1} (|x|^2+|y|^2+|z|^2) - C_{\Psi_2}, 
		\\ |\p_x \Psi(x,y,z)|, |\p_y \Psi(x,y,z)|, |\p_z \Psi(x,y,z)| &\leq C_{\Psi_3}(1+|x|+|y|+|z|) .\end{aligned}$$
	\end{enumerate}
	\label{Ass:Assumption}
\end{assumption}

Remarks on the assumptions: 
	\begin{itemize} \itemsep1em 
	\item[\cref{Ass:Source}] After a suitable reformulation of the source functions \cref{Eq:DerivationSourceCH} and \cref{Eq:DerivationSourceRD} with the cut-off operator $\C$, see \cref{Eq:Cutoff}, and replacing the Heaviside functions by the continuous Sigmoid function, the source functions can be brought into the form as stated in assumption \cref{Ass:Source}. Further, the assumption  $f_\beta \in \Lip(\R^{|\A|}) \cap PC^1(\R^{| \A|})$, $\beta \in \{N,\ECM\}$, ensures the validity of the chain rule if $f_\beta$ is composed with a vector-valued Sobolev function; see \cite{murat2003chain,leoni2007necessary}. In particular, we have for all $\alpha \in \A$,
	$$(\nabla f_\beta(\phib), \nabla \phi_\alpha)_H = \sum_{\gamma \in \A} ( \p_{\phi_\gamma} f_\beta(\phib) \nabla \phi_\gamma, \nabla \phi_\alpha)_H \leq f_\infty \sum_{\gamma \in \A} |\nabla \phi_\gamma|_H |\nabla \phi_\alpha|_H .$$
	\item[\cref{Ass:Jv}] We consider the unique, linear and continuous trace operator, see \cite{grosse2013sobolev}, 
	$$\tr_\Gamma: W \to W^{1/3,3/2}(\Gamma)\text{ such that }\tr_\Gamma u = u_{|\Gamma} \text{ for } u \in C^\infty(\Omega), $$ onto the two dimensional \blue{associated cylindrical} surface $\Gamma$ of the one-dimensional network $\Lambda$, \blue{see \cref{fig:1DNetwork}}. In two dimensions, we can apply the Sobolev embedding theorem to obtain $W^{1/3,3/2}(\Gamma) \con L^2(\Gamma)$, see \cref{Eq:SobolevInequality}. Note that this embedding does not hold in three dimensions. Consequently, we have
	\begin{align*}  \| \delta_\Gamma \|_{W'} = \sup_{|\varphi|_W \leq 1}| \langle \delta_\Gamma, \varphi \rangle_W |= \sup_{|\varphi|_W \leq 1} \Big|\int_{\Gamma} \tr_\Gamma\varphi(s) \ds \Big| &\leq \sup_{|\varphi|_W \leq 1}   |\tr_\Gamma\varphi|_{L^1(\Gamma)} \\ &\leq C_{W^{1/3,3/2}(\Gamma)}^{L^1(\Gamma)} |\tr_\Gamma|_{\mathcal{L}(W;W^{1/3,3/2}(\Gamma))},
	\end{align*}
	where $C_{W^{1/3,3/2}(\Gamma)}^{L^1(\Gamma)}$ denotes the embedding constant from $W^{1/3,3/2}(\Gamma) \hookrightarrow L^1(\Gamma)$. Therefore, we have $\delta_\Gamma \in W'$  and  in the following existence proof we often apply the estimate for $\varphi \in W$
	\begin{equation}\begin{aligned} \langle \delta_\Gamma, J_{\alpha v,\Gamma}\varphi \rangle_W = \int_\Gamma J_{\alpha v,\Gamma} \tr_\Gamma\varphi(s) \ds &\leq |J_{\alpha v,\Gamma}|_{L^2(\Gamma)} |\tr_\Gamma \varphi|_{L^2(\Gamma)} \leq C_\Gamma |J_{\alpha v,\Gamma}|_{L^2(\Gamma)} |\varphi|_W,
	\end{aligned} \label{Eq:EstimateDelta}
	\end{equation}
	for $\alpha \in \{\sigma,p\}$, where $$C_\Gamma= C_{W^{1/3,3/2}(\Gamma)}^{L^2(\Gamma)}|\tr_\Gamma|_{\mathcal{L}(W;W^{1/3,3/2}(\Gamma))}.$$
	Further, we can estimate the fluxes by
	\begin{equation} \label{Eq:EstimateJ}
	\begin{aligned}
	|J_{p v,\Gamma}|_{L^2(\Gamma)} &\leq L_p  (C_\Gamma |p|_W+ |\Pi_\Gamma|_{\mathcal{L}(Y;L^2(\Gamma))} |p_v|_Y), \\
	 |J_{\sigma v,\Gamma}|_{L^2(\Gamma)} &\leq f_\infty L_p  (C_\Gamma |p|_W+ |\Pi_\Gamma|_{\mathcal{L}(Y;L^2(\Gamma))} |p_v|_Y) +  L_\sigma (C_\Gamma|\phi_\sigma|_W+|\Pi_\Gamma|_{\mathcal{L}(Y;L^2(\Gamma))}|\phi_v|_Y).
 	\end{aligned}
	\end{equation}
	The assumption of smallness of $L_p$ and $L_\sigma$ is generally accepted in the analysis of very weak solution of the stationary Navier--Stokes equation. There, one also considers a distributional divergence, which should be sufficiently small, see \cite{kim2009existence}. Additionally, in \cite{dangelo2008coupling} the authors have shown well-posedness of an abstract stationary 3D-1D model if the prefactor of the Dirac delta functional is sufficiently small. 
	\item[\cref{Ass:Psi}] The assumption on the potential $\Psi$ is quite typical in the analysis of Cahn--Hilliard equations, see also \cite{fritz2018unsteady,fritz2019local}. 
	In order to take the fourth order polynomial $(x+y+z)^2(1-x-y-z)^2$, we have to extend it by a quadratic function outside of the interval $[0,1]$, i.e., 
	$$\Psi(x,y,z)=\begin{cases} (x+y+z)^2, &x+y+z <0, \\ (x+y+z)^2(1-x-y-z)^2, &x+y+z \in [0,1], \\(1-x-y-z)^2, & x+y+z >1, \end{cases}  $$
	and one can show that $\Psi \in C^2(\R^3;\R)$.
	
	We invoke from \cref{Ass:Psi} and the fundamental lemma of calculus the upper estimate
	\begin{equation}  \label{Eq:PsiGrowth}
	\begin{aligned}
	\Psi(x,y,z) &= \Psi(0,y,z)+ \int_0^x \partial_x \Psi(\tilde x,y,z) \,\dd \tilde x \\ 
	&= \Psi(0,0,0) + \int_0^x \p_x \Psi(\tilde x,y,z) \, \dd \tilde x + \int_0^y \p_y \Psi(0,\tilde y, z) \, \dd \tilde y + \int_0^z \p_z \Psi(0,0,\tilde z) \, \dd \tilde z \\ &\lesssim 1+|x|^2+|y|^2+|z|^2.
	\end{aligned}
	\end{equation}
\end{itemize}


We define a weak solution of the coupled 3D-1D system, see \cref{Eq:Model3D} and \cref{Eq:Model1D}, in the following way.

\begin{definition}[Weak solution] \label{Def:Weak} We call the tuple $( \phib, \mu_P, \mu_H, v, p, \phi_v, v_v, p_v)$ a weak solution of \cref{Eq:Model3D} and \cref{Eq:Model1D} with boundary data \cref{Eq:Initial3D} and \cref{Eq:Initial1D} if the functions $ \phib:(0,T) \times \Omega \to \R^{|\A|}$, $\mu_P, \mu_H, v, p, \phi_v, v_v, p_v: (0,T) \times \Omega \to \R$ have the regularity
	\begin{equation} \begin{aligned}
		  \phi_\alpha  &\in H^1(0,T;V') \cap L^\infty(0,T;V),      &&\alpha \in \{P,H\},              \\
		  \mu_\alpha    & \in L^2(0,T;V),            &&\alpha \in \{P,H\} ,       \\        
		  \phi_\beta  &\in H^1(0,T;H) \cap L^\infty(0,T;V),      &&\beta \in \{N,\ECM\},                               \\
		 \phi_\gamma    & \in H^1(0,T;V') \cap L^\infty(0,T;H) \cap L^2(0,T;V) , && \gamma \in \RD,  \\
		  (v,p-p_\infty) &\in L^2((0,T)\times \Omega;\R^3)\times L^2(0,T;V_0),  \\
		  \phi_v-\phi_{v,\infty}  &   \in H^1(0,T;X_0') \cap L^\infty(0,T;Y) \cap L^2(0,T;X_0) ,   \\
		  (v_v,p_v-p_{v,\infty}) & \in L^2(0,T;Y) \times L^2(0,T;X_0),  
	\end{aligned} \label{Eq:DefSolutionRegularity} \end{equation}
	fulfill the initial data $\phi_\alpha(0)=\phi_{\alpha,0}$, $\alpha \in \A$, $\phi_v(0)=\phi_{v,0}$, and satisfy the following variational form of \cref{Eq:Model3D},
	\begin{equation}  \label{Eq:ModelWeak3D}\begin{aligned}
			\langle \pt \phi_P , \varphi_1 \rangle_{W^{1,3}(\Omega)} -  (\mathcal{C}(\phi_P) v,\nabla \varphi_1)_H +(m_P \nabla \mu_P,\nabla \varphi_1)_H               & = (S_P,\varphi_1)_H,
			\\
			-(\mu_P,\varphi_2)_H + (\p_{\phi_P} \Psi,\varphi_2)_H + \eps^2_P (\nabla \phi_P,\nabla \varphi_2)_H                                         & = \chi_c (\phi_\sigma,\varphi_2)_H+\chi_h (\ecm,\varphi_2)_H,
			\\
			\langle \pt \phi_H , \varphi_3 \rangle_{W^{1,3}(\Omega)} - (\C(\phi_H) v,\nabla \varphi_3)_H +(m_H \nabla \mu_H,\nabla \varphi_3)_H                          & = (S_H,\varphi_3)_H       ,        \\
			-(\mu_H,\varphi_4)_H + (\p_{\phi_H}\Psi,\varphi_4)_H + \eps^2_H (\nabla \phi_H,\nabla \varphi_4)_H                                         & = \chi_c (\phi_\sigma,\varphi_4)_H+\chi_h (\ecm,\varphi_4)_H,
			\\
			(\pt \phi_N , \varphi_5)_H                                                                                                                & = (S_N,\varphi_5)_H,
			\\
			\langle \pt \phi_\sigma,\varphi_6 \rangle_{W^{1,3}(\Omega)}  - ( \C(\phi_\sigma)v,\nabla \varphi_6)_H +  D_\sigma(m_\sigma \nabla \phi_\sigma, \nabla \varphi_6)_H & = (S_\sigma,\varphi_6)_H + \langle \delta_\Gamma, J_{\sigma v,\Gamma} \varphi_6\rangle_W            \\ &\qquad -\chi_c (m_\sigma \nabla (\phi_P +\phi_H  ),\nabla \varphi_6)_H,
			\\
			\langle \pt \mde,\varphi_7 \rangle_V +  D_{\MDE} (m_{\MDE} \nabla \mde,\nabla \varphi_7)_H                                                   & = (S_{\MDE},\varphi_7)_H,
			\\
			\langle \pt \taf,\varphi_8 \rangle_V + D_{\TAF} (m_{\TAF} \nabla \taf,\nabla \varphi_8)_H                                                    & = (S_{\TAF},\varphi_8)_H,
			\\
			\langle \pt \ecm,\varphi_9 \rangle_V                                                                                                             & = (S_{\ECM},\varphi_9)_H,
			\\
			(v,\varphi_{10})_H &= -K (\nabla p, \varphi_{10})_H + K (S_p,\varphi_{10})_H, 
			\\
			K(\nabla  p,\nabla \varphi_{11})_H                                                                                                          & = \langle \delta_\Gamma, J_{pv,\Gamma} \varphi_{11} \rangle_{W}  +K(S_p,\nabla \varphi_{11})_H,
		\end{aligned}
	\end{equation}
	for all $\varphi_j \in V$, $j \in \{1,\dots,9\}$, $\varphi_{10} \in L^2(\Omega;\mathbb{R}^3)$, $\varphi_{11} \in V_0$, and the variational form of \cref{Eq:Model1D},
	\begin{equation}  \begin{aligned}
			\langle \pt \phi_v, \varphi_{12} \rangle_X -  (\C(\phi_v) v_v ,\nablala \varphi_{12} )_Y + D_v (m_v \nablala \phi_v,\nablala \varphi_{12})_Y & = -\wdr(J_{\sigma v},\varphi_{12})_Y  ,   \\
			(v_v, \varphi_{13})_Y &= - \wdk (\nablala p_v,\varphi_{13})_Y, \\
			\wdk( \nablala p_v, \nablala \varphi_{14} )_Y                                                                                                    & =  -\wdr(J_{pv},\varphi_{14})_Y,
		\end{aligned}
		\label{Eq:ModelWeak1D}
	\end{equation}
	for all $\varphi_j \in X_0$, $j\in \{12,14\}$, $\varphi_{13} \in Y$.
\end{definition}

The initial data $\phi_\alpha(0)=\phi_{\alpha,0}$, $\alpha \in \A$, are well-defined with assumption \cref{Ass:Initial} on the regularity of the initial data. Indeed, from the regularity given in \cref{Eq:DefSolutionRegularity}, we achieve, by the embeddings \cref{Eq:EmbeddingCont}, the continuity-in-time regularity
	$$\begin{aligned}
		\phi_\alpha   &\in C^0([0,T];H) \cap C_w([0,T];V),  && \alpha \in \CH \cup \{\ECM\}, \\
		\phi_\beta &\in C^0([0,T];V') \cap C_w([0,T];H),  &&\beta \in \RD, \\
		\phi_v &\in C^0([0,T];X_0') \cap C_w([0,T];Y),  &&
   \end{aligned}$$
and therefore, $\phi_\alpha(0)$ is well-defined in $H$, $\phi_\beta(0)$ in $V'$ and $\phi_v(0)$ in $X_0'$.

We use a mixed boundary approach for $p, \phi_v, p_v$, e.g., for the pressure $p$ we define $\tilde p= p-p_\infty$ with $\tilde p|_{\p \Omega_D}=0$ and $(\p_n \tilde p+ \p_n p_\infty)_{\p \Omega \backslash \p \Omega_D}=0$. Hence, we consider the partial differential equation
$$-\div (K \nabla \tilde p) - \div (K \nabla p_\infty) = \delta_\Gamma J_{pv,\Gamma} - \div S_p,$$
with the weak form with the test function $q \in V_0$
$$K (\nabla \tilde p + \nabla p_\infty, \nabla q)_H - K (\p_n \tilde p + \p_n p_\infty ,q)_{L^2(\p \Omega)} = \langle \delta_\Gamma, J_{pv,\Gamma} q \rangle_W + (K S_p, \nabla q)_H - (K S_p \cdot n,q)_{L^2(\p \Omega)},   $$
or, after the cancellation of the boundary terms,
$$K (\nabla p , \nabla q)_H  = \langle \delta_\Gamma, J_{pv,\Gamma} q \rangle_W + (K S_p, \nabla q)_H.   $$

The main result of this paper involves stating the existence of a weak solution of the 3D-1D model, see \cref{Eq:Model3D} and \cref{Eq:Model1D},  in the sense of \cref{Def:Weak}.

\begin{theorem}[Existence of a global weak solution] Let \cref{Ass:Assumption} hold.
	Then there exists a weak solution tuple $( \phib, \mu_P, \mu_H, p, \phi_v, p_v)$ to the 3D-1D model in the sense of \cref{Def:Weak}, which additionally satisfies the energy inequality 
	\begin{equation} \label{Eq:EnergySolution}
		\begin{aligned}
				& \quad \|\Psi\|_{L^\infty(0,T;L^1(\Omega))} +\sum_{\mathclap{\alpha \in \CH \cup \{\ECM\}}} ~ \| \phi_\alpha\|^2_{L^\infty(0,T;V)} + ~\sum_{\mathclap{\beta \in \{P,H\}}} ~ \| \mu_\beta \|_{L^2(0,T;V)}^2 + \sum_{\gamma \in \RD}  
			   \|\phi_\gamma\|^2_{L^\infty(0,T;H) \cap L^2(0,T;V)}   \\[-.1cm] &\qquad +\|v\|_{L^2(0,T;H)}^2 + \|p\|_{L^2(0,T;V)}^2 +\|\phi_v\|^2_{L^\infty(0,T;Y) \cap L^2(0,T;X)}    +\|v_v\|_{L^2(0,T;Y)}^2+ \|p_v\|_{L^2(0,T;X)}^2  \\ &\quad \lesssim  1+ |\phi_{v,0}|_{Y}^2+  \sum_{\mathclap{\alpha \in \CH \cup \{\ECM\}}} ~ | \phi_{\alpha,0}|^2_{V} + \sum_{\beta \in \RD}  
			   |\phi_{\beta,0}|^2_{H}  + \|p_\infty \|^2_{L^2(0,T;V)}  + |\phi_{v,\infty}|^2_{H^1(0,T)} + |p_{v,\infty}|^2_{L^2(0,T)}  .
		   \end{aligned} \end{equation}
	\label{Thm:Existence}
\end{theorem}

\section{Proof of Theorem 1} \label{Sec:Proof}

To prove the existence of a weak solution, we use the Faedo--Galerkin method \cite{evans2010partial} and semi-discretize the original problem in space. The discretized model can be formulated as an ordinary differential equation system and by the Cauchy--Peano theorem \cite{walter1998ordinary}, we conclude the existence of a discrete solution, see \cref{Sec:Discrete}.  Having derived energy estimates in \cref{Sec:Energy}, we deduce from the Banach--Alaoglu theorem the existence of limit functions which eventually form a weak solution, see \cref{Sec:Limit}. This method is by now standard in the analysis of tumor growth models, e.g.,see \cite{frigeri2015diffuse,jiang2015well,lowengrub2013analysis, frigeri2018on}. Nevertheless, the novel nonlinear coupling of the equations requires a thorough proof of the existence of a solution to the system.

\subsection{Faedo--Galerkin discretization} \label{Sec:Discrete} 
We introduce the discrete spaces
\begin{align*}
	H_k & =\text{span}\{ h_1,\dots,h_k\}, \\
	H_k^0 &= \text{span}\{ h_1^0,\dots,h_k^0\}, \\
	Y_k & =\text{span}\{ y_1,\dots,y_k\},
\end{align*}
where $h_j: \Omega \to \R$, $h_j^0: \Omega \to \R$, $y_j : \Lambda \to \R$, $j \in \{1,\dots,k\}$, are the eigenfunctions to the eigenvalues $\lambda_{h,j},  \lambda_{h^0,j},\lambda_{y,j} \in \R$ of the following respective problems
$$\begin{aligned}
(\nabla h_j,\nabla v)_H &= \lambda_{h,j} (h_j,v)_H &&\forall  v \in V,  \\
(\nabla h_j^0,\nabla v)_H &= \lambda_{h^0,j} (h_j^0,v)_H &&\forall  v \in V_0,  \\
(\nabla y_j,\nabla v)_Y &= \lambda_{y,j} (y_j,v)_Y &&\forall  v \in X_0. 
\end{aligned}$$
Since the inverse Neumann--Laplace operator is a compact, self-adjoint, injective, positive operator on $L_0^2(\Omega)$, we conclude by the spectral theorem, see e.g.,\cite[12.12 and 12.13]{alt2016linear}, that
\begin{alignat*}{3}
	 & \{h_j\}_{j \in \mathbb{N}} &  & \text{ is an orthonormal basis in } H &  & \text{ and orthogonal in } V.
\end{alignat*}
Therefore, $\blue{\cup_{k\in\mathbb{N}}} H_k$ is dense in $V$. Additionally, $\{h_j\}_{j \in \N}$ is a basis in $H^2_N(\Omega)=\{u \in H^2(\Omega) : \partial_n u=0 \text{ on } \p\Omega\}$, see \cite{garcke2016global}.

Next, we investigate the inverse Dirichlet--Neumann Laplacian $(-\Delta)^{-1}|_H : H \to H$, see, e.g., \cite{saedpanah2014well} for the consideration of the Dirichlet--Neumann Laplacian in a Faedo--Galerkin approach. According to the Lax--Milgram theorem, for all $f \in H$ there exists a unique solution $u_f\in V_0$ to the problem 
$$(\nabla u_f,\nabla v)_H = (u_f,f)_H \quad \forall v \in V_0.$$
Additionally, it holds $|u_f|_{V_0} \lesssim |f|_H$ for all $f\in H$ and we can construct an operator $T \in \L(H;V_0)$ with $Tf=u_f$. Since $V_0$ is compactly embedded in $H$, we conclude the compactness of $T \in \L(H;H)$. Taking the test function $v=Tg$ for an arbitrary element $g \in H$, we obtain the self-adjointness of $T$,
$$(Tg,f)_H=(\nabla Tf,\nabla Tg)_H=(g,Tf)_H,$$
and taking $g=f$ yields the positivity of $T$,
$$(Tf,f)_H=|\nabla Tf|_H^2 \geq 0.$$
Additionally, $T$ is injective, since $Tf=0$ yields $(f,v)_H=0$ for all $v \in H$ and hence, $f=0$ almost everywhere. Similarly, we can derive the same results for an operator $\widetilde T \in \L(Y;Y)$ corresponding to the eigenvalue problem on $Y$. Hence, by the spectral theorem we conclude 
\begin{alignat*}{3}
& \{h_j^0\}_{j \in \mathbb{N}} &  & \text{ is an orthonormal basis in } H &  & \text{ and orthogonal in } V_0, \\
& \{y_j\}_{j \in \mathbb{N}} &  & \text{ is an orthonormal basis in } Y &  & \text{ and orthogonal in } X_0.
\end{alignat*}
Additionally, we deduce that $\blue{\cup_{k\in\mathbb{N}}}H_k^0$ is dense in $V_0$ and $\blue{\cup_{k\in\mathbb{N}}}Y_k$ is dense in $X_0$.

We consider the Faedo--Galerkin approximations, $\alpha \in \A$, $\beta \in \{P,H\}$,
\begin{equation}\begin{gathered}
		\phi_\alpha^k (t) = \sum_{j=1}^k \xi_{\alpha,j}(t) h_j,
		\quad \mu_\beta^k (t) = \sum_{j=1}^k \zeta_{\beta,j}(t) h_j,
		\quad \phi_v^k(t) = \phi_{v,\infty}(t)+ \sum_{j=1}^k \xi_{v,j}(t) y_j,  \\
		p^k(t) =p_\infty(t) + \sum_{j=1}^k \zeta_{p,j}(t) h^0_j, \quad p_v^k(t) =  p_{v,\infty}(t) + \sum_{j=1}^k \zeta_{p_v,j}(t) y_j,
	\end{gathered}
	\label{Eq:GalerkinAnsatzFunctions}
\end{equation}
where $(\xi_{\alpha,j})_{\alpha \in \A}:(0,T) \to \R^{|\A|}$,  $(\zeta_{\beta,j})_{\beta \in \{P,H\}}:(0,T) \to \R^{2}$ and 
$\xi_{v,j},\, \zeta_{p,j}, \, \zeta_{p_v,j} : (0,T) \to \R$
 are coefficient functions for all $j \in \{1,\dots,k\}$.
To simplify the notation, we set $\phib^k = (\phi_\alpha^k)_{\alpha \in \A}$, and
$$\begin{gathered}
S_\alpha^k = S_\alpha( \phib^k), \quad
m_\beta^k = m_\beta( \phib^k), \quad
\Psi^k = \Psi(\phi_P^k,\phi_H^k,\phi_N^k), \\  J_{p v}^k = J_{p v}( \ov p^k, p_v^k), \quad J_{p v,\Gamma}^k = J_{p v}(  p^k,\Pi_\Gamma p_v^k), \quad J_{\sigma v}^k = J_{\sigma v}(\ov \phi_\sigma^k, \ov p^k, \phi_v^k,p_v^k), \quad J_{\sigma v,\Gamma}^k = J_{\sigma v}(\phi_\sigma^k, p^k, \Pi_\Gamma \phi_v^k,\Pi_\Gamma p_v^k),
\end{gathered} $$
where $\alpha \in \A$ and $\beta \in \A \backslash \{N,\ECM\}$. The Faedo--Galerkin system of the model then reads
\begin{equation}  \label{Eq:ModelGalerkin3D}\begin{aligned}
		(\pt \phi_P^k , \varphi_1)_H -  (\C(\phi_P^k) v^k,\nabla \varphi_1)_H +(m_P^k \nabla \mu_P^k,\nabla \varphi_1)_H                             & = (S_P^k,\varphi_1)_H,
		\\
		-(\mu_P^k,\varphi_2)_H + (\p_{\phi_P^k}\Psi^k,\varphi_2)_H + \eps^2_P (\nabla \phi_P^k,\nabla \varphi_2)_H                                           & = \chi_c (\phi_\sigma^k,\varphi_2)_H +\chi_h (  \ecm^k,\varphi_2)_H,
		\\
		( \pt \phi_H^k , \varphi_3 )_H -  (\C(\phi_H^k) v^k,\nabla \varphi_3)_H +(m_H^k \nabla \mu_H^k,\nabla \varphi_3)_H                             & = (S_H^k,\varphi_3)_H           ,           \\
		-(\mu_H^k,\varphi_4)_H + (\p_{\phi_H^k}\Psi^k,\varphi_4)_H + \eps^2_H (\nabla \phi_H^k,\nabla \varphi_4)_H                                           & = \chi_c (\phi_\sigma^k,\varphi_4)_H+ \chi_h ( \ecm^k,\varphi_4)_H,
		\\
		(\pt \phi_N^k , \varphi_5)_H                                                                                                                        & = (S_N^k,\varphi_5)_H,
		\\
		(\pt \phi_\sigma^k,\varphi_6)_H - (\C(\phi_\sigma^k) v^k,\nabla \varphi_6)_H+  D_\sigma(m_\sigma^k \nabla \phi_\sigma^k, \nabla \varphi_6)_H & = (S_\sigma^k,\varphi_6)_H + \langle \delta_\Gamma, J_{\sigma v,\Gamma}^k \varphi_6\rangle_W             \\ &\qquad -\chi_c (m_\sigma^k \nabla (\phi_P^k +\phi_H^k  ),\nabla \varphi_6)_H,
		\\
		(\pt \mde^k,\varphi_7 )_H +  D_{\MDE} (m_{\MDE}^k \nabla \mde^k,\nabla \varphi_7)_H                                                       & = (S_{\MDE}^k,\varphi_7)_H,
		\\
		( \pt \taf^k,\varphi_8)_H + D_{\TAF} (m_{\TAF}^k \nabla \taf^k,\nabla \varphi_8)_H                                                        & = (S_{\TAF}^k,\varphi_8)_H,
		\\
		(\pt \ecm^k,\varphi_9 )_H                                                                                                                     & = (S_{\ECM}^k,\varphi_9)_H,
		\\
		K(\nabla p^k,\nabla \varphi_{10})_H                                                                                                                  & =   \langle \delta_\Gamma, J_{pv,\Gamma}^k\varphi_{10} \rangle_W + K(S_p^k,\nabla \varphi_{10})_H,
	\end{aligned}
\end{equation}
for all $\varphi_i \in H_k$, $i \in \{1,\dots,9\}$, $\varphi_{10} \in H_k^0$, and
\begin{equation} \label{Eq:ModelGalerkin1D}  \begin{aligned}
		(\pt \phi_v^k, \varphi_{11})_Y + D_v (m_v^k \nablala \phi_v^k,\nablala \varphi_{11})_Y & = (\C( \phi_v^k) v_v^k ,\nablala \varphi_{11} )_Y   -\wdr(J_{\sigma v}^k,\varphi_{11})_Y,
		\\
		\wdk( \nablala p_v^k, \nablala \varphi_{12} )_Y                                                                                             & = -\wdr (J_{pv}^k,\varphi_{12})_Y ,
	\end{aligned}
\end{equation}
for all $\varphi_j \in Y_k$, $j\in \{11,12\}$, where we define the Faedo--Galerkin ansatz for the velocities $v^k, v_v^k$ by
\begin{equation} \label{Eq:ModelGalerkinVelocity}
	\begin{aligned}
		v^k &= -K ( \nabla p^k - S_p^k), \\
		v_v^k &= - \wdk \nablala p_v^k. 
	\end{aligned}
\end{equation}
We equip the system with the initial data,
\begin{equation} \begin{aligned}
	\phi_\alpha^k(0) &= \Pi_{H_k} \phi_{\alpha,0}, &&\alpha \in \A, \\
	\phi_v^k(0) &=\phi_{v,\infty}(0)  + \Pi_{Y_k} \phi_{v,0}, 
\end{aligned}
\label{Eq:GalerkinInitial}\end{equation}
where  $\Pi_{H_k}: H \to H_k$ and $\Pi_{Y_k} : Y \to Y_k$ are the orthogonal projections onto the finite dimensional spaces, which can be written as 
$$\begin{aligned}
	\Pi_{H_k} h = \sum_{j=1}^k (h,h_j)_H h_j, ~\text{ and }~
	\Pi_{Y_k} y = \sum_{j=1}^k (y,y_j)_Y y_j.
\end{aligned}$$

After inserting the Faedo--Galerkin ansatz functions \cref{Eq:GalerkinAnsatzFunctions} into the system \cref{Eq:ModelGalerkin3D}--\cref{Eq:ModelGalerkin1D}, one can see that the Faedo--Galerkin system is equivalent to a system of nonlinear ordinary differential equations in the unknowns $( (\xi_{\alpha,j})_{\alpha \in \A \cup \{ v\}}, (\zeta_{\beta,j})_{\beta \in \CH \cup \{p,p_v\}})_{1\leq j \leq k}$ with the initial data, 
$$\begin{aligned} \xi_{\alpha,j}(0)  & = (\phi_{\alpha,0},h_j)_H, &&\alpha \in \A, \\
		\xi_{v,j}(0)       & = (\phi_{v,0},y_j)_Y. 
\end{aligned}$$
Due to the continuity of the involved nonlinear functions the existence of solutions to \cref{Eq:ModelGalerkin3D}--\cref{Eq:ModelGalerkin1D} with the initial data \cref{Eq:GalerkinInitial} follows from the standard theory of ordinary differential equations, according to the Cauchy--Peano theorem \cite{walter1998ordinary}. We thus have local-in-time existence of a continuously differentiable solution,  $$\begin{aligned} (\phib^k,\mu_P^k,\mu_H^k,p^k-p_\infty,\phi_v^k-\phi_{v,\infty},p_v^k-p_{v,\infty}) &\in [C^1([0,T_k];H_k)]^{|\A|} \times [C^0([0,T_k];H_k)]^{2} \times C^0([0,T_k];H_k^0) \\ & \qquad  \times C^1([0,T_k];Y_k)  \times C^0([0,T_k];Y_k), \end{aligned}$$ to the Faedo--Galerkin problem \cref{Eq:ModelGalerkin3D}--\cref{Eq:ModelGalerkin1D} on some sufficiently short time interval $[0, T_k]$. Further,  we obtain $\div S_p^k \in H$ \blue{by the representation of $S_p^k$, see \cref{Ass:Source},} and therefore, $v^k \in H$ with $\div v^k = -K (\Delta p^k-\div S_p^k)=J_{pv,\Gamma}^k \delta_\Gamma$. Similarly, $v_v^k \in Y$ with $\div v_v^k = - RJ_{pv}^k$.

\subsection{Energy estimates} \label{Sec:Energy}
Next, we extend the existence interval to $[0,T]$ by deriving $T_k$-independent estimates. In particular, these estimates allow us to deduce that the solution sequences converge to some limit functions as $k \to \infty$. It will turn out that exactly these limit functions will form a weak solution to our 3D-1D model \cref{Eq:Model3D}--\cref{Eq:Model1D} in the sense of \cref{Def:Weak}.

\subsubsection*{Step 1 (Testing)}
We derive energy estimates of the model \cref{Eq:ModelGalerkin3D}--\cref{Eq:ModelGalerkin1D} by choosing suitable test functions in the variational form. For the Cahn--Hilliard type equations, we choose $\varphi_1 = \mu_P^k+\chi_c \phi_\sigma^k+\chi_h \ecm^k$, $\varphi_2 = \partial_t \phi_P^k-\mu_P^k$, $\varphi_3 = \mu_H^k+\chi_c \phi_\sigma^k+\chi_h \ecm^k$, $\varphi_4 = \partial_t \phi_H^k-\mu_H^k$, $\varphi_5 =\Pi_{H_k} \p_{\phi_N^k}\Psi^k - \eps_N^2 \Delta \phi_N^k$, and we arrive at the system of equations, 

\begin{equation} \label{Eq:TestingCH}  \begin{aligned}
		(\pt \phi_P^k,\mu_P^k+\chi_c \phi_\sigma^k+\chi_h \ecm^k)_H  + \Big|\sqrt{m_P^k}\nabla \mu_P^k\Big|_H^2
		                                                                                                                                                      & =  (\C(\phi_P^k) v^k, \nabla \mu_P^k+\chi_c \nabla \phi_\sigma^k+\chi_h \nabla \ecm^k)_H \\[-.15cm] &\qquad -(m_P^k \nabla \mu_P^k,\chi_c \nabla \phi_\sigma^k+\chi_h \nabla \ecm^k)_H \\[-.05cm]  &\qquad + (S_P^k,\mu_P^k + \chi_c\phi_\sigma^k+\chi_h \ecm^k)_H,
		\\
		(\p_{\phi_P^k}\Psi^k,\pt \phi_P^k)_H + \frac{\eps^2_P}{2} \ddt |\nabla \phi_P^k|_H^2 + |\mu_P^k|_H^2
		                                                                                                                                                      & = (\mu_P^k+\chi_c \phi_\sigma^k + \chi_h \ecm^k,\pt \phi_P^k)_H+ (\p_{\phi_P^k}\Psi^k,\mu_P^k)_H \\[-.1cm] &\qquad -  (\chi_c \phi_\sigma^k+\chi_h \ecm^k,\mu_P^k)_H  + \eps_P^2 (\nabla \phi_P^k,\nabla \mu_P^k)_H,
		\\
		(\pt \phi_H^k,\mu_H^k+\chi_c \phi_\sigma^k+\chi_h \ecm^k)_H +  \Big|\sqrt{m_H^k}\nabla \mu_H^k\Big|_H^2 
		                                                                                                                                                      & =(\C(\phi_H^k) v^k,\nabla \mu_H^k+\chi_c \nabla \phi_\sigma^k+\chi \nabla \ecm^k )_H  \\[-.15cm] &\qquad -(m_H^k \nabla \mu_H^k,\chi_c\nabla \phi_\sigma^k+\chi_h \nabla \ecm^k)_H  \\&\qquad + (S_H^k,\mu_H^k + \chi_c\phi_\sigma^k+\chi_h \ecm^k)_H,
		\\
		(\p_{\phi_H^k}\Psi^k,\pt \phi_H^k)_H + \frac{\eps^2_H}{2} \ddt |\nabla \phi_H^k|_H^2   + |\mu_H^k|_H^2
		                                                                                                                                                      & = (\mu_H^k+\chi_c \phi_\sigma^k+\chi_h \ecm^k,\pt \phi_H^k)_H + (\p_{\phi_H^k}\Psi^k,\mu_H^k)_H \\[-.1cm] &\qquad -(\chi_c\phi_\sigma^k+\chi_h \ecm^k,\mu_H^k)_H + \eps_H^2 (\nabla \phi_H^k,\nabla \mu_H^k)_H,
		\\
		(\pt \phi_N^k,\Pi_{H_k}\p_{\phi_N^k}\Psi^k)_H  + \frac{\eps_N^2}{2} \ddt |\nabla \phi_N^k|_H^2  & =    (S_N^k,\Pi_{H_k}\p_{\phi_N^k}\Psi^k)_H  + \eps_N^2 (\nabla S_N^k, \nabla \phi_N^k)_H.
	\end{aligned} \end{equation}
We exploit that the time derivative operator is invariant under the adjoint of the orthogonal projection. 

Further, for the reaction-diffusion type equations, we choose $\varphi_6 = C_\sigma \phi_\sigma^k$, $C_\sigma>0$ to be determined, $\varphi_7 = \mde^k$, $\varphi_8 = \taf^k$, $\varphi_9 = \ecm^k -\Delta \ecm^k$, which yields the system,
\begin{equation} \label{Eq:TestingRD}  \begin{aligned}
		\frac{C_\sigma}{2} \ddt |\phi_\sigma^k|_H^2 + C_\sigma D_\sigma \Big|\sqrt{m_\sigma^k} \nabla \phi_\sigma^k \Big|_H^2
		 & =  \chi_c C_\sigma  (m_\sigma^k  \nabla (\phi_P^k +\phi_H^k  ),\nabla \phi_\sigma^k )_H+C_\sigma (S_\sigma^k,\phi_\sigma^k)_H  \\[-.1cm] &\qquad +C_\sigma  (\C(\phi_\sigma^k)  v^k , \nabla \phi_\sigma^k)_H  + C_\sigma\langle \delta_\Gamma,J_{\sigma v,\Gamma}^k\phi_\sigma^k \rangle _W,
		\\
		 \frac12 \ddt |\mde^k|_H^2 + D_{\MDE} \Big|\sqrt{m_{\MDE}^k} \nabla \mde^k\Big|_H^2 
		 & = (S_{\MDE}^k,\mde^k)_H,
		\\
		\frac12 \ddt |\taf^k|_H^2  + D_{\TAF} \Big|\sqrt{m_{\TAF}^k} \nabla \taf^k\Big|_H^2 
		 & = (S_{\TAF}^k,\taf^k )_H,
		\\
		\frac12 \ddt |\ecm^k|_H^2 +\frac12 \ddt |\nabla \ecm^k|_H^2 
		 & =(S_{\ECM}^k,\ecm^k)_H+ (\nabla S_{\ECM}^k,\nabla \ecm^k)_H,
	\end{aligned} \end{equation}
and for the equations in $Y_k$, we choose $\varphi_{12}=C_v(\phi_v^k-\phi_{v,\infty})$, $C_v>0$ to be determined, giving
\begin{equation} \label{Eq:Testing1D}  \begin{aligned}
		\frac{C_v}{2} \ddt |\phi_v^k-\phi_{v,\infty}|_Y^2   +  C_vD_v \Big|\sqrt{m_v^k} \nablala \phi_v^k\Big|^2_Y 
		                                  & =  C_v(\C(\phi_v^k)  v_v^k,\nablala \phi_v^k )_Y
		                                  \! -\! C_v(RJ_{\sigma v}^k+\phi_{v,\infty}',\phi_v^k-\phi_{v,\infty})_Y, 
	\end{aligned} \end{equation}
		Similarly, we test the equations \cref{Eq:ModelGalerkinVelocity}$_2$ by $\frac{1}{\wdk}v_v^k$ to obtain
				\begin{align}
		\frac{1}{\wdk} |v_v^k|_H^2 &= - (\nablala (p_v^k-p_{v,\infty}),v_v^k)_Y  =   -\wdr(J_{pv}^k, p_v^k-p_{v,\infty})_Y . 
		\label{Eq:TestingVel1D}  
		\end{align}
		We test \cref{Eq:ModelGalerkinVelocity}$_1$ by $\frac{1}{K} v^k$ and simplify the first term on the right hand side by comparing it with the equation \cref{Eq:ModelGalerkinVelocity}$_1$ for the velocity $v^k$ and the pressure equation \cref{Eq:ModelGalerkin3D}, which is tested by $\varphi_{10}=p^k-p_\infty$. This procedure yields
		\begin{equation} \begin{aligned}
		\frac{1}{K} |v^k|_H^2  &= - (\nabla (p^k-p_\infty),v^k)_H + (\nabla p_\infty,v^k)_H + (S_p^k,v^k)_H
		\\  &= K (\nabla p^k-S_p^k,\nabla (p^k-p_\infty))_H  + (\nabla p_\infty,v^k)_H + (S_p^k,v^k)_H
		 \\ & =    \langle \delta_\Gamma, J_{pv,\Gamma}^k (p^k-p_\infty) \rangle_W - (\nabla p_\infty,v^k)_H+(S_p^k,v^k)_H . 
		\end{aligned} \label{Eq:TestingVel3D} \end{equation}

	\subsubsection*{Step 2 (Estimates)} We separate this step into three sub-steps by deriving the energy estimates separately for \cref{Eq:TestingCH}, \cref{Eq:TestingRD} and \cref{Eq:Testing1D}.

	\subsubsection*{Step 2.1 (Estimate for \cref{Eq:TestingCH})}

	Adding the equations in \cref{Eq:TestingCH} and \cref{Eq:TestingVel3D} gives
	\begin{equation} \label{Eq:AddingCH} \begin{aligned}
			 &     \ddt  |\Psi^k|_{L^1(\Omega)}   +  \frac{\eps^2_P}{2} \ddt |\nabla \phi_P^k|_H^2
			+\frac{\eps^2_H}{2} \ddt |\nabla \phi_H^k|_H^2 +  \frac{\eps_N^2}{2} \ddt |\nabla \phi_N^k|_H^2 +\Big|\sqrt{m_P^k}\nabla \mu_P^k\Big|_H^2 +   \Big|\sqrt{m_H^k}\nabla \mu_H^k\Big|_H^2  \\
			& \qquad +|\mu_P^k|_H^2 + |\mu_H^k|_H^2  + \frac{1}{K} |v^k|_H^2 \\
			 & = (\C(\phi_P^k)v^k,\nabla \mu_P^k + \chi_c \nabla \phi_\sigma^k + \chi_h \nabla \ecm^k)_H -(m_P^k \nabla \mu_P^k,\chi_c\nabla \phi_\sigma^k+\chi_h \nabla \ecm^k)_H 
			 \\ &\qquad  + (S_P^k,\mu_P^k + \chi_c\phi_\sigma^k+\chi_h \ecm^k)_H +   (\p_{\phi_P^k}\Psi^k-\chi_c \phi_\sigma^k-\chi_h \ecm^k,\mu_P^k)_H + \eps_P^2 (\nabla \phi_P^k,\nabla \mu_P^k)_H 
			 \\ &\qquad+ (\C(\phi_H^k)v^k,\nabla \mu_H^k + \chi_c \nabla \phi_\sigma^k + \chi_h \nabla \ecm^k)_H-(m_H^k \nabla \mu_H^k,\chi_c\nabla \phi_\sigma^k + \chi_h \nabla \ecm^k)_H  
			\\ &\qquad+(S_H^k,\mu_H^k + \chi_c\phi_\sigma^k+ \chi_h \ecm^k)_H + (\p_{\phi_H^k}\Psi^k-\chi_c \phi_\sigma^k - \chi_h \ecm^k,\mu_H^k)_H    + \eps_H^2 (\nabla \phi_H^k,\nabla \mu_H^k)_H
	   \\ &\qquad 	+ (S_N^k,\Pi_{H_k}\p_{\phi_N^k}\Psi^k)_H + \eps_N^2 (\nabla S_N^k,\nabla \phi_N^k )_H  
	   +\langle \delta_\Gamma, J_{pv,\Gamma}^k (p^k-p_\infty) \rangle_W - (\nabla p_\infty,v^k)_H+(S_p^k,v^k)_H \\
			& = \RHS_{\CH}.
		\end{aligned} \end{equation}
		We note that the two convection terms cancel together with the the last term $(S_p^k,v^k)_H$ on the right hand side.
		We apply the H\"older inequality on the terms on the right hand side, and use the assumptions \cref{Ass:Source} and \cref{Ass:Mobility}, which gives
	\begin{equation} \label{Eq:HolderCH} \begin{aligned}
		\RHS_\CH & \leq m_\infty  | \nabla \mu_P^k|_H (\chi_c |\nabla \phi_\sigma^k|_H + \chi_h |\nabla \ecm^k|_H)  + |S_P^k|_H ( |\mu_P^k|_H +  \chi_c |\phi_\sigma^k|_H+  \chi_h |\ecm^k|_H) 
		\\ &\qquad  +|\mu_P^k|_H (|\p_{\phi_P^k}\Psi^k|_H+ \chi_c |\phi_\sigma^k|_H + \chi_H |\ecm^k|_H)  + \eps_P^2 |\nabla \phi_P^k|_H |\nabla \mu_P^k|_H 
	    \\ &\qquad + m_\infty  |\nabla \mu_H^k|_H (\chi_c |\nabla \phi_\sigma^k|_H+ \chi_h  |\nabla \ecm^k|_H) + |S_H^k|_H ( |\mu_H^k|_H + \chi_c |\phi_\sigma^k|_H+ \chi_h |\ecm^k|_H  )
		  \\ &\qquad +  |\mu_H^k|_H (|\p_{\phi_H^k}\Psi^k|_H+\chi_c |\phi_\sigma^k|_H + \chi_h |\ecm^k|_H) + \eps_H^2 |\nabla \phi_H^k|_H |\nabla \mu_H^k|_H 
		  \\ &\qquad + |S_N^k|_{H} |\Pi_{H_k}\p_{\phi_N^k}\Psi^k|_H + \eps_N^2 |\nabla S_N^k|_H |\nabla \phi_N^k |_H + \langle \delta_\Gamma, J_{pv,\Gamma}^k (p^k-p_\infty) \rangle_W + |\nabla p_\infty|_H |v^k|_H.
	\end{aligned} \end{equation}
	We note that the norm of the orthogonal projection is bounded by $1$. We use a similar argument as in \cref{Eq:EstimateDelta} and \cref{Eq:EstimateJ} to estimate the term involving the Dirac delta functional $\delta_\Gamma$, i.e., with the assumption on the form of $J_{pv,\Gamma}^k$, see \cref{Ass:Jv}, we obtain
	$$\begin{aligned} &\langle \delta_\Gamma, J_{pv,\Gamma}^k (p^k-p_\infty) \rangle_W\\ &\leq  C_\Gamma |p^k-p_\infty|_W|J_{pv,\Gamma}^k|_{L^2(\Gamma)}  \\ &\leq C_\Gamma L_p  |p^k-p_\infty|_W (C_\Gamma |p^k-p_\infty|_W + C_\Gamma |p_\infty|_W + |\Pi_\Gamma|_{\mathcal{L}(Y;L^2(\Gamma))} |p_v^k-p_{v,\infty}|_Y +|\Pi_\Gamma|_{\mathcal{L}(Y;L^2(\Gamma))} |p_{v,\infty}|_Y) \\
	&\leq C_1 L_p ( |\nabla p^k|_H^2 + |\nablala p_v^k|_Y^2 + |p_\infty|^2_V + |p_{v,\infty}|^2),
	\end{aligned} $$
	where we also applied the Poincar\'e inequality on $p^k-p_\infty \in V_0$ and $p_v^k-p_{v,\infty} \in X_0$ with the Poincar\'e constants $C_{P,\Omega}$ and $C_{P,\Lambda}$, giving the constant
	$$C_1 =  \max\{2C_\Gamma^2 (C_{W}^{V})^2 (C_{P,\Omega}^2+1);|\Pi_\Gamma|_{\mathcal{L}(Y;L^2(\Gamma))}^2 C_{P,\Lambda}^2 ;|\Pi_\Gamma|_{\mathcal{L}(Y;L^2(\Gamma))}^2 |\Lambda|  \}.$$
	 Further, using the form on $v^k$ and $v_v^k$ gives
		$$\begin{aligned} \langle \delta_\Gamma, J_{pv,\Gamma}^k (p^k-p_\infty) \rangle_W \leq C_1 L_p (K^{-2} |v^k|_H^2 +|S_p^k|_H^2  + \wdk^{-2} |v_v^k|_Y^2 + |p_\infty|^2_V + |p_{v,\infty}|^2).
	\end{aligned} $$
	
	We apply Young's inequality on the norm products to separate the terms. Here, the goal is to make the terms involving $|\mu_P^k|_V$, $|\mu_H^k|_V$, $|\nabla \phi_\sigma^k|_H$ small, since we cannot absorb them with the Gr\"onwall--Bellman lemma later on. We only track the important constants, which are used to absorb the terms on the right hand side with the left hand side, the other ones we simply denote by the generic constant $C$. We have 
	\begin{equation} \label{Eq:YoungCH} \begin{aligned}
		\RHS_\CH & \leq \frac{m_0}{4} | \nabla \mu_P^k|_H^2 + \frac{2m_\infty^2 \chi_c^2}{m_0} |\nabla \phi_\sigma^k|_H^2 + \frac{2m_\infty^2 \chi_h^2}{m_0} |\nabla \ecm^k|_H^2  + 3 |S_P^k|_H^2  +\frac14 \big( |\mu_P^k|_H^2 + \chi_c^2|\phi_\sigma^k|_H^2
		\\ &\qquad + \chi_h^2|\ecm^k|_H^2\big) 
		  + \frac14 |\mu_P^k|_H^2 + 3 ( |\p_{\phi_P^k}\Psi^k|_H^2+ \chi_c^2 |\phi_\sigma^k|_H^2 +\chi_h^2 |\ecm^k|_H^2) + \frac{\eps_P^4}{m_0} |\nabla \phi_P^k|_H^2
		    \\ &\qquad		 + \frac{m_0}{4} |\nabla \mu_P^k|_H^2  + \frac{m_0}{4} |\nabla \mu_H^k|_H^2    + \frac{2m_\infty^2 \chi_c^2}{m_0} |\nabla \phi_\sigma^k|_H^2 + \frac{2m_\infty^2 \chi_h^2}{m_0} |\nabla \ecm^k|_H^2   + 3|S_H^k|_H^2 
		 \\ &\qquad+\frac14 ( |\mu_H^k|_H^2 + \chi_c^2|\phi_\sigma^k|_H^2 + \chi_h^2|\ecm^k|_H^2) + \frac14 |\mu_H^k|_H^2    + 3 ( |\p_{\phi_H^k}\Psi^k|_H^2+ \chi_c^2 |\phi_\sigma^k|_H^2 +\chi_h^2 |\ecm^k|_H^2) 
		  \\ &\qquad + \frac{\eps_H^4}{m_0} |\nabla \phi_H^k|_H^2 + \frac{m_0}{4} |\nabla \mu_H^k|_H^2+ \frac12 |S_N^k|_H^2 + \frac12 |\p_{\phi_N^k}\Psi^k|_H^2 + \eps|\nabla S_N^k|_H^2 + \frac{\eps_N^4}{4\eps} |\nabla \phi_N^k |_H^2 
		  \\ &\qquad + C_1 L_p\big(K^{-2} |v^k|_H^2 +|S_p^k|_H^2 + \wdk^{-2} |v_v^k|_Y^2 + |p_\infty|^2_V + |p_{v,\infty}|^2\big) + \frac{K}{2} |\nabla p_\infty|_H^2 + \frac{1}{2K} |v^k|_H^2,
	\end{aligned} \end{equation}
	where $\eps>0$ is a constant, which will be determined later on, see \cref{Eq:Epsilon} below for more details. We estimate the terms involving the potential $\Psi$ via assumption \cref{Ass:Psi} and afterwards, we collect the terms with the same norms, which yields
	\begin{equation} \label{Eq:CollectCH} \begin{aligned}
		\RHS_\CH & \leq \frac{m_0}{2} \big( | \nabla \mu_P^k|_H^2 +  | \nabla \mu_H^k|_H^2 \big) +\frac12 \big( |\mu_P^k|_H^2 + |\mu_H^k|_H^2 \big) + \frac{4m_\infty^2 \chi_c^2}{m_0}  |\nabla \phi_\sigma^k|_H^2    + \eps|\nabla S_N^k|_H^2 
				\\ &\qquad+ \frac{1}{2K} |v^k|_H^2 + C_1 L_p \big(K^{-2} |v^k|^2_H  + |S_p^k|_H^2 + \wdk^{-2} | v_v^k|_Y^2+ |p_\infty|_V^2 + |p_{v,\infty}|^2 \big) 
				\\ &\qquad+ C \big( |\phi_\sigma^k|_H^2 +  |\phi_P|_V^2 + |\phi_H|_V^2 + |\phi_N|_V^2 +  |\ecm^k|_V^2   +  |p_\infty|_V^2+ |S_P^k|_H^2 + |S_H^k|_H^2 +  |S_N^k|_H^2 \big) .
	\end{aligned} \end{equation}
	We insert this estimate into \cref{Eq:AddingCH}, and absorb the terms involving the chemical potentials, and arrive at the upper bound
	\begin{equation} \label{Eq:FinalCH} \begin{aligned}
		&    \ddt \bigg[  |\Psi^k|_{L^1(\Omega)}   +  \frac{\eps^2_P}{2} |\nabla \phi_P^k|_H^2
	   +\frac{\eps^2_H}{2}  |\nabla \phi_H^k|_H^2 +  \frac{\eps_N^2}{2}  |\nabla \phi_N^k|_H^2 \bigg] +\frac{m_0}{2} |\nabla \mu_P^k|_H^2 +   \frac{m_0}{2} |\nabla \mu_H^k|_H^2 +\frac12 |\mu_P^k|_H^2   \\
	   &\qquad + \frac12 |\mu_H^k|_H^2+\left(  \frac{1}{2K}- \frac{C_1 L_p}{K^2} \right) |v_k|_H^2 \\
	   &\leq  \frac{4m_\infty^2 \chi_c^2}{m_0} |\nabla \phi_\sigma^k|_H 
	   +  \eps|\nabla S_N^k|_H^2 + C_1 L_p \big(|S_p^k|_H^2 + \wdk^{-2} |v_v^k|_Y^2  \big)
	    + C \big( |\phi_\sigma^k|_H^2 +  |\phi_P^k|_V^2 + |\phi_H^k|_V^2  \\ &\qquad+ |\phi_N^k|_V^2+  |\ecm^k|_V^2 + |\phi_v^k|_Y^2  +  |p_\infty|_V^2 + |p_{v,\infty}|^2 + |S_P^k|_H^2 + |S_H^k|_H^2 +  |S_N^k|_H^2 \big).
   \end{aligned} \end{equation}

   \subsubsection*{Step 2.2 (Estimate for \cref{Eq:TestingRD})}

   Adding the equations in \cref{Eq:TestingRD} gives
	\begin{equation} \label{Eq:AddingRD} \begin{aligned}
			 & \frac{C_\sigma }{2} \ddt |\phi_\sigma^k|_H^2
			+ C_\sigma D_\sigma \Big|\sqrt{m_\sigma^k} \nabla \phi_\sigma^k\Big|_H^2
			+  \frac12 \ddt |\mde^k|_H^2 + D_{\MDE} \Big|\sqrt{m_{\MDE}^k} \nabla \mde^k\Big|_H^2+  \frac12 \ddt |\taf^k|_H^2
			\\ &\qquad  + D_{\TAF} \Big|\sqrt{m_{\TAF}^k} \nabla \taf^k\Big|_H^2  +\frac12 \ddt |\ecm^k|_H^2 +\frac12 \ddt |\nabla \ecm^k|_H^2 
			\\
			 & =  \chi_c C_\sigma  (m_\sigma^k  \nabla (\phi_P^k +\phi_H^k ),\nabla \phi_\sigma^k)_H
			+C_\sigma (S_\sigma^k,\phi_\sigma^k)_H + C_\sigma (\C(\phi_\sigma^k) v^k,\nabla \phi_\sigma^k)_H
			+C_\sigma\langle \delta_\Gamma, J_{\sigma v,\Gamma}^k \phi_\sigma^k \rangle_W
			\\ &\qquad + (S_{\MDE}^k,\mde^k)_H + (S_{\TAF}^k,\taf^k)_H +(S_{\ECM}^k,\ecm^k)_H + (\nabla S_{\ECM}^k,\nabla \ecm^k)_H 
			\\ & = \RHS_\RD.
		\end{aligned} \end{equation}
		We estimate the term involving the Dirac delta functional as before, i.e., we use assumption \cref{Ass:Jv} and the inequalities \cref{Eq:EstimateDelta} and \cref{Eq:EstimateJ} to obtain 
		$$\begin{aligned}&C_\sigma\langle \delta_\Gamma, J_{\sigma v,\Gamma}^k \phi_\sigma^k \rangle_W \\ &\leq C_\sigma C_\Gamma |\phi_\sigma^k|_W \Big(f_\infty L_p (C_\Gamma |p^k|_W + |\Pi_\Gamma|_{\mathcal{L}(Y;L^2(\Gamma))} |p_v^k|_Y) + L_\sigma (C_\Gamma |\phi_\sigma^k|_W + |\Pi_\Gamma|_{\mathcal{L}(Y;L^2(\Gamma))} |\phi_v^k|_Y  ) \Big) \\
		&\leq C_2 C_\sigma \max\{L_p;L_\sigma\} \big( |\phi_\sigma^k|_V^2 +|\phi_v^k|_Y^2 +  K^{-2} |v^k|_H^2 + |S_p^k|_H^2 + \wdk^{-2} |v_v^k|_Y^2 + |p_\infty|_V^2 + |p_{v,\infty}|^2   \big),
		\end{aligned}$$
		where
		$$C_2 = \max\{2^4 C_\Gamma^2 (C_W^V)^2; f_\infty^2 C_\Gamma^2 (C_W^{V})^2 (C_{P,\Omega}^2+1); f_\infty^2 |\Pi_\Gamma|_{\mathcal{L}(Y;L^2(\Gamma))}^2 C_{P,\Lambda}^2 ; f_\infty^2 |\Pi_\Gamma|_{\mathcal{L}(Y;L^2(\Gamma))}^2 |\Lambda|; |\Pi_\Gamma|_{\mathcal{L}(Y;L^2(\Gamma))}^2 \}. $$
		In a similar way to the estimates before, we apply H\"older's and Young's inequalities on the terms on the right hand side, which results in
\begin{equation} \label{Eq:HolderRD} \begin{aligned}
	\RHS_\RD & \leq  \frac{C_\sigma m_\infty^2 \chi_c^2}{D_\sigma m_0} |\nabla (\phi_P^k+\phi_H^k)|_H^2 + \frac{C_\sigma D_\sigma m_0}{4} |\nabla \phi_\sigma^k|_H^2 + C( |S_\sigma^k|_H^2 +  |\phi_\sigma^k|_H^2 ) 
	\\ &\qquad  + \frac{C_\sigma D_\sigma m_0}{4} |\nabla \phi_\sigma^k|_H^2+ \frac{C_\sigma}{D_\sigma m_0} |v^k|_H^2 +C_2 C_\sigma  \max\{L_p;L_\sigma\} \big(|\phi_\sigma^k|_V^2 +|\phi_v^k|_Y^2 +  K^{-2} |v^k|_H^2 \\ &\qquad + |S_p^k|_H^2 + \wdk^{-2} |v_v^k|_Y^2 + |p_\infty|_V^2 + |p_{v,\infty}|^2  \big)     +  C \big( |S_{\MDE}^k|_H^2 +  |\mde^k|_H^2  +  |S_{\TAF}^k|_H^2 
	\\ &\qquad+  |\taf^k|_H^2   +  |S_{\ECM}^k|_H^2 + |\ecm^k|_H^2 \big)+ \eps |\nabla S_{\ECM}^k|_H^2 + \frac{1}{4 \eps} |\nabla \ecm^k|_H^2 ,
\end{aligned} \end{equation}
where we used the same constant $\eps$ as before in \cref{Eq:AddingCH} and applied the assumption on the form of $J_{\sigma v,\Gamma}^k$, see \cref{Ass:Jv}. Again, collecting the terms on the right hand side and absorbing the terms with their counterparts, we have
\begin{equation} \label{Eq:FinalRD} \begin{aligned}
	& \frac12 \ddt \bigg[ C_\sigma |\phi_\sigma^k|_H^2  +|\mde^k|_H^2+|\taf^k|_H^2 + |\ecm^k|_H^2 + |\nabla \ecm^k|_H^2 \bigg]
   + D_{\MDE} m_0| \nabla \mde^k|_H^2  
   \\ &\qquad+ D_{\TAF} m_0 |\nabla \taf^k|_H^2+  \frac{C_\sigma}{2} \big(  D_\sigma m_0 - 2C_2 \max\{L_p;L_\sigma\}\big) |\nabla \phi_\sigma^k|_H^2
    \\	& \leq  \left(\frac{C_\sigma}{D_\sigma m_0} + \frac{C_2 C_\sigma \max\{L_p;L_\sigma\} }{K^2}  \right) |v^k|_H^2 + C_2 C_\sigma \max\{L_p;L_\sigma\} \big(|S_p^k|_H^2 + \wdk^{-2} |v_v^k|_Y^2  \big)+ \eps |\nabla S_{\ECM}^k|_H^2  
	\\ &\qquad   + C \big( |\phi_P|_V^2+ |\phi_H|_V^2  + |\phi_\sigma^k|_H^2  + |\mde^k|_H^2 + |\taf^k|_H^2 + |\ecm^k|_V^2 + |\phi_v^k|_Y^2 + |p_\infty|_V^2
	\\ &\qquad  + |p_{v,\infty}|^2+|S_\sigma^k|_H^2 + |S_{\MDE}^k|_H^2
	 + |S_{\TAF}^k|_H^2 + |S_{\ECM}^k|_H^2 \big).
\end{aligned} \end{equation}

\subsubsection*{Step 2.3 (Estimate for \cref{Eq:Testing1D})}	
		Lastly, adding the equations in \cref{Eq:Testing1D} and \cref{Eq:TestingVel1D} gives
		\begin{equation} \label{Eq:Adding1D} \begin{aligned}
			 &\frac{C_v}{2} \ddt |\phi_v^k-\phi_{v,\infty}|_Y^2
		   + C_vD_v \Big|\sqrt{m_v^k} \nablala \phi_v^k\Big|_Y^2 + \frac{1}{\wdk} |v_v^k|_Y^2 \\
			&=C_v (\C(\phi_v^k) v_v^k,\nablala \phi_v^k)_Y -C_v(RJ_{\sigma v}^k+\phi_{v,\infty}',\phi_v^k-\phi_{v,\infty})_Y  -  \wdr (J_{pv}^k,p_v^k-p_{v,\infty})_Y. \end{aligned}
	  \end{equation}
We estimate the last term on the right hand side with the Poincar\'e--Wirtinger inequality \cref{Eq:SobolevInequality} with constant $C_P$, the Darcy law \cref{Eq:ModelGalerkinVelocity} and the Young inequality as follows
$$-\wdr (J_{pv}^k,p_v^k-p_{v,\infty})_Y \leq \wdr C_P |J_{pv}^k|_Y |\nabla p_v^k|_Y = \frac{R^2 C_P^2}{\wdk} |J_{pv}^k|_Y |v_v^k|_Y \leq \frac{R^2 C_P^2}{\wdk} |J_{pv}^k|_Y^2 + \frac{1}{4\wdk} |v_v^k|_Y^2  .  $$
Additionally, repeating the steps from before and using the assumption on the forms on $J_{pv}^k$ and $J_{\sigma v}^k$, see \cref{Ass:Jv}, we arrive at
\begin{equation*} \begin{aligned} &\frac{C_v}{2} \ddt |\phi_v^k-\phi_{v,\infty}|_Y^2
   + C_vD_v m_0 |\nablala \phi_v^k|_Y^2 + \frac{1}{\wdk} |v_v^k|_Y^2
     \\ &\leq \frac{C_v}{D_v m_0} |v_v^k|_Y^2 + \frac{C_vD_vm_0}{4} |\nablala \phi_v^k|_Y^2 +  R^2|J_{\sigma v}^k|_Y^2 + C\big( |\phi_{v,\infty}'|^2 +  |\phi_v^k-\phi_{v,\infty}|^2_Y\big)  + \frac{R^2C_P^2}{\wdk} |J_{pv}^k|_Y^2   + \frac{1}{4\wdk} |v_v^k|_Y^2
          \\ &\leq \frac{C_v}{D_v m_0} |v_v^k|_Y^2 + \frac{C_vD_vm_0}{4} |\nablala \phi_v^k|_Y^2 + C\big(|\phi_{v,\infty}'|^2 +  |\phi_v^k-\phi_{v,\infty}|^2_Y \big)    + \frac{1}{4\wdk} |v_v^k|_Y^2
          \\ &\qquad + (1+C_P^2\wdk^{-1})  R^2C_2 \max\{L_p;L_\sigma\}   \big( |\phi_\sigma^k|_V^2 +|\phi_v^k|_Y^2 +  K^{-2} |v^k|_H^2 + |S_p^k|_H^2 + \wdk^{-2} |v_v^k|_Y^2 + |p_\infty|_V^2 + |p_{v,\infty}|^2   \big),
\end{aligned}\end{equation*}
which gives after choosing $C_v > \frac{4 \wdk}{D_v m_0}$ and absorbing 
\begin{equation} 
\label{Eq:Final1D}
\begin{aligned} &\frac{C_v}{2} \ddt |\phi_v^k-\phi_{v,\infty}|_Y^2
   + \frac{3C_v D_v m_0}{4} |\nablala \phi_v^k|_Y^2 + \left(\frac{1}{2\wdk}  - \frac{(1+C_P^2\wdk^{-1})R^2 C_2 \max\{L_p;L_\sigma\}}{\wdk^2} \right) |v_v^k|_Y^2\\
    &\leq (1+C_P^2\wdk^{-1})R^2C_2 \max\{L_p;L_\sigma\} \big( |\phi_\sigma^k|_V^2 +|\phi_v^k-\phi_{v,\infty}|_Y^2 + |\phi_{v,\infty}|^2 +  K^{-2} |v^k|_H^2 + |S_p^k|_H^2 
\\ &\qquad + |p_\infty|_V^2 + |p_{v,\infty}|^2  \big) +C\big( |\phi_{v,\infty}'|^2 +  |\phi_v^k-\phi_{v,\infty}|^2_Y    \big).
    \end{aligned}
\end{equation}

	\subsubsection*{Step 3 (Adding)}

We add the equations \cref{Eq:FinalCH}, \cref{Eq:FinalRD} and \cref{Eq:Final1D} to arrive at
\begin{equation} \label{Eq:Adding} \begin{aligned}
	&   \frac12 \ddt \bigg[2  |\Psi^k|_{L^1(\Omega)}   +  \eps^2_P |\nabla \phi_P^k|_H^2
   +\eps^2_H  |\nabla \phi_H^k|_H^2 + \eps_N^2 |\nabla \phi_N^k|_H^2 + C_\sigma |\phi_\sigma^k|_H^2  +|\mde^k|_H^2+|\taf^k|_H^2 \\
   &\qquad + |\ecm^k|_V^2 + C_v|\phi_v^k-\phi_{v,\infty}|_Y^2  \bigg]  +\frac{m_0}{2} |\nabla \mu_P^k|_H^2 +   \frac{m_0}{2} |\nabla \mu_H^k|_H^2 +\frac12 |\mu_P^k|_H^2 + \frac12 |\mu_H^k|_H^2   \\
   &\qquad +  \left( \frac{1}{2K} - \frac{C_1L_p}{K^2}  - \frac{C_\sigma}{D_\sigma m_0} - \frac{ (C_\sigma+R^2+R^2C_P^2 \wdk^{-1}) C_2 \max\{L_p;L_\sigma\}}{K^2} \right) |v^k|_H^2
   \\ &\qquad  + \left(\frac{C_\sigma D_\sigma m_0}{2} - \frac{4m_\infty^2 \chi_c^2}{m_0} - (C_\sigma+\wdk+1) C_2 \max \{L_p;L_\sigma\}  \right)  |\nabla \phi_\sigma^k|_H^2 
  \\ &\qquad  + D_{\MDE} m_0| \nabla \mde^k|_H^2  + D_{\TAF} m_0 |\nabla \taf^k|_H^2  + \frac{3C_vD_v m_0}{4} |\nablala \phi_v^k|_Y^2   \\ &\qquad    
   +  \left( \frac{1}{2\wdk} - \frac{C_1L_p}{\wdk^2} - \frac{(C_\sigma +R^2+R^2C_P^2\wdk^{-1}) C_2 \max\{L_p;L_\sigma\}}{\wdk^2}  \right) |v_v^k|_Y^2
   \\ &\leq 
    \eps \big(|\nabla S_N^k|_H^2+  |\nabla S_{\ECM}^k|_H^2 \big) + (C_1 + (C_\sigma +R^2+R^2C_P^2 \wdk^{-1}) C_2) \max\{L_p;L_\sigma\} |S_p^k|_H^2
	+ C \big( 1+|\phi_P|_V^2\\ &\qquad  + |\phi_H^k|_V^2 + |\phi_N^k|_V^2 + \blue{|\phi_v^k|_Y^2}+ |\phi_\sigma^k|_H^2  + |\mde^k|_H^2 
	  + |\taf^k|_H^2 
   + |\ecm^k|_V^2 + |\phi_v - \phi_{v,\infty}|_Y^2 + |p_\infty|_V^2 + |p_{v,\infty}|^2  \\ &\qquad  + |\phi_{v,\infty}|^2+|\phi_{v,\infty}'|^2+|S_P^k|_H^2 + |S_H^k|_H^2+  |S_N^k|_H^2 +|S_\sigma^k|_H^2 + |S_{\MDE}^k|_H^2
	+ |S_{\TAF}^k|_H^2 + |S_{\ECM}^k|_H^2 \big).
\end{aligned} \end{equation}
By assumption \cref{Ass:Source} on the source functions, we have the three estimates
\begin{itemize} \itemsep1em
	\item $\displaystyle \sum_{\alpha \in \A} |S_\alpha^k|_H^2
	\lesssim  \sum_{\alpha \in \A} |\phi_\alpha|_H^2$,
	\item $\displaystyle|\nabla S_N^k|_H^2 +  |\nabla S_{\ECM}^k|_H^2 \leq 2^{|\A|} f_\infty^2 |\A|^2 \sum_{\alpha \in \A} |\nabla \phi_\alpha|_H^2$,
	\item $\displaystyle|S_p^k|_H^2 \leq 8 \big( |\nabla \mu^k_P|^2_H +  |\nabla \mu^k_H|^2_H+  2\chi_c^2 |\nabla \phi_\sigma^k|_H^2 + 2\chi_h^2 |\nabla \ecm^k|_H^2  \big)$,
	\end{itemize}
and insert these estimates into \cref{Eq:Adding}. Further, in order to treat the factor $|\nabla S_N^k|_H^2 +  |\nabla S_{\ECM}^k|_H^2$, we choose the constant
\begin{equation} \label{Eq:Epsilon} \eps=\frac{m_0}{2^{|\A|+2} f_\infty^2 |\A|^2} \min\{C_\sigma D_\sigma,D_{\MDE},D_{\TAF}\}, \end{equation}
so that we can conclude
\begin{equation*}
	\begin{aligned}
	\eps \big(|\nabla S_N^k|_H^2 +  |\nabla S_{\ECM}^k|_H^2 \big)&\leq \frac{m_0}{4}  \min\{C_\sigma D_\sigma,D_{\MDE},D_{\TAF}\} \sum_{\alpha \in \A} |\nabla \phi_\alpha|_H^2 \\ &\leq \frac{C_\sigma D_\sigma m_0}{4} |\nabla \phi_\sigma^k|_H^2 + \frac{D_{\MDE} m_0}{4} |\nabla \mde^k|_H^2 + \frac{D_{\TAF} m_0}{4} |\nabla \taf^k|_H^2 \\
	&\qquad + C (|\phi_P^k|_V^2+|\phi_H^k|^2+|\phi_N^k|^2+|\ecm^k|_V^2).
	\end{aligned}
\end{equation*}
We absorb, collect and summarize the constants, giving
\begin{equation} \label{Eq:Collect} \begin{aligned}
	&  \frac12  \ddt \bigg[ 2 |\Psi^k|_{L^1(\Omega)}   +  \eps^2_P |\nabla \phi_P^k|_H^2
   + \eps^2_H  |\nabla \phi_H^k|_H^2 +  \eps_N^2  |\nabla \phi_N^k|_H^2 +C_\sigma |\phi_\sigma^k|_H^2  +|\mde^k|_H^2+|\taf^k|_H^2 
   \\ &\qquad + |\ecm^k|_V^2 + |\phi_v^k-\phi_{v,\infty}|_Y^2 \bigg]+\frac12 |\mu_P^k|_H^2 + \frac12 |\mu_H^k|_H^2
   \\ &\qquad  +\left( \frac{m_0}{2} - 8(C_1 + (C_\sigma +R^2+R^2 C_P^2 \wdk^{-1}) C_2) \max\{L_p;L_\sigma\} \right) \big(|\nabla \mu_P^k|_H^2 +   |\nabla \mu_H^k|_H^2 \big) 
   \\  &\qquad +  \left( \frac{1}{2K} - \frac{C_1L_p}{K^2}  - \frac{C_\sigma}{D_\sigma m_0} - \frac{(C_\sigma +R^2+R^2 C_P^2 \wdk^{-1}) C_2 \max\{L_p;L_\sigma\}}{K^2} \right) |v^k|_H^2
      \\ &\qquad  + \left(\frac{C_\sigma D_\sigma m_0}{2} - \frac{4m_\infty^2 \chi_c^2}{m_0} - (C_\sigma +R^2+R^2 C_P^2 \wdk^{-1}) C_2 \max \{L_p;L_\sigma\}  \right)  |\nabla \phi_\sigma^k|_H^2 
 \\ &\qquad + \frac{3m_0}{4}\big(D_{\MDE}  | \nabla \mde^k|_H^2  + D_{\TAF} |\nabla \taf^k|_H^2  + D_v  |\nablala \phi_v^k|_Y^2\big)
 \\ &\qquad +\left( \frac{1}{2\wdk} - \frac{C_1L_p}{\wdk^2} - \frac{(C_\sigma +R^2+R^2 C_P^2 \wdk^{-1}) C_2 \max\{L_p;L_\sigma\}}{\wdk^2}  \right) |v_v^k|_Y^2 
   \\ &\leq  C \big( 1+    |\phi_P^k|_V^2 + |\phi_H^k|_V^2 + |\phi_N^k|_V^2 +\blue{|\phi_v^k|_Y^2}+ |\phi_\sigma^k|_H^2   +|\mde^k|_H^2+|\taf^k|_H^2 + |\ecm^k|_V^2  +|\phi_v^k-\phi_{v,\infty}|_Y^2    
   \\ &\qquad +  |p_\infty|_V^2 + |p_{v,\infty}|^2 + |\phi_{v,\infty}|^2 + |\phi_{v,\infty}'|^2 \big),
\end{aligned} \end{equation}
and we choose $C_\sigma$ and $L_p$, $L_\sigma$, $K$ such that the prefactors are positive, see also assumption \cref{Ass:Jv}. In particular, we have to ensure the condition
$$\frac{8m_\infty^2 \chi_c^2}{m_0^2 D_\sigma} < C_\sigma < \frac{D_\sigma m_0}{2K}.$$

	\subsubsection*{Step 4 (Gr\"onwall--Bellman lemma)}
\blue{First, we eliminate the prefactors on the left hand side of the energy inequality \cref{Eq:Collect}  by estimating it with the minimum of all prefactors and bringing it to the right hand side to the generic constant $C$.} Afterwards, we integrate the inequality over the time interval $(0,t)$ with $t \in (0,T_k)$, apply the growth assumption \cref{Ass:Psi}, and obtain
\begin{equation*} \begin{aligned}
	&    |\Psi^k(t)|_{L^1(\Omega)}   +  |\phi_P^k(t)|_V^2
	+  |\phi_H^k(t)|_V^2 +    |\phi_N^k(t)|_V^2 + |\phi_\sigma^k(t)|_H^2  +|\mde^k(t)|_H^2+|\taf^k(t)|_H^2+ |\ecm^k(t)|_V^2  
	\\ &\qquad  + |\phi_v^k(t)-\phi_{v,\infty}|_Y^2  +  \|\phi_\sigma^k\|_{L^2(0,T_k;V)}^2 + \| \mu_P^k\|_{L^2(0,T_k;V)}^2 +   \|\mu_H^k\|_{L^2(0,T_k;V)}^2  +  \| \mde^k\|_{L^2(0,T_k;V)}^2   \\ & \qquad + \| \taf^k\|_{L^2(0,T_k;V)}^2  
	+  \|\phi_v^k - \phi_{v,\infty} \|_{L^2(0,T_k;X_0)}^2 
	+ \|v^k\|_{L^2(0,T_k;H)}^2 + \|v_v^k\|_{L^2(0,T_k;Y)}^2 \\ & \qquad - C \big( \|\phi_P^k\|_{L^2(0,T_k;V)}^2 + \|\phi_H^k\|_{L^2(0,T_k;V)}^2 + \|\phi_N^k\|_{L^2(0,T_k;V)}^2 +\blue{\|\phi_v^k\|_{L^2(0,T_k;Y)}^2}+ \|\phi_\sigma^k\|_{L^2(0,T_k;H)}^2 \\ &\qquad  +\|\mde^k\|_{L^2(0,T_k;H)}^2 +\|\taf^k\|_{L^2(0,T_k;H)}^2   +\|\ecm^k\|_{L^2(0,T_k;V)}^2  +\|\phi_v^k-\phi_{v,\infty}\|_{L^2(0,T_k;Y)}^2   \big)
	\\ &\leq C(T_k) \cdot \big (  1+   |\Psi^k(0)|_{L^1(\Omega)}   +  |\nabla \phi_{P,0}^k|_H^2
	+   |\nabla \phi_{H,0}^k|_H^2 +    |\nabla \phi_{N,0}^k|_H^2 + |\phi_{\sigma,0}^k|_H^2  +|\phi_{\MDE,0}^k|_H^2+|\phi_{\TAF,0}^k|_H^2  
	  \\ &\qquad + |\phi_{\ECM,0}^k|_V^2+ |\phi_{v,0}^k|_Y^2 + \|p_\infty\|_{L^2(0,T;V)}^2 + |p_{v,\infty}|_{L^2(0,T)}^2+ |\phi_{v,\infty}|_{H^1(0,T)}^2  \big)  .
\end{aligned} \end{equation*}

By applying the Gr\"onwall--Bellman lemma, see \cref{Lem_Gronwall}, we obtain
\begin{equation} \label{Eq:Gronwall} \begin{aligned} 
		 & \|\Psi^k\|_{L^\infty(0,T_k;L^1(\Omega))} +  \sum_{\mathclap{\alpha \in \CH \cup \{\ECM\}}} ~ \| \phi_\alpha^k\|^2_{L^\infty(0,T_k;V)} + ~\sum_{\mathclap{\alpha \in \{P,H\}}} ~ \| \mu_\alpha^k\|^2_{L^2(0,T_k;V)} + \sum_{\beta \in \RD}  
		\|\phi_\beta^k\|^2_{L^\infty(0,T_k;H) \cap L^2(0,T_k;V)}  
		\\ &\qquad + \|v^k\|_{L^2(0,T_k;H)}^2  +\|\phi_v^k-\phi_{v,\infty}\|^2_{L^\infty(0,T_k;Y) \cap L^2(0,T_k;X_0)}  + \|v_v^k\|_{L^2(0,T_k;Y)}^2  \\ &\leq C(T_k) \cdot \Big( 1+|\phi_{v,0}^k|^2_{Y}+ \sum_{\mathclap{\alpha \in \CH \cup \{\ECM\}}} ~ | \phi_{\alpha,0}^k|^2_{V} + \sum_{\beta \in \RD}  
		|\phi_{\beta,0}^k|^2_{H}   +  |\Psi(\phi_{P,0}^k,\phi_{H,0}^k,\phi_{N,0}^k)|_{L^1(\Omega)} + \|p_\infty\|_{L^2(0,T;V)}^2
		\\[-.2cm] &\qquad  + |p_{v,\infty}|_{L^2(0,T)}^2 + |\phi_{v,\infty}|_{H^1(0,T)}^2  \Big).
	\end{aligned} \end{equation}

We have chosen the initial values of the Faedo--Galerkin approximations as the orthogonal projections of the initial values of their counterpart, see \cref{Eq:GalerkinInitial}. The operator norm of an orthogonal projection is bounded by $1$ and, therefore, uniform estimates are obtained in \cref{Eq:Gronwall}; for example
$$|\phi_{P,0}^k|_V^2 = |\Pi_{H_k} \phi_{P,0}|_V^2 \leq |\phi_{P,0}|_V^2.$$
Using the upper bound \cref{Eq:PsiGrowth} of $\Psi$, we treat the term involving the potential function on the right hand side in the following way: 
$$\begin{aligned} 
	|\Psi(\phi_{P,0}^k,\phi_{H,0}^k,\phi_{N,0}^k)|_{L^1(\Omega)} &\lesssim 1+  |\phi_{P,0}^k|_H^2+|\phi_{H,0}^k|_H^2+|\phi_{N,0}^k|_H^2 \\  &= 1+  |\Pi_{H_k} \phi_{P,0}|_H^2|+|\Pi_{H_k} \phi_{H,0}|_H^2+|\Pi_{H_k} \phi_{N,0}|_H^2 \\ &\leq 1+  |\phi_{P,0}|_H^2+|\phi_{H,0}|_H^2+|\phi_{N,0}|_H^2.
\end{aligned}$$
Now, the $k$-independent right hand side in the estimate allows us to extend the time interval by setting $T_k=T$ for all $k \in \mathbb{N}$.
Therefore, we have the final uniform energy estimate,
\begin{equation}\label{Eq:FinalEnergy}
	\begin{aligned}
			& \|\Psi^k\|_{L^\infty(0,T;L^1(\Omega))} + \sum_{\mathclap{\alpha \in \CH \cup \{\ECM\}}} ~ \| \phi_\alpha^k\|^2_{L^\infty(0,T;V)} + ~\sum_{\mathclap{\alpha \in \{P,H\}}} ~ \| \mu_\alpha^k\|^2_{L^2(0,T;V)} + \sum_{\mathclap{\beta \in \RD}}  
		   \|\phi_\beta^k\|^2_{L^\infty H \cap L^2(0,T;V)}   \\ &\qquad + \|v^k\|_{L^2(0,T;H)}^2+\|\phi_v^k-\phi_{v,\infty}\|^2_{L^\infty(0,T;Y) \cap L^2(0,T;X_0)}  + \|v_v^k\|_{L^2(0,T;Y)}^2  \\ &\leq C(T)\! \cdot \! \Big(1+|\phi_{v,0}|^2_{Y}+ \sum_{\mathclap{\alpha \in \CH \cup \{\ECM\}}} ~ | \phi_{\alpha,0}|^2_{V} + \sum_{\mathclap{\beta \in \RD}}  
		   |\phi_{\beta,0}|^2_{H} +\|p_\infty\|_{L^2(0,T;V)}^2
		     + |p_{v,\infty}|_{L^2(0,T)}^2 + |\phi_{v,\infty}|_{H^1(0,T)}^2   \Big).
	   \end{aligned} \end{equation}
	   From this energy inequality and \cref{Eq:ModelGalerkinVelocity} we also get bounds for the pressures $p^k$ and $p_v^k$ in the following way $$\|p^k-p_\infty\|_{L^2(0,T;V_0)}  + \|p_v^k - p_{v,\infty} \|_{L^2(0,T;X_0)}  \leq C.$$

\subsection{Limit process} \label{Sec:Limit}

\subsubsection*{Weak convergence}
Next, we prove that there are subsequences of $\phib^k,\mu_P^k,\mu_H^k,p^k,\phi_v^k,p_v^k$, which converge to a weak solution of our model \cref{Eq:Model3D}--\cref{Eq:Model1D} in the sense of \cref{Def:Weak}. From the energy estimate \cref{Eq:FinalEnergy} we deduce that
\begin{equation}\begin{aligned}
		\{\phi_\alpha^k\}_{k \in \mathbb{N}}      & \text{ is bounded in } L^\infty(0,T;V),        &&\alpha \in \CH \cup \{\ECM\}            ,            \\
		\{\mu_\alpha^k\}_{k \in \mathbb{N}}         & \text{ is bounded in } L^2(0,T;V),       &&\alpha \in \{P,H\},                                    \\
		\{\phi_\beta^k\}_{k \in \mathbb{N}} & \text{ is bounded in } L^\infty(0,T;H) \cap L^2(0,T;V), &&\beta \in \RD,     \\
		\{v^k\}_{k \in \mathbb{N}}           & \text{ is bounded in }  L^2(0,T;L^2(\Omega;\R^3)), \\
				\{p^k\}_{k \in \mathbb{N}}           & \text{ is bounded in }  (p_\infty + L^2(0,T;V_0)), \\
		\{\phi_v^k\}_{k \in \mathbb{N}} & \text{ is bounded in } L^\infty(0,T;Y) \cap (\phi_{v,\infty} + L^2(0,T;X_0)),     \\
				\{v_v^k\}_{k \in \mathbb{N}}           & \text{ is bounded in }  L^2(0,T;Y), \\
		\{p_v^k\}_{k \in \mathbb{N}}           & \text{ is bounded in }  (p_{v,\infty} + L^2(0,T;X_0)),
	\end{aligned} \label{Eq:Bounds1} \end{equation}
and, by the Banach--Alaoglu theorem, these bounded sequences have weakly/weakly-$*$ convergent subsequences. By a standard abuse of notation, we drop the subsequence index. Consequently, there are functions $\phib:(0,T)\times \Omega \to \R^{|\A|}$, $\mu_P, \mu_H, p :(0,T)\times \Omega \to \R$, $v:(0,T) \times \Omega \to \R^3$, $ \phi_v, v_v,p_v:(0,T)\times \Lambda \to \R$ such that, for $k \to \infty$,
\begin{equation}
	\begin{alignedat}{3}
		\phi_\alpha^k &\longweak \phi_\alpha &&\text{ weakly-$*$ in } L^\infty(0,T;V), &&\quad \alpha \in \CH \cup \{\ECM\},        
		\\
		\mu_\alpha^k &\longweak \mu_\alpha &&\text{ weakly\phantom{-*} in } L^2(0,T;V), &&\quad \alpha \in \CH \backslash \{N\},      \\
		\phi_\beta^k &\longweak \phi_\beta &&\text{ weakly-$*$ in } L^\infty(0,T;H) \cap L^2(0,T;V), &&\quad \beta \in \RD, \\
		v^k &\longweak v &&\text{ weakly\phantom{-*} in } L^2(0,T;L^2(\Omega;\R^3)) , \\
		p^k &\longweak p &&\text{ weakly\phantom{-*} in } (p_\infty +L^2(0,T;V_0)) , \\
		\phi_v^k &\longweak \phi_v &&\text{ weakly-$*$ in } L^\infty(0,T;Y) \cap (\phi_{v,\infty} + L^2(0,T;X_0)), \\
				v_v^k &\longweak v_v &&\text{ weakly\phantom{-*} in }  L^2(0,T;Y), \\
		p_v^k &\longweak p_v &&\text{ weakly\phantom{-*} in } (p_{v,\infty} + L^2(0,T;X_0)).
	\end{alignedat}
	\label{Eq:WeakConv}
\end{equation}

\subsubsection*{Strong convergence} We now consider taking the limit $k \to \infty$ in the Faedo--Galerkin system \cref{Eq:ModelGalerkin3D}--\cref{Eq:ModelGalerkin1D} in the hope to attain the initial variational system \cref{Eq:ModelWeak3D}--\cref{Eq:ModelWeak1D}. Since the equations in \cref{Eq:ModelGalerkin3D}--\cref{Eq:ModelGalerkin1D} are nonlinear in $\phib^k$ and $\phi_v^k$, we want to achieve strong convergence of these sequences before we take the limit in \cref{Eq:ModelGalerkin3D}--\cref{Eq:ModelGalerkin1D}. Therefore, our goal is to bound their time derivatives and to apply the Aubin--Lions--Simon compactness lemma \cref{Eq:EmbeddingComp}.

Let $(\varphi,\hat \varphi, \tilde \varphi)$ be such that $\varphi \in L^2(0,T;V)$, $\hat \varphi \in L^2(0,T;H)$, $\tilde \varphi  \in L^2(0,T;X_0)$, and
$$\Pi_{H_k}\varphi = \sum_{j=1}^k \varphi_j^k h_j, \quad \Pi_{H_k}\hat \varphi = \sum_{j=1}^k \hat \varphi_j^k h_j, \quad \Pi_{Y_k} \tilde \varphi = \sum_{j=1}^k \tilde \varphi_j^k y_j,$$
with time-dependent coefficient functions $\varphi_j^k, \hat \varphi_j^k, \tilde \varphi_j^k : (0,T) \to \R$, $j \in \{1,\dots,k\}$.  We multiply the equations \cref{Eq:ModelGalerkin3D} and \cref{Eq:ModelGalerkin1D} by $\tilde \varphi_j^k$ by the appropriate coefficient functions, sum up each equation from $j=1$ to $k$ and integrate in time over $(0,T)$, to obtain the equation system,
\begin{equation} \label{Eq:Derivative}  \begin{aligned}
	\int_0^T \langle \pt \phi_\alpha^k , \varphi \rangle_V \dt &= \int_0^T (\C(\phi_\alpha^k) v^k,\nabla \Pi_{H_k} \varphi)_H -(m_\alpha^k \nabla \mu_\alpha^k,\nabla \Pi_{H_k} \varphi)_H                             + (S_\alpha^k,\Pi_{H_k} \varphi)_H \dt,
	\\
	\int_0^T \langle \pt \phi_\beta^k , \hat \varphi \rangle_V     \dt                                                                                                                   & = \int_0^T (S_\beta^k,\Pi_{H_k} \hat \varphi)_H \dt,
	\\
	\int_0^T \langle \pt \phi_\sigma^k,\varphi \rangle_V \dt &=\int_0^T (\C(\phi_\sigma^k)v^k,\nabla\Pi_{H_k}\varphi)_H -  D_\sigma(m_\sigma^k \nabla \phi_\sigma^k, \nabla \Pi_{H_k}\varphi)_H + (S_\sigma^k,\Pi_{H_k}\varphi)_H \\ &\qquad+ \langle \delta_\Gamma, J_{\sigma v,\Gamma}^k \Pi_{H_k}\varphi\rangle_W              -\chi_c (m_\sigma^k \nabla (\phi_P^k +\phi_H^k +\phi_N^k ),\nabla \Pi_{H_k} \varphi)_H \dt,
	\\
	\int_0^T \langle \pt \phi_\gamma^k,\varphi\rangle_V \dt &= \int_0^T -  D_\gamma (m_\gamma^k \nabla \phi_\gamma^k,\nabla \Pi_{H_k} \varphi)_H                                                       + (S_\gamma^k,\Pi_{H_k} \varphi)_H \dt,
	\\
	\int_0^T \langle \pt \phi_v^k, \tilde \varphi\rangle_X \dt &= \int_0^T  (\C(\phi_v^k) v_v^k ,\blue{\nablala} \Pi_{Y_k}\tilde \varphi )_Y - D_v (m_v^k \nablala \phi_v^k,\nablala \Pi_{Y_k} \tilde \varphi)_Y  -\wdr(J_{\sigma v}^k,\Pi_{Y_k} \tilde \varphi)_Y \dt,
\end{aligned}
\end{equation}
where $\alpha \in \{P,H\}$, $\beta \in \{N,\ECM\}$, $\gamma \in \{\MDE,\TAF\}$. Each equation in \cref{Eq:Derivative} can be treated using standard inequalities \blue{and the estimate involving the trace operator, see \cref{Eq:EstimateDelta},} the boundedness of the orthogonal projection and the energy estimate \cref{Eq:FinalEnergy}, e.g., we find
\begin{equation*} \begin{aligned} \int_0^T \langle \pt \phi_\sigma^k , \varphi \rangle_V \dt &\lesssim \int_0^T  |v^k|_H |\varphi|_V +  | \phi_\sigma^k|_V |\varphi|_V + \sum_{\alpha \in \A} |\phi_\alpha^k|_H |\varphi|_H + |J_{\sigma v,\Gamma}^k|_{L^2(\Gamma)} |\varphi|_V + \sum_{\beta \in \CH } |\phi_\beta^k|_V |\varphi|_V  \dt
	\\ & \lesssim  \|\varphi\|_{L^2(0,T;V) }.
	\end{aligned} \label{Eq:DerivativeBoundCH}
\end{equation*}
From this inequality and the bounds derived earlier, see \cref{Eq:Bounds1}, we conclude that
\begin{alignat*}{3}
	 & \{\phi_\alpha^k\}_{k \in \mathbb{N}}          &  & \text{ is bounded in }  H^1(0,T;V') \cap L^\infty(0,T;V) ,                          && \qquad \alpha \in \{P,H\},                                     \\
	 & \{\phi_\beta^k\}_{k \in \mathbb{N}}          &  & \text{ is bounded in } H^1(0,T;H) \cap L^\infty(0,T;V),                          && \qquad \beta \in \{N,\ECM\},                                     \\
	 & \{\phi_\gamma^k\}_{k \in \mathbb{N}}     &  & \text{ is bounded in } H^1(0,T;V') \cap L^\infty(0,T;H) \cap L^2(0,T;V) ,     && \qquad \gamma \in \RD,              \\
	 & \{\phi_v^k\}_{k \in \mathbb{N}}     &  & \text{ is bounded in } H^1(0,T;X_0') \cap L^\infty(0,T;Y) \cap (\phi_{v,\infty} + L^2(0,T;X_0)) .               
\end{alignat*}
We apply the Aubin--Lions--Simon compactness lemma \cref{Eq:EmbeddingComp},
yielding the strong convergences as $k \to \infty$
\begin{equation} \label{Eq:StrongConv}
	\begin{alignedat}{3}
		\phi_\alpha^k &\longrightarrow \phi_\alpha &&\text{ strongly in } C^0([0,T];H), && \qquad \alpha \in \CH \cup \{\ECM\},\\
		\phi_\beta^k &\longrightarrow \phi_\beta &&\text{ strongly in } L^2(0,T;H)\cap C^0([0,T];V'),&& \qquad \beta \in \RD, \\
		\phi_v^k &\longrightarrow \phi_v &&\text{ strongly in } L^2(0,T;Y)\cap C^0([0,T];X_0').
	\end{alignedat}
\end{equation}
The strong convergence $\phi_\alpha^k \to \phi_\alpha$ in $C^0([0,T];H)$ implies $\phi_\alpha(0)=\phi_{\alpha,0}$ in $H$ and similarly $\phi_\beta(0)=\phi_{\beta,0}$ in $V'$ and $\phi_v(0)=\phi_{v,0}$ in $X'_0$. Therefore, the limit functions $(\phib,\phi_v)$ of the Faedo--Galerkin approximations fulfill the initial conditions for the system \cref{Eq:Model3D}--\cref{Eq:Initial1D}.

\subsubsection*{Limit process} We show that the limit functions also satisfy the variational form \cref{Eq:ModelWeak3D}--\cref{Eq:ModelWeak1D}, as defined in \cref{Def:Weak}. Multiplying the Faedo--Galerkin system \cref{Eq:ModelGalerkin3D}--\cref{Eq:ModelGalerkin1D} by $\eta \in C_c^\infty(0,T)$ and integrating from $0$ to $T$, gives
\begin{equation}\begin{aligned}
	\int_0^T \langle \pt \phi_\alpha^k , h_j \rangle_V \eta(t) \dt &= \int_0^T \big( (\C(\phi_\alpha^k) v^k,\nabla h_j)_H -(m_\alpha^k \nabla \mu_\alpha^k,\nabla  h_j)_H                             + (S_\alpha^k,h_j)_H\big) \eta(t) \dt,
	\\
	\int_0^T (\mu_\alpha^k,h_j)_H \eta(t) \dt &= \int_0^T \big( (\p_{\phi_\alpha} \Psi^k - \chi_c \phi_\sigma^k - \chi_h \ecm^k,h_j)_H + \eps_\alpha^2 (\nabla \phi_\alpha^k, \nabla h_j)_H \big) \eta(t) \dt, 
	\\
	\int_0^T ( \pt \phi_\beta^k , h_j )_H   \eta(t) \dt                                                                                                                    & = \int_0^T (S_\beta^k,h_j)_H \eta(t) \dt,	\end{aligned}
\end{equation}
 and 
\begin{equation}  \begin{aligned}
	\int_0^T \langle \pt \phi_\sigma^k,h_j \rangle_V \dt &=\int_0^T \big( (\C(\phi_\sigma^k) v^k,\nabla h_j)_H -  D_\sigma(m_\sigma^k \nabla \phi_\sigma^k, \nabla h_j)_H + (S_\sigma^k,h_j)_H \\ &\qquad+ \langle \delta_\Gamma, J_{\sigma v,\Gamma}^k h_j\rangle_W              -\chi_c (m_\sigma^k \nabla (\phi_P^k +\phi_H^k +\phi_N^k ),\nabla h_j)_H \big) \eta(t) \dt,
	\\
	\int_0^T \langle \pt \phi_\gamma^k,h_j\rangle_V  \eta(t) \dt &= \int_0^T \big( -  D_{\gamma} (m_{\gamma}^k \nabla \phi_\gamma^k,\nabla h_j)_H                                                       + (S_{\gamma}^k,h_j)_H \big) \eta(t) \dt,
	\\
		\int_0^T (v^k,h_j)_H \eta(t) \dt &= \int_0^T \big(-K(\nabla p,h_j)_H + (S_p,\nabla h_j)_H \big) \eta(t) \dt,  \\
	\int_0^T K(\nabla p^k,\nabla h_j)_H \eta(t) \dt &= \int_0^T \big( \langle \delta_\Gamma, J_{pv,\Gamma}^k h_j \rangle_W + (S_p,\nabla h_j)_H \big) \eta(t) \dt, 
\end{aligned}
\end{equation}
and 
\begin{equation}  \begin{aligned}
	\int_0^T \langle \pt \phi_v^k, y_j\rangle_X \eta(t)  \dt &= \int_0^T \big(   (\C( \phi_v^k) v_v^k ,\nabla y_j )_Y - D_v (m_v^k \nablala \phi_v^k,\nablala y_j)_Y  -\wdr(J_{\sigma v}^k,y_j)_Y \big) \eta(t) \dt, \\
		\int_0^T  (v_v^k, y_j)_Y \eta(t) \dt &= \int_0^T -  \wdk (\nablala p_v^k,y_j)_Y \eta(t) \dt, \\
	\int_0^T \wdk (\nablala p_v^k, \nablala y_j)_Y \eta(t) \dt &= \int_0^T -  \wdr(J_{pv}^k,y_j)_Y \eta(t) \dt,
\end{aligned}
\end{equation}
for each $j\in \{1,\dots,k\}$, $\alpha \in \{P,H\}$, $\beta \in \{N,\ECM\}$, $\gamma \in \{\MDE,\TAF\}$. We take the limit $k \to \infty$ in each equation. The linear terms can be treated directly in the limit process since they can be justified via the weak convergences \cref{Eq:WeakConv}, e.g.,the functional
$$\mu_P^k \mapsto \int_0^T (\mu_P^k,h_j)_H \eta(t) \, \dt \leq \|\mu_P^k\|_{L^2(0,T;H)} |h_j|_{H} |\eta|_{L^2(0,T)}, $$
is linear and continuous on $L^2(0,T;H)$ and hence, as $k\to \infty$,
$$\int_0^T (\mu_P^k,h_j)_H \eta(t) \, \dt  \longrightarrow \int_0^T (\mu_P,h_j) \eta(t) \, \dt.$$
Thus, it remains to examine the nonlinear terms. We do so in the steps (i)--(v) as follows. \\[0.1cm]
\begin{enumerate}[label=(\roman*), wide, labelwidth=!, labelindent=0pt] \itemsep1em
\item We have derived the convergence, see \cref{Eq:StrongConv}, $$\begin{aligned}\phi_\alpha^k \longrightarrow \phi_\alpha &\text{ in } L^2(0,T;H) \cong L^2((0,T)\times \Omega), &&\alpha \in \A,\end{aligned}$$ for $k \to \infty$ and, consequently, we have by the continuity and boundedness of $m_\alpha$, 
	$$m_\alpha^k=m_\alpha\big(\phib^k(t,x)\big) \longrightarrow m_\alpha\big(\phib(t,x)\big)=:m_\alpha \text{ a.e. in } (0,T)\times \Omega \text{ for } k \to \infty.$$ 
Applying the Lebesgue dominated convergence theorem, gives for $k \to \infty$
$$m_\alpha^k  \nabla h_j \eta \longrightarrow m_\alpha \nabla h_j \eta \text{ in } L^2((0,T) \times \Omega;\R^d),$$
and, together with $\nabla \mu_\alpha^k \weak \nabla \mu_\alpha$ weakly in $L^2((0,T)\times \Omega;\R^d)$ as $k\to\infty$, we have for $k \to \infty$
$$m_\alpha(\phib^k) \eta \nabla h_j \cdot \nabla \mu_\alpha^k \longrightarrow m_\alpha(\phib) \eta \nabla h_j \cdot \nabla \mu_\alpha \text{ in } L^1((0,T)\times \Omega).$$
We use here the fact that the product of a strongly and a weakly converging sequence in $L^2$ converges strongly in $L^1$.  \\[0.1cm]
\item By \cref{Eq:StrongConv}, we have $\phi_\alpha^k \to \phi_\alpha$ in $L^2((0,T)\times \Omega)$ and $v^k \weak v$ in $L^2((0,T)\times \Omega;\R^d)$ as $k\to \infty$, hence for $k \to \infty$ $$\C(\phi_\alpha^k) v^k  \cdot \nabla h_j \eta \longrightarrow \C(\phi_\alpha) v \cdot \nabla h_j \eta \text{ in } L^1((0,T)\times \Omega).$$
\item By the continuity and the growth assumptions on $\p_{\phi_\alpha} \Psi$, we have for $k \to \infty$
\begin{gather*} \p_{\phi_\alpha} \Psi\big(\phi_P^k(t,x),\phi_H^k(t,x),\phi_N^k(t,x)\big) \longrightarrow \p_{\phi_\alpha}\Psi\big(\phi_P(t,x),\phi_H(t,x),\phi_N(t,x)\big) \text{ a.e. in } (0,T)\times \Omega, \\
	|\p_{\phi_\alpha}\Psi(\phi_P^k,\phi_H^k,\phi_N^k) \eta h_j| \leq C(1+|\phi^k_P|+|\phi^k_H|+|\phi^k_N|)|\eta h_j|,
\end{gather*}
and the Lebesgue dominated convergence theorem yields for $k \to \infty$
$$\p_{\phi_\alpha}\Psi(\phi_P^k,\phi_H^k,\phi_N^k) \eta h_j \longrightarrow \p_{\phi_\alpha}\Psi(\phi_P,\phi_H,\phi_N) \eta h_j \text{ in } L^1((0,T)\times \Omega).$$
\item  We have the strong convergence of $\phi_P^k$ and $\phi_H^k$ in $L^2((0,T)\times \Omega)$ and the continuity and boundedness of $\C$. Together with the weak convergence of $\nabla \phi_\sigma^k$ and $\nabla \mu^k$ in $L^2((0,T)\times \Omega;\R^d)$ it is enough to conclude the convergence of the term involving $$S_p^k=-\C(\phi_P^k) (\nabla \mu_P^k+\chi_c \nabla \phi_\sigma^k)-\C(\phi_H^k) (\nabla \mu_H^k+\chi_c \nabla \phi_\sigma^k).$$
\item We have $\phi_v^k \to \phi_v$ in $L^2((0,T)\times \Lambda)$ and $\phi_\sigma^k \to \phi_\sigma$ in $L^2((0,T) \times \Omega)$ as $k \to \infty$ and therefore, also $\Pi_\Gamma \phi_v^k \to \Pi_\Gamma \phi_v$ in $L^2((0,T)\times \Gamma)$. Since $f_{\sigma,v}$ is a continuous and bounded function, we conclude
$$\int_0^T \int_\Gamma |f_{\sigma,v}(\phi_v^k,\Pi_\Gamma \phi_\sigma^k) \tr_\Gamma h_j \eta(t) |^2 \dS \dt \lesssim  \|f_{\sigma,v}(\Pi_\Gamma \phi_v^k,\phi_\sigma^k)\|_{L^\infty((0,T)\times \Gamma)}^2 |h_j|_{V}^2 |\eta|_{L^2(0,T)}^2,$$
and the Lebesgue dominated convergence theorem gives for $k\to \infty$
$$f_{\sigma,v}(\phi_v^k,\Pi_\Gamma \phi_\sigma^k) \tr_\Gamma h_j \eta(t) \longrightarrow f_{\sigma,v}(\phi_v,\Pi_\Gamma \phi_\sigma) \tr_\Gamma h_j \eta(t) \text{ in } L^2((0,T) \times \Gamma).  $$
Together with the weak convergence of $\Pi_\Gamma p_v^k$ and $p^k$ we have for $k \to \infty$
$$\begin{aligned} \int_0^T \langle \delta_\Gamma,J_{\sigma v,\Gamma}^k h_j \rangle_W \eta(t) \dt &=\int_0^T \int_\Gamma \Big( f_{\sigma,v}(\phi_\sigma^k,\Pi_\Gamma \phi_v^k) L_p (\Pi_\Gamma p_v^k-p^k) + L_\sigma (\Pi_\Gamma \phi_v^k - \phi_\sigma^k) \Big) \tr_\Gamma h_j  \eta(t) \dS \dt \\ &\to \int_0^T \langle \delta_\Gamma,J_{\sigma v,\Gamma} h_j \rangle_W \eta(t) \dt .
\end{aligned}$$
\end{enumerate}

Using the densities of $\cup_{k \in \mathbb{N}} H_k$ in $V$, $\cup_{k \in \mathbb{N}} H_k^0$ in $V_0$ and $\cup_{k \in \mathbb{N}} Y_k$ in $X$, and the fundamental lemma of the calculus of variations, we obtain a solution tuple $(\phib,\mu_P,\mu_H,v,p,\phi_v,v_v,p_v)$ to our model \cref{Eq:Model3D} and \cref{Eq:Model1D} in the weak sense as defined in \cref{Def:Weak}. 

\subsubsection*{Energy inequality} It remains to prove that found solution tuple satisfies the energy inequality \cref{Eq:EnergySolution}. First, we note that norms are weakly/weakly-$*$ lower semicontinuous, e.g.,we have $\mu_P^k \weak \mu_P$ in $L^2(0,T;V)$ and therefore, we infer
$$\| \mu_P \|_{L^2(0,T;V)} \leq \liminf_{k \to \infty} \|\mu_P^k\|_{L^2(0;T;V)}.$$  
We apply the Fatou lemma on the continuous and non-negative function $\Psi$ to obtain
$$\int_\Omega \Psi(\phi_P,\phi_H,\phi_N) \, \dd x  \leq \liminf_{k \to \infty} \int_\Omega \Psi(\phi_P^k,\phi_H^k,\phi_N^k) \, \dd x.$$ 
Consequently, passing the limit $k \to \infty$ in the discrete energy inequality \cref{Eq:FinalEnergy} leads to  \cref{Eq:EnergySolution}. 
\qed

\section{Numerical simulations}
\label{Sec:Simulation}


We present in this section two applications of the theory presented earlier. The first is the simple scenario in which two straight and idealized blood vessels are considered, one representing an artery and the other a vein. This means that the tissue block containing the two vessels is supplied nutrients by a single artery and drained of nutrients by a single vein. In between the two vessels, a tumor core is present which accepts nutrients injected through an inlet of the artery. The second scenario deals with a small blood vessel network described by data  in \cite{secomb2000theoretical} and on the following web page: \url{https://physiology.arizona.edu/sites/default/files/brain99.txt}. At four inlets of this network, nutrients are injected and transported through the network. As in the first scenario, the impact of the nutrients on a small tumor core surrounded by the network is investigated. We note that in all cases we consider a stationary vessel structure; for evolving and bifurcating vessels we refer to the further work \cite{fritz2020part2}.

To solve the one-dimensional partial differential eqautions \cref{Eq:Model1D} numerically, the vascular graph model (VGM) is employed \cite{reichold2009vascular,vidotto2019hybrid}. This method corresponds in principle to a vertex centered finite volume method with a two-point flux approximation. For the three-dimensional partial differential equations presented in \cref{sec:3Dmodel} a mixed finite volume-finite element discretization method is employed. The equations governing the pressure and nutrients in tissue, which are directly coupled with the one-dimensional system, are solved by a standard cell-centered finite volume scheme. Since the permeability of the cancerous and healthy tissue is given by a scalar field, a two-point flux approximation of the fluxes is used. 
For the remaining species, we consider a continuous and piecewise linear finite element approximation over a uniform cubic mesh. The coupled nonlinear partial differential equations are discretized in time using the semi-implicit Euler method. To solve the nonlinear system of equations arising in each time step, a fixed point iteration method is applied. We consider following double-well potential in free energy functional in \cref{Eq:GinzburgLandau} 
\begin{align}
	\Psi(\phi_P,\phi_H,\phi_N)&=C_{\Psi_T} \phi_T^2(1-\phi_T)^2.
\end{align}

\subsection{Tumor between two straight vessels} 
\label{Sec:Simulation2Straight}
We consider a tissue domain $\Omega = (0,2)^3$ containing two blood vessels aligned along the $z$-axis. The center lines of the vessels are located diagonally opposite to each other. 
The center line of the Vessel $1$ and $2$ pass through $(0.2, 0.2, 1)$ and $(1.8,1.8,1)$, respectively. 
We choose a radius of $R=0.08$ and $R = 0.1$ for Vessel $1$ and $2$. At the inlets of Vessel $1$, located at $(0.2, 0.2, 0)$ and  $(0.2, 0.2, 2)$, pressure values of $10000$ and $5000$ are prescribed. The inlets of Vessel $2$ are located at $(1.8, 1.8, 0)$ and $(1.8, 1.8, 2)$. Here, we consider the pressure values $1000$ and $2000$, respectively. Thus, Vessel $2$ will act as a vein taking up nutrients and blood plasma from the tissue domain. On the other hand Vessel $1$ has the function of an artery transporting nutrients into the tissue block $\Omega$.  \blue{We note that we choose the boundary values such that the velocities are sufficiently large in order to ensure that the transport processes are visible in the simulations.} Based on the pressure boundary conditions, we choose the boundary conditions for the nutrients as follows:
\ \\
\begin{itemize} \itemsep0em
	\item $\phi_v = \phi_{v,inlet} = 1$ at $(0.2, 0.2, 0)$.
	\ \\
	\item $\phi_v = 0$ at $(1.8, 1.8, 2)$.
	\ \\
	\item At all the remaining boundaries, we consider free flow boundary conditions.
\end{itemize}
\ \\
The initial tumor core is given by a ball of radius $0.3$ and centered at $(1,1,1)$. Within the tumor core, the total tumor volume fraction, $\phi_T$, decays smoothly  from $1$ in the center to $0$ on the boundary of the ball. Thereby, the necrotic and hypoxic volume fractions, $\phi_N$ and $\phi_H$, are set to zero. In the rest of the domain all the volume fractions for the tumor species are set to $0$ at $t=0$. The nutrient volume fraction, in the tissue domain $\Omega$, is initially fixed to a constant initial value of $0.6$, which is below the threshold values for prolific-to-hypoxic transition and above the threshold value for hypoxic-to-prolific transition. We also set $\phi_{\ECM} = 1$ at $t=0$. According to \cref{Eq:Initial3D} homogeneous Neumann boundary conditions are prescribed on $\partial \Omega$. 

As a simulation time period, the interval $\left(0,T \right)$ with $T=5$ is considered and the size of the time step is given by $\Delta t = 0.025$. The spatial domain $\Omega$ is discretized by cubic elements with an edge length of $h=0.025$.  We choose the parameters as listed in \cref{tab:params}.
\begin{table}[H]
	\begin{center}
		\caption{List of parameters and their values for the numerical experiments described in \cref{Sec:Simulation2Straight} and \cref{Sec:SimulationNetwork}. Parameters not mentioned below are set to zero.}\label{tab:params} \vspace{-.2cm}
				\begin{tabular}{|l | l || l | l || l | l |}
			\hline
			\textbf{Parameter} &
			\textbf{Value} & \textbf{Parameter} &
			\textbf{Value} & \textbf{Parameter} &
			\textbf{Value}  \\ \hline
			$\lambda_P$ & 5 & $\lambda_{P_h}$ & 0.5 & $\lambda_A$ & 0.005   \\
			$\lambda_{A_h}$ & 0.005  & $\lambda_{P\!H}$ & 1 & $\lambda_{H\!P}$ & 1 \\
			$\lambda_{H\!N}$ & 1 & $\sigma_{P\!H}$ & 0.55 & $\sigma_{H\!P}$ & 0.65  \\
			$\sigma_{H\!N}$ & 0.44 & $\varepsilon_P$ & 0.005 & $\varepsilon_H$ & 0.005\\
			$M_P$ & 50 & $M_H$ & 25 &  $\lambda_{\TAF_{\!P}}$ & 10 \\
			$D_{\TAF}$ & 0.5 & $m_{\TAF}$ & 1 & $L_p$  & $10^{-7}$ \\
			$D_\sigma$ & 1 & $m_\sigma$ & 1 & $K$ & $10^{-9}$ \\
			$D_{\MDE}$ & 0.5 & $m_{\MDE}$ & 1 &  $D_v$  & $0.1$ \\
			$\mu_{\text{bl}}$ & 1 & $L_\sigma $ & 10 & $\lambda_{\ECM_{\!D}}$ & 5 \\
			$\lambda_{\ECM_{\!P}}$ & 0.01 & $\lambda_{\MDE_{\!D}}$ & 1 & $\lambda_{\MDE_{\!D}}$ & 1 \\
			$\phi_{\ECM_{\!P}}$ & 0.5 & $C_{\Psi_T}$ & 0.045 & - & - \\
			\hline
		\end{tabular}
	\end{center}
\end{table}
\vspace{-.3cm}

\begin{figure}[H]
	\begin{center}		
		\includegraphics[width=.8\textwidth]{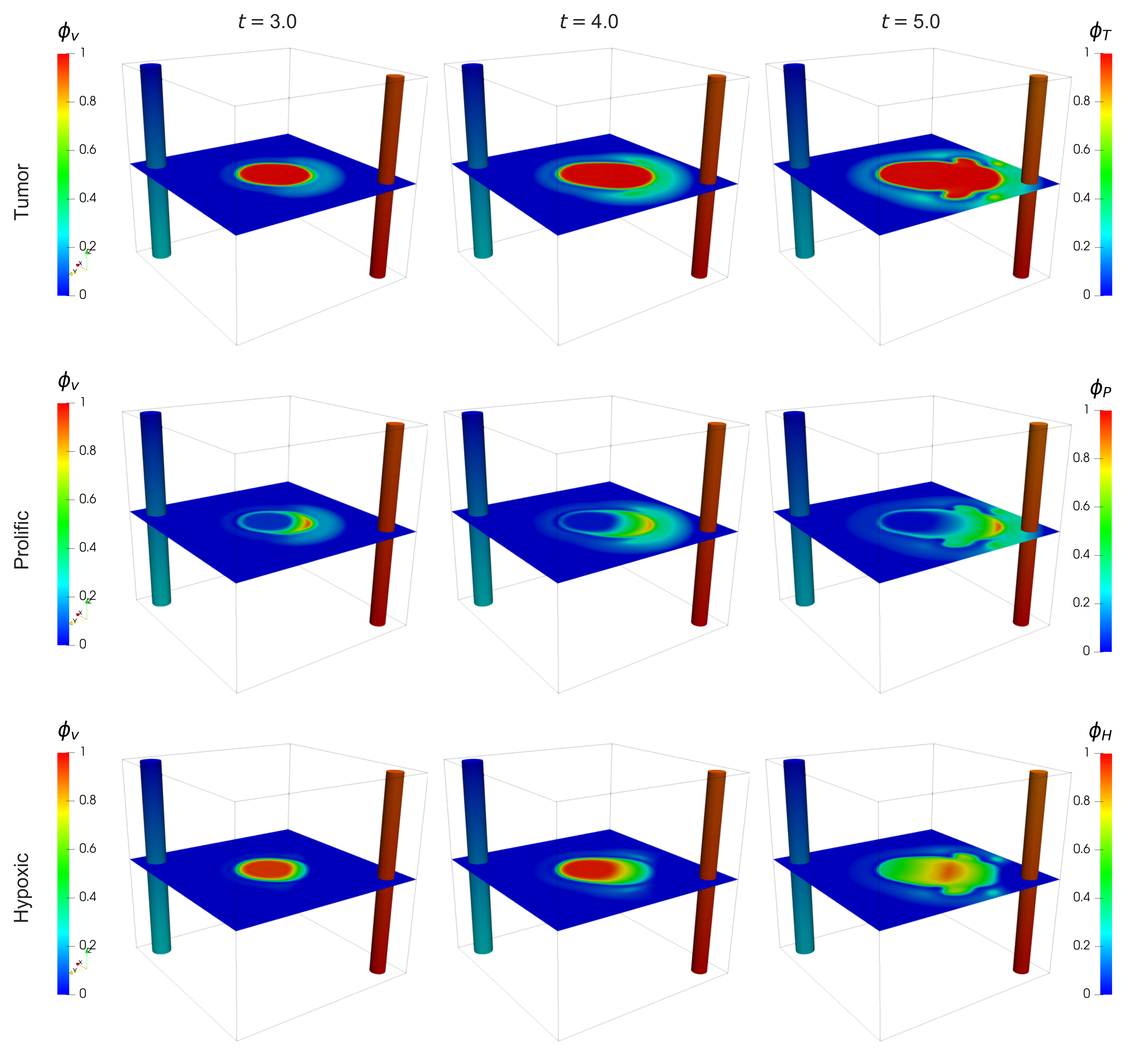}
	\end{center}
	\vspace{-0.4cm}
	\caption{\label{fig:two-vessels-vars} Evolution of the total tumor $\phi_T$ (top), prolific $\phi_P$ (middle), and hypoxic $\phi_H$ (bottom) volume fractions at the times $t \in \{3,4,5\}$ (left, middle,right) in the $z=1$ plane of the domain $\Omega$. On the two vessels the nutrients are described by the 1D constituent $\phi_v$. }
\end{figure}

\begin{figure}[H]
	\begin{center}		
	\hfill \includegraphics[width=.28\textwidth]{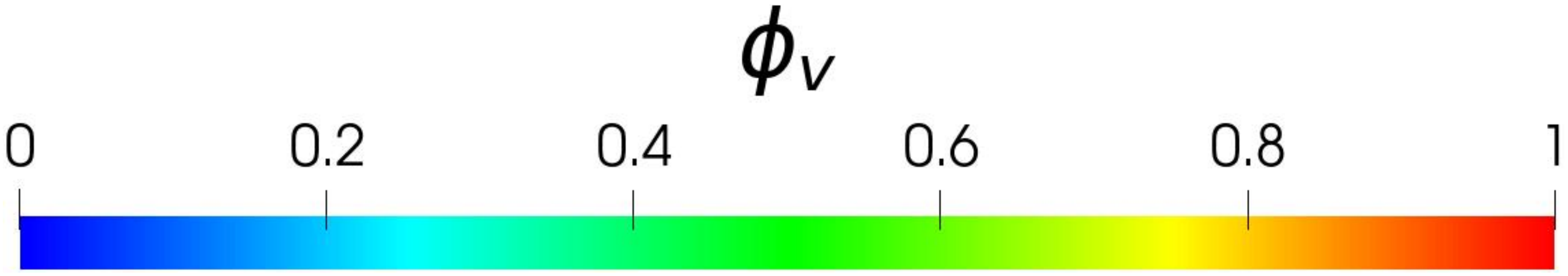} \hspace{1cm} ~  \\
	\includegraphics[page=1,height=.22\textheight,valign=t]{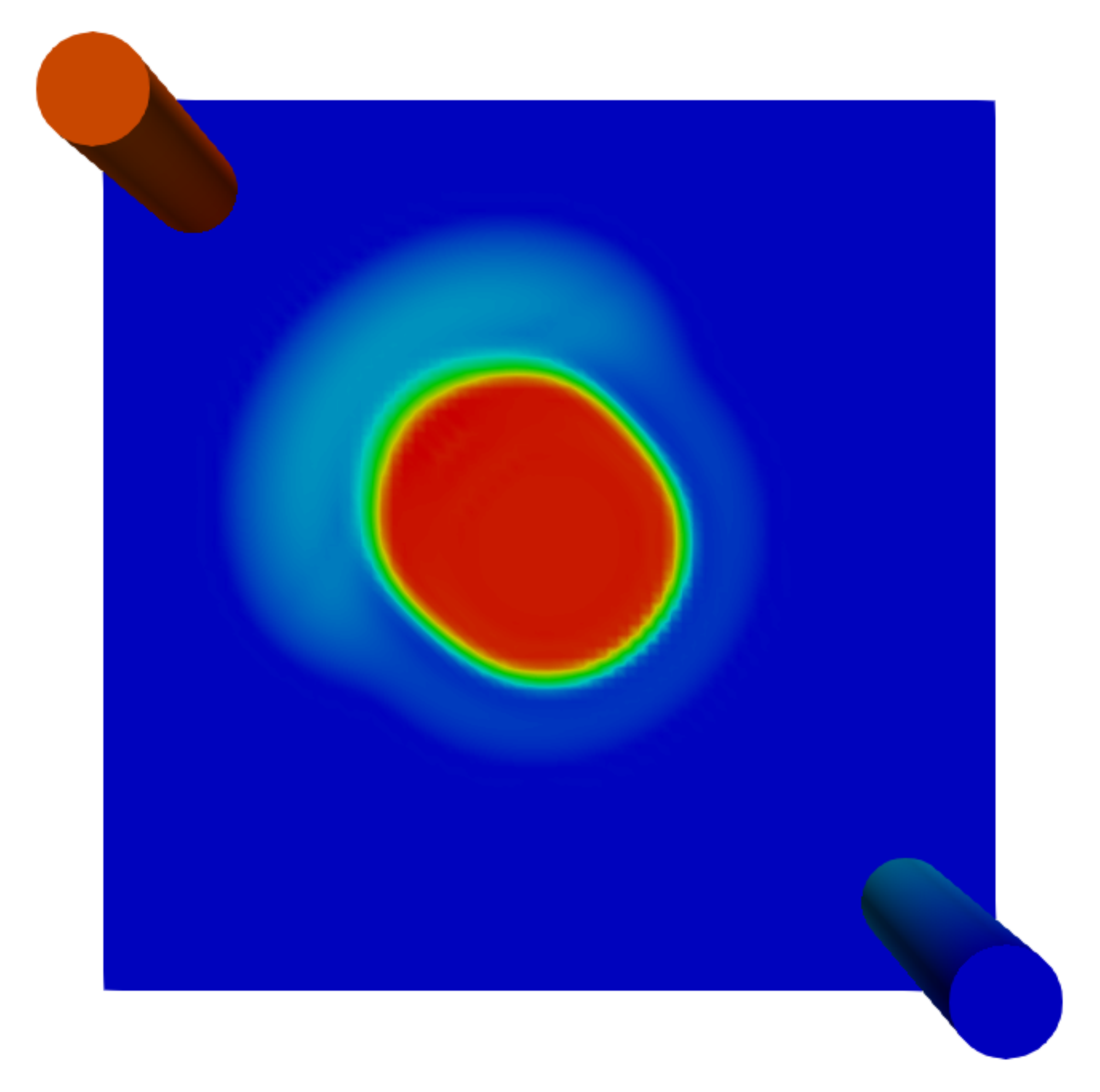} \hspace{-0.38cm}
	\includegraphics[page=2,height=.22\textheight,valign=t]{simulation2.pdf}
	\hspace{-0.38cm}
	\includegraphics[page=3,height=.22\textheight,valign=t]{simulation2.pdf}
		\includegraphics[height=.17\textheight,valign=t]{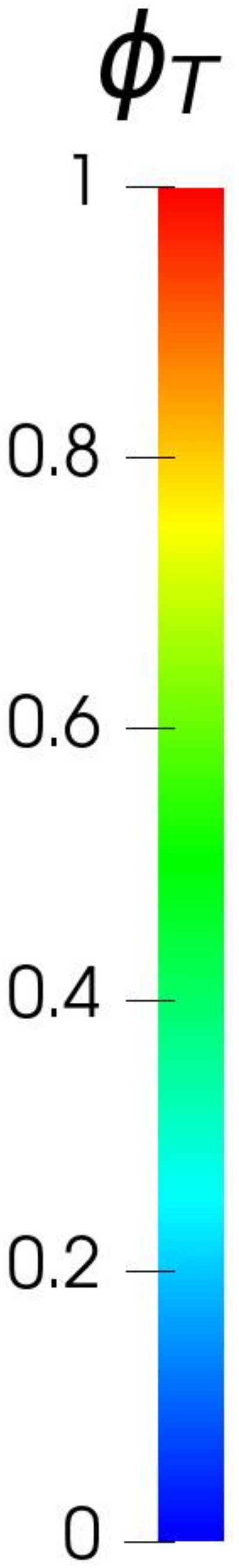}  \hfill
		\\
\includegraphics[width=\textwidth]{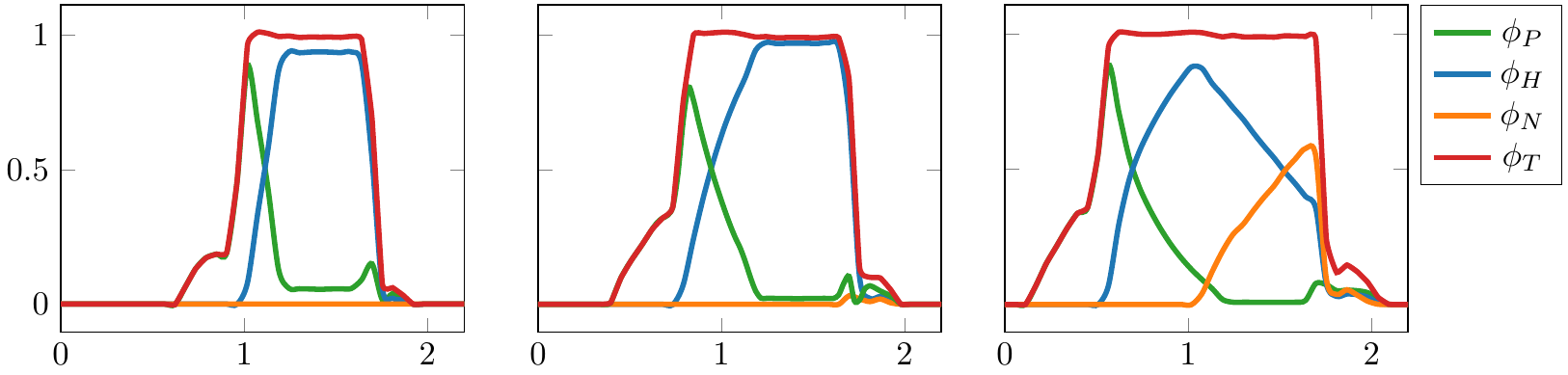}
	\end{center}
	\vspace{-0.4cm}
	\caption{\label{fig:two-vessels-line-plot} Left: Evolution of the tumor volume fraction ($\phi_T$) in the $z=1$ plane of the domain $\Omega$. Right: Plots of tumor ($\phi_T$), prolific ($\phi_P$), hypoxic ($\phi_H$) and necrotic ($\phi_N$) volume fractions along the line passing through the points $(0,0,1)$ and $(2,2,1)$ in the domain $\Omega$. From top row to bottom row, the plots correspond to time $t\in \{3,4,5\}$. }
\end{figure}

\begin{figure}[H]
	\begin{center}
		\includegraphics[width=.28\textwidth]{color1.pdf} \\
				\includegraphics[page=1,width=.3\textwidth]{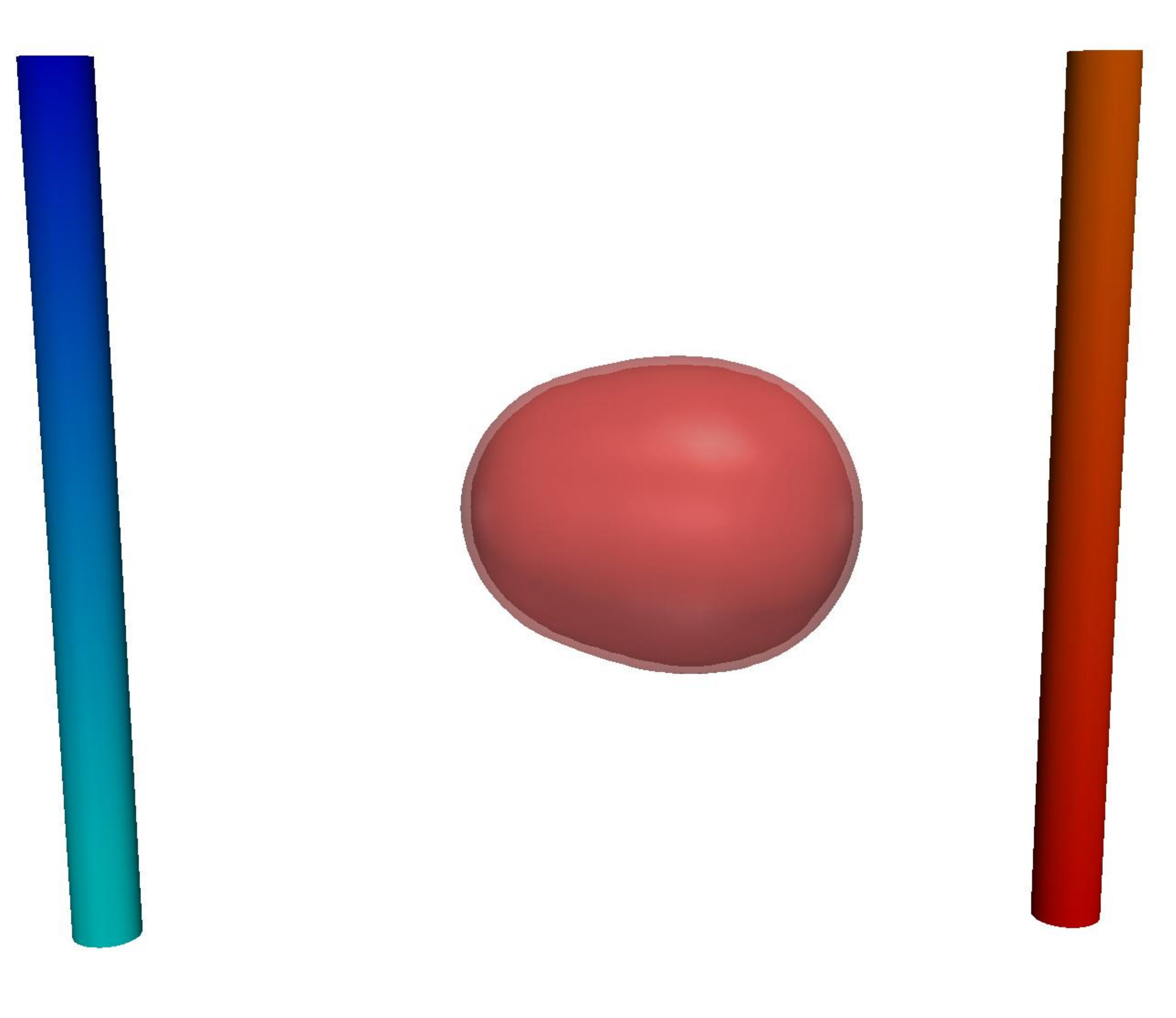}
				\includegraphics[page=2,width=.3\textwidth]{simulation4.pdf}	
				\includegraphics[page=3,width=.3\textwidth]{simulation4.pdf} \\[-.5cm]	
				\includegraphics[page=1,width=.3\textwidth]{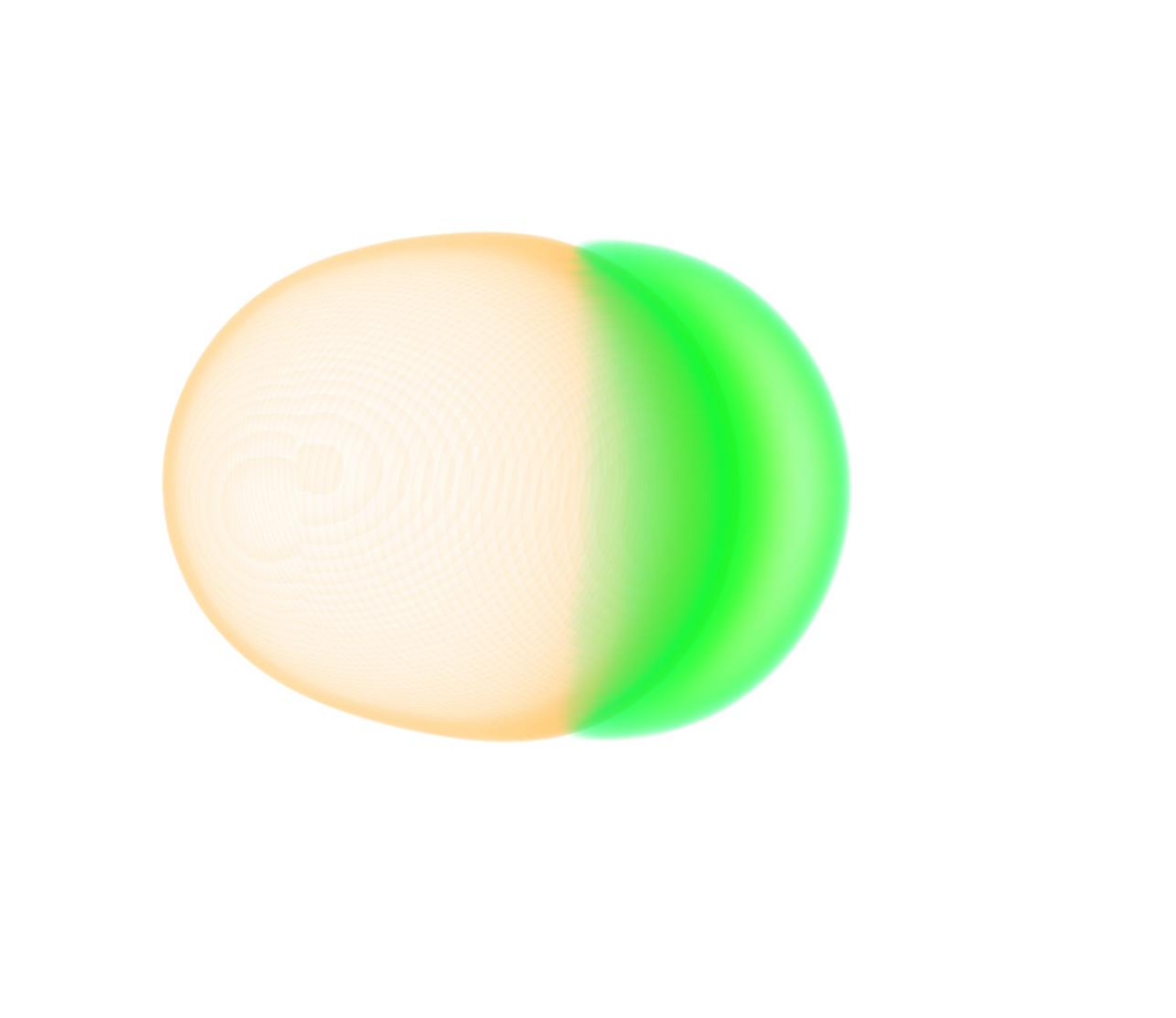}
\includegraphics[page=2,width=.3\textwidth]{simulation5.pdf}	
\includegraphics[page=3,width=.3\textwidth]{simulation5.pdf} 
%
	\end{center}
	\vspace{-0.4cm}
	\caption{\label{fig:two-vessels-multi-contour} Evolution of contour plots of the tumor volume fraction $\phi_T$ with vessels (top) at $0.8$ (light red) and $0.95$ (red) at the time $t \in \{4.25,4.75,5\}$ (left to right); the necrotic core is plotted at the contour line $\phi_N=0.42$ (black). Contour plots of the tumor phases without vessels (bottom) at $\phi_P=0.5$ (green), $\phi_H=0.45$ (orange), $\phi_N=0.4$ (dark red).}
\end{figure}

\begin{figure}[H]
	\begin{center}		
		\includegraphics[width=.9\textwidth]{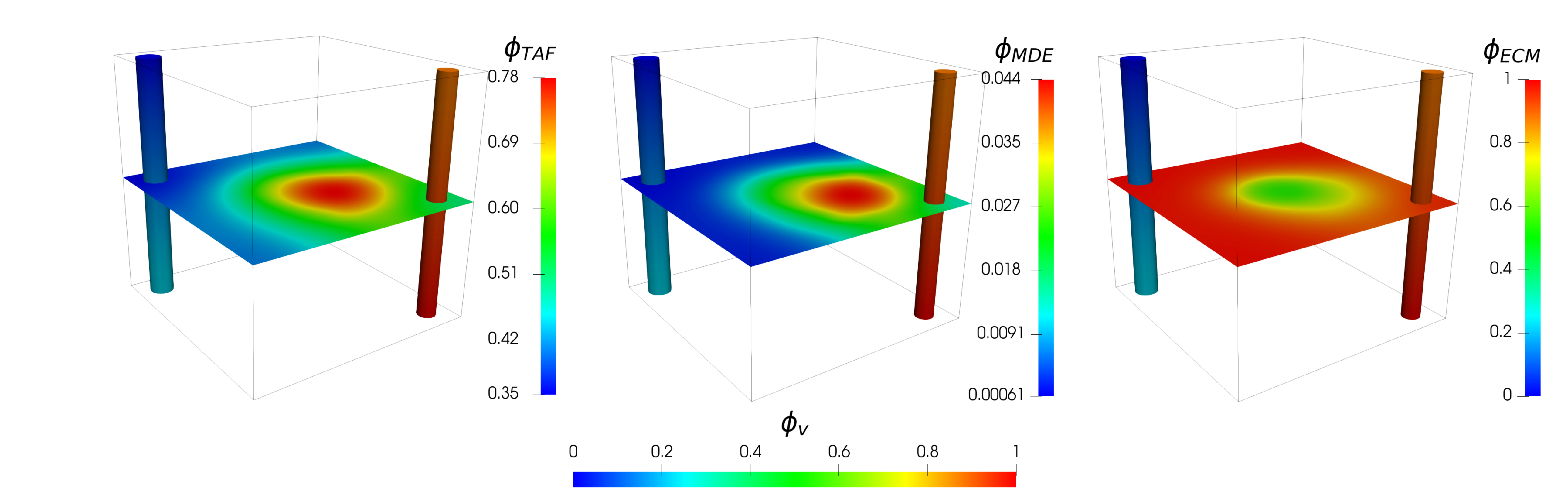}
	\end{center}
	\vspace{-0.4cm}
	\caption{\label{fig:two-vessels-taf} Plots of $\taf$ (left), $\mde$ (middle), and $\ecm$ (right) at time $t=5$ in the $z=1$ plane of the domain $\Omega$. The colors (horizontal color bar) on the vessels show the transport of the 1D nutrient ($\phi_v$). The production of TAF and MDE is maximal where the hypoxic tumor phase is located. The decay of ECM happens in the regions where MDE is produced.}
\end{figure}

Plots of the tumor species $\phi_T, \phi_P, \phi_H$ at the $z=1$ plane in the domain $\Omega$ at time points $t\in \{3,4,5\}$ are shown in \cref{fig:two-vessels-vars}. In \cref{fig:two-vessels-line-plot}, the tumor phases are shown along a one-dimensional line. It can be observed that the tumor separates into its three phases and moves towards the nutrient-rich regions of the domain. Moreover, the contour lines of the total tumor and its phases is presented at different times in \cref{fig:two-vessels-multi-contour}. We see that tumor is growing towards the artery. As expected, the proliferation is higher near the artery. In \cref{fig:two-vessels-taf}, plots of TAF, MDE, ECM at the time point $t=5$ in the $z=1$ plane of $\Omega$ are shown.  \vspace{-.3cm}

\vspace{-.5cm}

\begin{figure}[H]
	\begin{center}		
		\includegraphics[width=0.7\textwidth]{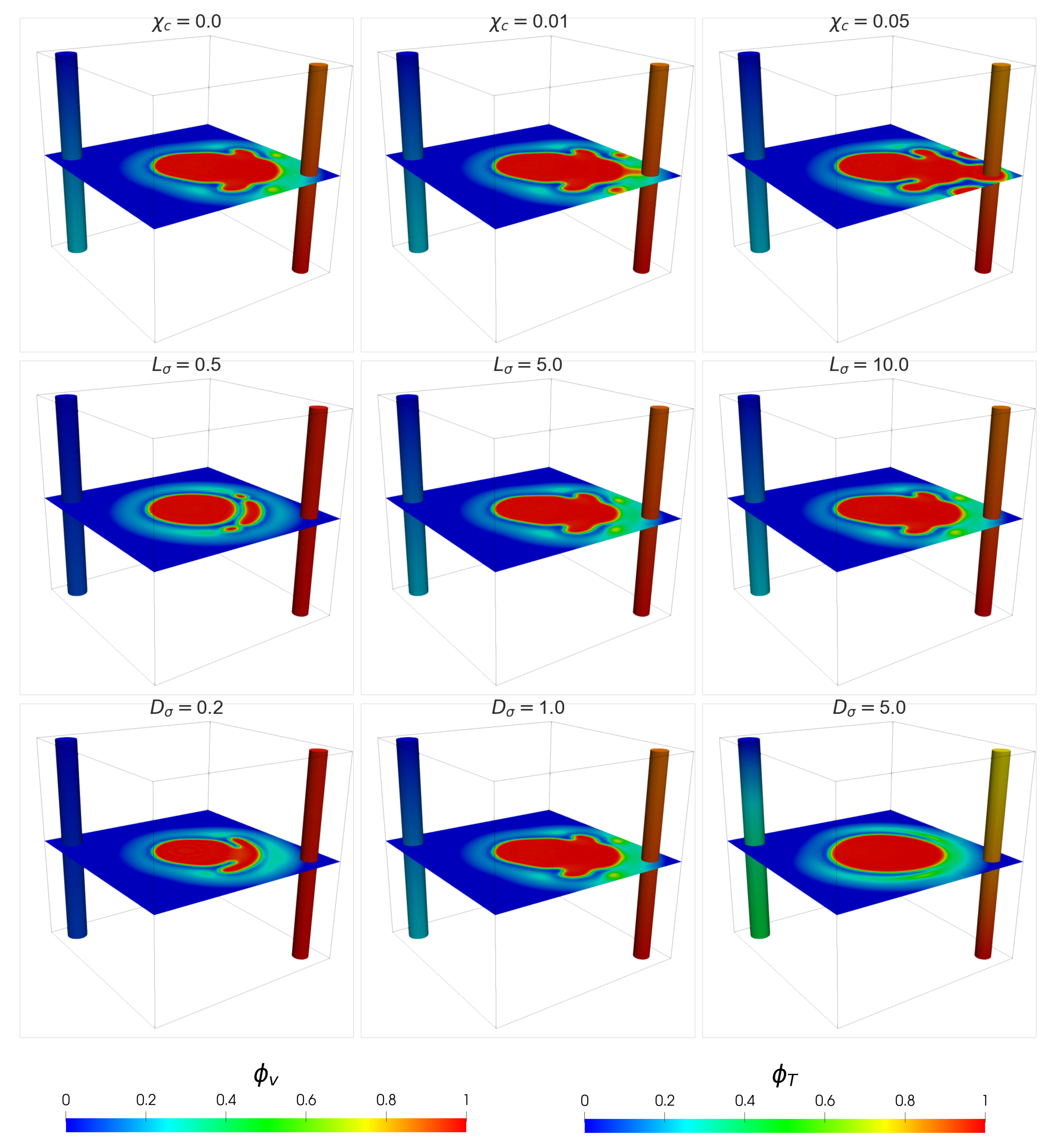}
	\end{center}
	\vspace{-0.4cm}
	\caption{\label{fig:two-vessels-chi} \label{fig:two-vessels-Ls} \label{fig:two-vessels-Dsigma} Effect of the chemotactic constant $\chi_c$ (top) for $\chi_c \in \{0,0.01,0.05\}$ (left to right), permeability of the vessel wall $L_\sigma$ (middle) for $L_\sigma \in \{0.5,5,10\}$ (left to right), and diffusivity constant $D_\sigma$ (bottom) for $D_\sigma \in \{0.2,1,5\}$ (left to right) on the growth of the tumor volume. 
		The color bar for 1D nutrients (left) and total tumor volume fraction (right) in included at the bottom.}
\end{figure}

\cref{fig:two-vessels-Dsigma}  contains simulation results for the different values of parameters at $t=T=5$. Among the numerous model parameters, we focused on the chemotactic constant $\chi_c$, mobility $M_P$, proliferation rate $\lambda_P$, nutrient diffusion coefficient $D_\sigma$, and permeability constant $L_\sigma$. We vary one of these parameters while keeping other parameters fixed to their respective values listed in \cref{tab:params}.

It can be observed that all selected parameters strongly affect the tumor growth. Except for the parameter $D_\sigma$, the larger the remaining parameters, the faster the tumor cells move away from the vein and towards the nutrient rich artery. This means that for the chosen parameter values, the fluxes $J_\alpha,\;\alpha \in \CH$
given by \cref{Eq:Flux}, dominate the corresponding convective terms in \cref{Eq:Model3D} so that the tumor cells can move against the velocity field. It can be stated that for these parameter choices, the model simulates the migration of  tumor cells towards the nutrient sources in the vicinity. 

\begin{figure}[H]
	\begin{center}		
		\includegraphics[width=.7\textwidth]{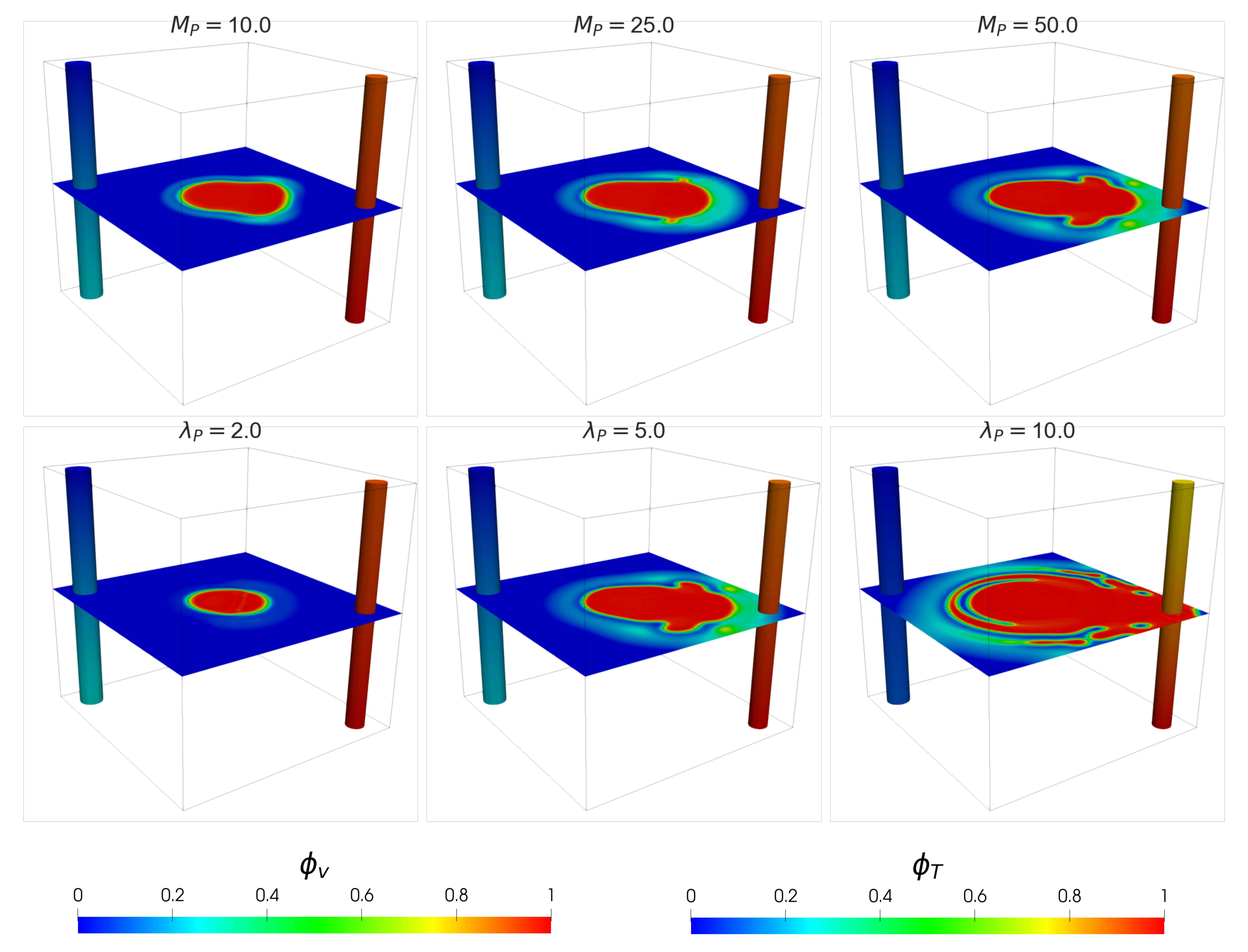}
	\end{center}
	\vspace{-0.4cm}
	\caption{\label{fig:two-vessels-mob}\label{fig:two-vessels-lamp} Effect of the mobility $M_P$ (top) for $M_P \in \{10,25,50\}$ (left to right) and proliferation rate $\lambda_P$ (bottom) for $\lambda_P \in \left\{2,5,10 \right\}$ (left to right) on the growth of the tumor volume. 
		The color bar for 1D nutrients (left) and total tumor volume fraction (right) in included at the bottom.}
\end{figure}

\begin{figure}[H]
	\begin{center}		
		\includegraphics[width=0.66\textwidth]{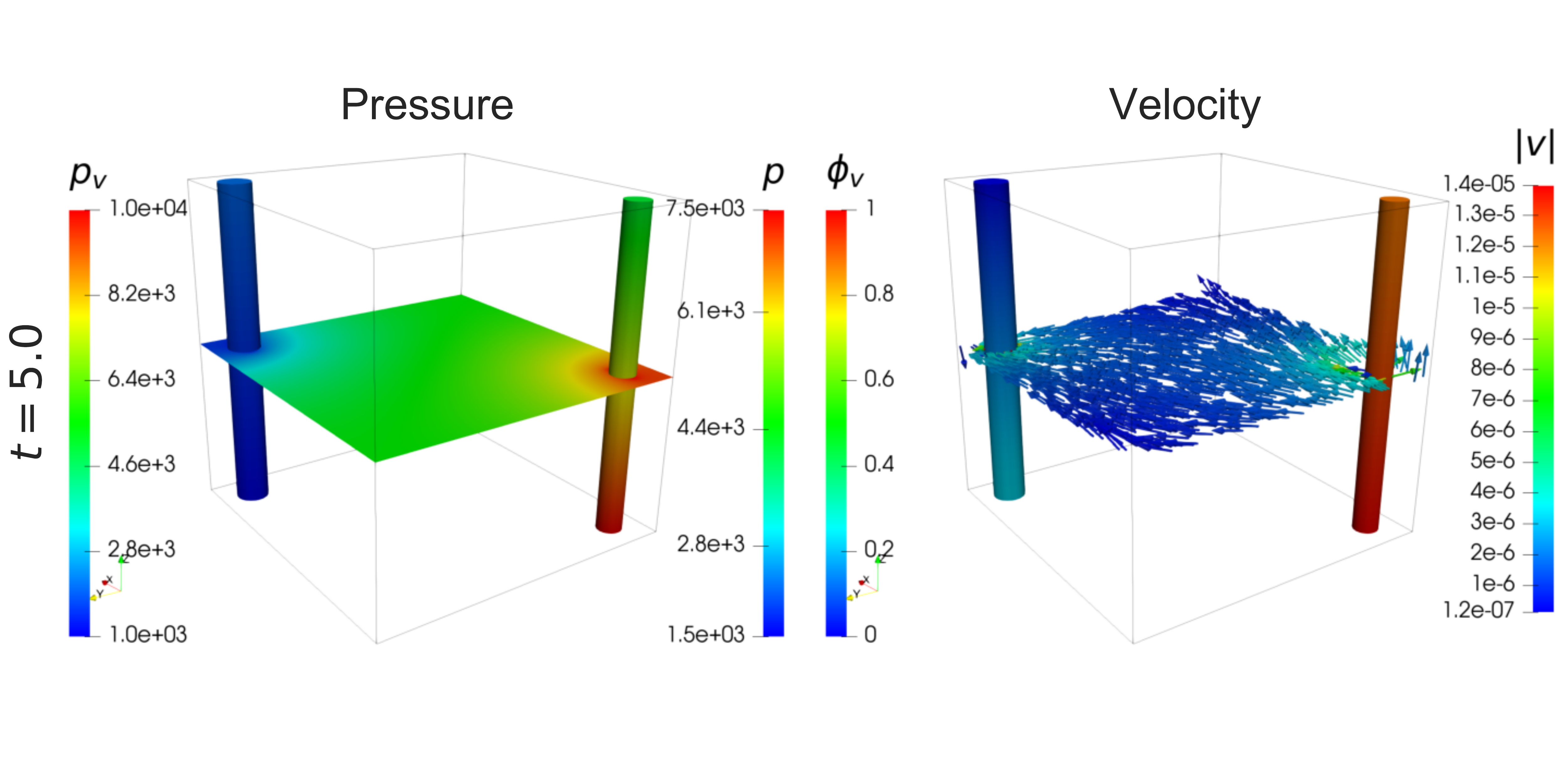}
	\end{center}
	\vspace{-1cm}
	\caption{\label{fig:two-vessels-pres-vel} Pressure (left) and velocity field (right) in a plane perpendicular to the $z$-axis $\left(\text{at } z = 1 \right)$ of the domain $\Omega$. The artery is located in the right corner with a pressure decay from $10000$ to $5000$. The vein is located in the left corner with a pressure decay from $2000$ to $1000$. The velocity field induced by the pressure distribution is directed from the artery toward the vein.}
\end{figure}
\cref{fig:two-vessels-pres-vel} shows the pressure distribution in the vessels as well as the tissue pressure and velocity field within a plane that is perpendicular to the $z$-axis and located at $z=1$. The tissue pressure ranges between $1500$ and $7500$, which means that it is bounded by the extreme pressures in the vascular system. Furthermore, a gradient in the tissue pressure can be detected pointing from the artery to the vein. As a result, the velocity field is orientated from the artery to the vein.

\subsection{Tumor surrounded by a network} 
\label{Sec:SimulationNetwork}
\begin{figure}[H]
	\begin{center}		
		\includegraphics[width=0.45\textwidth]{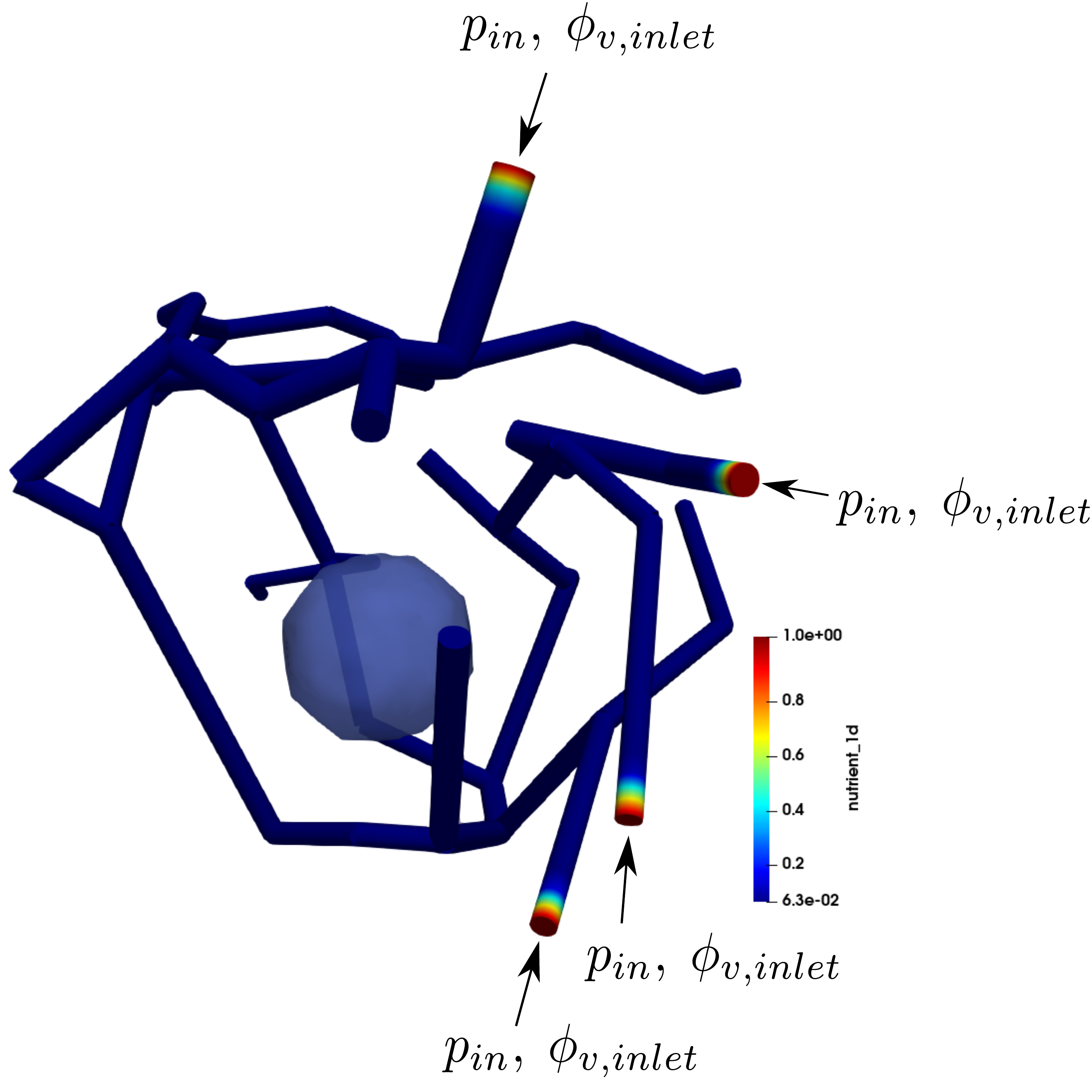}
	\end{center}
	\caption{\label{fig:BoundariesNetwork} Outline of the scaled blood vessel network with an initial tumor core. The tumor core is represented by a contour surface with respect to $\phi_T$ (at $\phi_T = 0.1$). At the four inlets, indicated by an arrow, nutrients are injected i.e. at these boundaries, we set $\phi_v = \phi_{v,inlet} =1$ and $p_v = p_{in}$.}
\end{figure}

In the second subsection, we consider a small capillary network given by the data in \cite{secomb2000theoretical}. To keep the same computational domain as in the previous subsection, the network is scaled such that it fits into $\Omega = \left(0,2\right)^3$. After scaling, the resulting network has maximum, minimum, and mean vessel radius $0.0613$, $0.0307$, $0.0418$ respectively. At the inlets that are marked by an arrow, see \cref{fig:BoundariesNetwork}, we prescribe the pressure $p_{in}=25000$, while for all the other inlets, we use $p_{out}=10000$ as a boundary value. \blue{Again, we choose the boundary values synthetically in order to ensure that the transport processes are visible in the simulations.} Further, the boundary condition for $\phi_v$ is given by $\phi_v = \phi_{v,inlet} = 1$ if it holds $p_v = p_{in}$, and free outflow boundary if $p_v = p_{out}$.

Contrary to the previous subsection, the spherical tumor core has a radius of $0.25$ and the center $(1.3, 0.9, 0.7)$. The same model parameters are employed. The domain $\Omega$ is discretized using cubic elements of mesh size $0.025$ and final time and time step of the simulation are $T = 5$ and $\Delta t = 0.025$.

\begin{figure}[H]
\begin{center}		
\includegraphics[width=.8\textwidth]{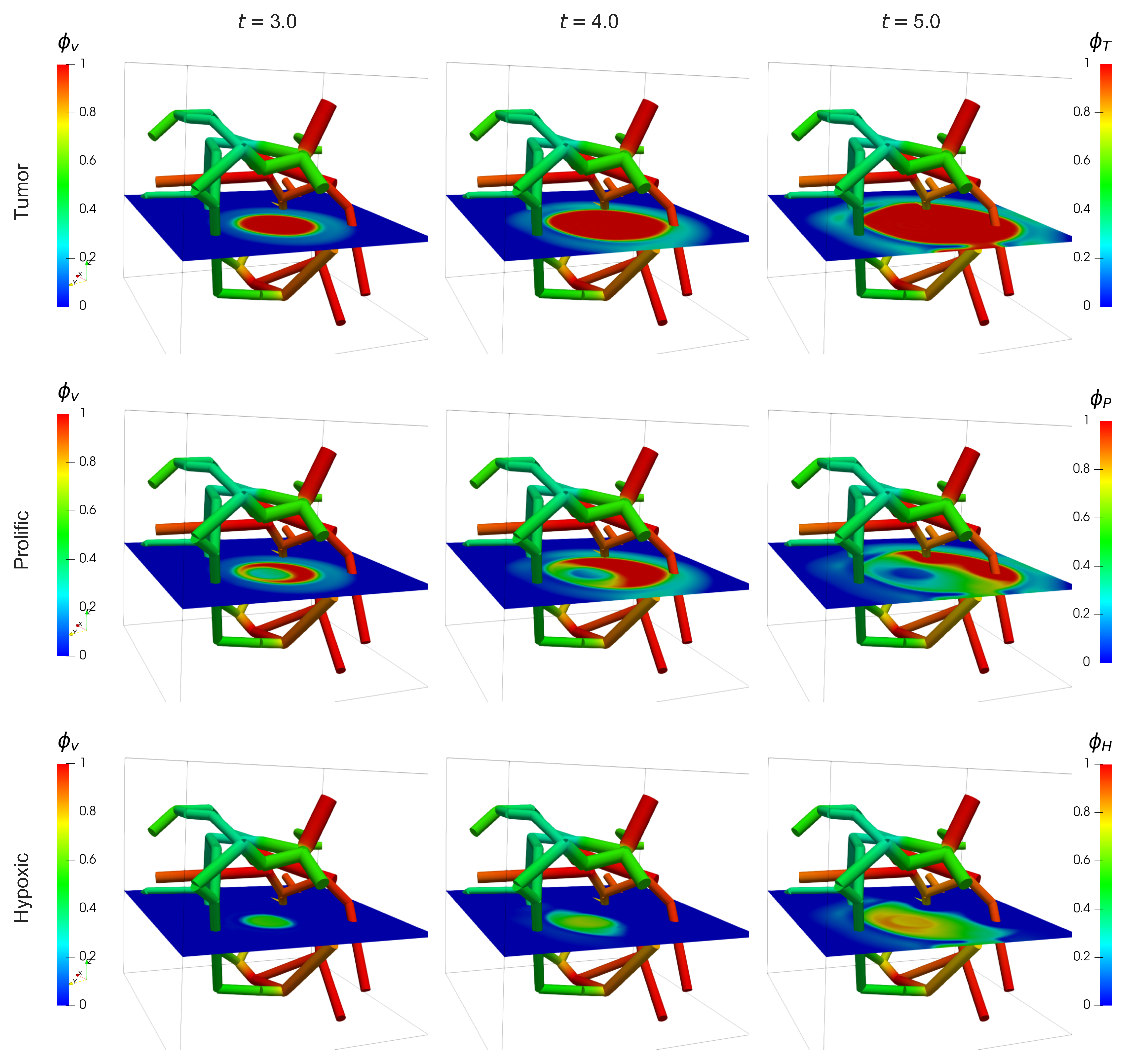}
\end{center}
\vspace{-0cm}
\caption{\label{fig:complex-network-t2-times} Distribution of the tumor cell volume fraction $\phi_T$ (top), prolific cell volume fraction $\phi_P$ (middle), and hypoxic cell volume fraction $\phi_H$ (bottom)  for $t \in \left\{3,4,5\right\}$. The tumor cells migrate towards to nutrient rich vessels.}
\end{figure}

\begin{figure}[H]
	\begin{center}		
		\includegraphics[width=.28\textwidth]{color1.pdf} \\
		\includegraphics[page=1,width=.31\textwidth]{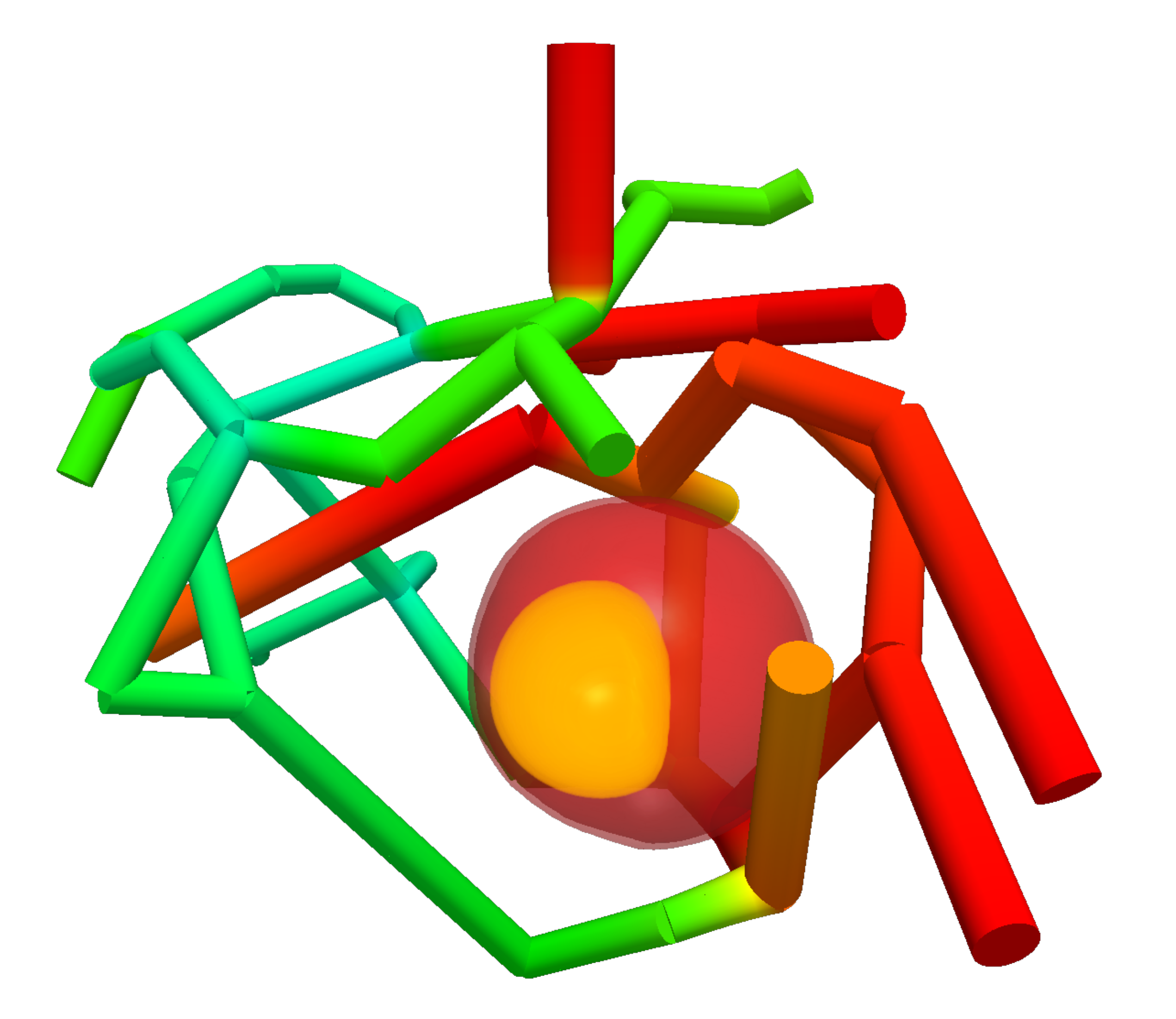}
		\includegraphics[page=2,width=.31\textwidth]{network3.pdf}
		\includegraphics[page=3,width=.31\textwidth]{network3.pdf}
	\end{center}
	\caption{\label{fig:complex-network-t2-times-contour-multi}
		Evolution of contour plots of the tumor volume fraction $\phi_T$  at the values $0.8$ (light red) and $0.95$ (red), and of the hypoxic phase $\phi_H$ at $0.35$ at times $t \in \{3,3.5,4\}$ (left to right). The tumor growth is directed towards the nutrient rich vessels. }
\end{figure}

\begin{figure}[H]
	\begin{center}		
		\includegraphics[width=.87\textwidth]{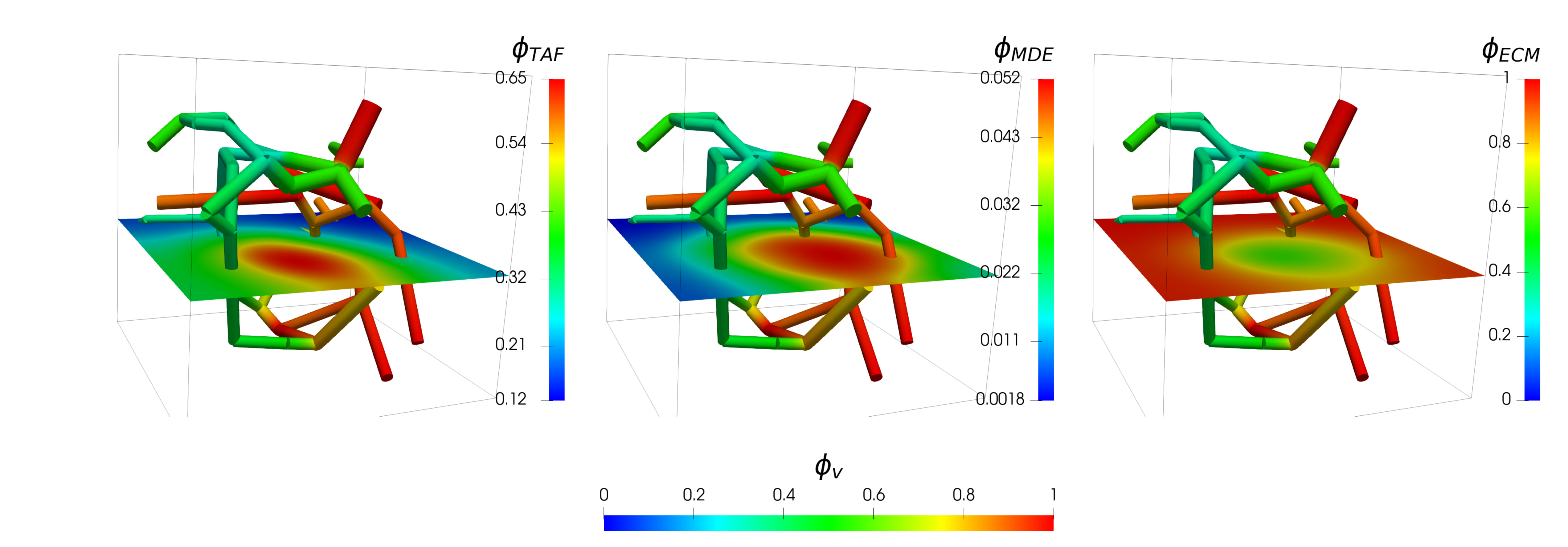}
	\end{center}
	\vspace{-0.4cm}
	\caption{\label{fig:complex-network-taf} Plots of $\taf$ (left), $\mde$ (middle), and $\ecm$ (right) at time $t=5$ in the $z=0.8$ plane of the domain $\Omega$. The colors (horizontal color bar) on the vessels show the transport of the 1D nutrient ($\phi_v$). As in the case of two-vessels setting, the production of TAF and MDE is maximal where the hypoxic region is located.}
\end{figure}

 In \cref{fig:complex-network-t2-times}, the tumor cell volume fraction $\phi_T$, prolific cell volume fraction $\phi_P$, and hypoxic cell volume fraction $\phi_H$ are shown at $z=0.8$ plane and at time points $t \in \left\{3,4,5\right\}$. Finally in \cref{fig:complex-network-t2-times-contour-multi}, the contour plots for $\phi_T = 0.8$ and $\phi_T=0.95$ are presented. Further, the hypoxic phase is shown inside the tumor. In \cref{fig:complex-network-taf}, plots of TAF, MDE, ECM at $t=5$ in $z=0.8$ plane are shown. 

The behavior of the tumor cells is similar to the two-vessel scenario. It seems that for the given parameter set the tumor cells are attracted by the nutrient rich blood vessels of the network. As can be observed in \cref{fig:complex-network-t2-times} (last row), the chemical potential of the tumor exhibits high gradients at the interface between tumor and healthy tissue. Therefore, the corresponding flux of the chemical potential given by \cref{Eq:Flux} is potentially high at this location. As a result the tumor cells are pulled towards the interface between cancerous and healthy tissue. Apparently, the flux is particularly high near the nutrient rich vessels such that the tumor cells move preferably towards the nutrient rich vessels.

\begin{figure}[H]
	\begin{center}	
		\includegraphics[width=.75\textwidth,clip,trim={0 0 0 95}]{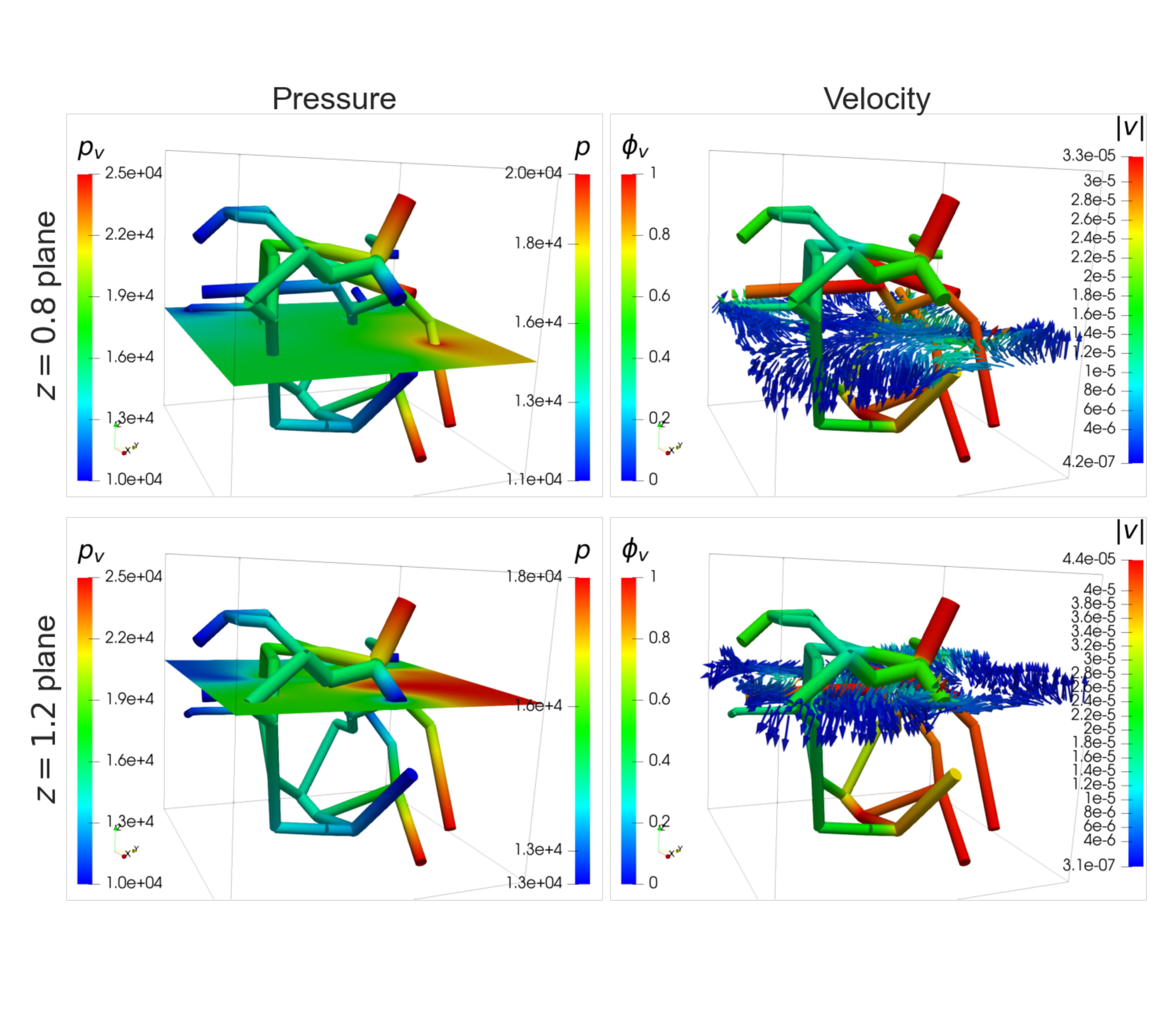}
	\end{center}
	\vspace{-1cm}
	\caption{\label{fig:complex-network-t2-pres-vel} Plots of pressure (left) and velocity field (right) in the $z=0.8$ (top) and $z=1.2$ (bottom) planes in the domain $\Omega$. The plots correspond to the simulation time $t=5$. As in the two-vessel case, the velocity field is pointing from the high pressure region to the low pressure region.}
\end{figure}

\cref{fig:complex-network-t2-pres-vel} shows the tissue pressure and the corresponding velocity fields in $z=0.8$ plane. Just as in the two-vessel scenario, the pressure distribution induces a velocity field that goes from the high pressure region to the low pressure region.

\pagebreak
\section{Summary and outlook}
In this work, we have presented a 3D-1D coupled multispecies model for tumor growth including the influence of nutrient transport in a vascular system that is located in the vicinity of a solid tumor. Flow and transport of nutrients within the vascular system are governed by one-dimensional partial differential equations. The corresponding flow and transport processes in the healthy and tumor tissue are based on Darcy's law as well as a standard convection-diffusion equation. Coupling of the three-dimensional equations with their one-dimensional counterparts is done via filtration laws and source terms. In the source terms of the three-dimensional partial differential equations, Dirac measures occur. They are concentrated on the vessel surfaces of the vascular system, since there the exchange processes between the tissue and the vascular system take place. The remaining three-dimensional equations governing the distribution of the tumor cells are of Cahn--Hilliard type. The evolution of matrix degrading enzymes and the tumor angiogensis factor are modeled by convection-diffusion equations. Lastly, the extracellular matrix density is described by an abstract ordinary differential equation.

The centerpiece of our work is a mathematical analysis of this model with a focus on the existence of solutions. We have shown the existence of weak solutions. Our proof is based on the Faedo--Galerkin method. Thereby, the system of partial differential equations is semi-discretized in space and reduced to a system of ordinary differential equations. Using the Cauchy--Peano theorem we show that the system of ordinary differential equations exhibits a solution. In a next step, the existence of weak solutions with respect to the partial differential equations is derived by means of the Banach--Alaoglu theorem. Finally, we present some simulation results for two different settings, illustrating the performance of our model. Our simulation results indicate that the tumor cells sense the vessels with an increased nutrient concentration and move towards them. Furthermore, the impact of several model parameters on the solution variables is discussed.

Among extensions and applications of the models described here are the simulation and optimal control of
chemotherapy drug and radiation as well as modeling of the onset of metastasis. Simulation of these protocols and phenomena represent challenging goals for future work.


\section*{Acknowledgements}
We gratefully acknowledge the support of the German Science Foundation (DFG) for funding part of this work through grant WO 671/11-1. The work of PKJ and JTO was supported by the U.S. Department of Energy, Office of Science, Office of Advanced Scientific Computing Research, Mathematical Multifaceted Integrated Capability 
Centers (MMICCS), under Award Number DE-SC0019303.

\bibliography{literature}

\begin{thebibliography}{10}

\bibitem{alt2016linear}
{\sc H.~W. Alt}, {\em Linear Functional Analysis: An Application-Oriented
  Introduction}, Springer Science \& Business Media, 2016.

\bibitem{anderson1998continuous}
{\sc A.~R. Anderson and M.~A.~J. Chaplain}, {\em Continuous and discrete
  mathematical models of tumor-induced angiogenesis}, Bulletin of Mathematical
  Biology, 60 (1998), pp.~857--899.

\bibitem{boyer2012mathematical}
{\sc F.~Boyer and P.~Fabrie}, {\em Mathematical Tools for the Study of the
  Incompressible {N}avier--{S}tokes Equations and Related Models}, Springer
  Science \& Business Media, 2012.

\bibitem{brezis2010functional}
{\sc H.~Brezis}, {\em Functional Analysis, {S}obolev Spaces and Partial
  Differential Equations}, Springer Science \& Business Media, 2010.

\bibitem{byrne2003modelling}
{\sc H.~Byrne and L.~Preziosi}, {\em Modelling solid tumour growth using the
  theory of mixtures}, Mathematical Medicine and Biology: A Journal of the IMA,
  20 (2003), pp.~341--366.

\bibitem{chaplain2011mathematical}
{\sc M.~A. Chaplain, M.~Lachowicz, Z.~Szyma{\'n}ska, and D.~Wrzosek}, {\em
  Mathematical modelling of cancer invasion: the importance of cell--cell
  adhesion and cell--matrix adhesion}, Mathematical Models and Methods in
  Applied Sciences, 21 (2011), pp.~719--743.

\bibitem{chaplain2005mathematical}
{\sc M.~A. Chaplain and G.~Lolas}, {\em Mathematical modelling of cancer cell
  invasion of tissue: The role of the urokinase plasminogen activation system},
  Mathematical Models and Methods in Applied Sciences, 15 (2005),
  pp.~1685--1734.

\bibitem{cherfils2011cahn}
{\sc L.~Cherfils, A.~Miranville, and S.~Zelik}, {\em The cahn-hilliard equation
  with logarithmic potentials}, Milan Journal of Mathematics, 79 (2011),
  pp.~561--596.

\bibitem{cristini2009nonlinear}
{\sc V.~Cristini, X.~Li, J.~S. Lowengrub, and S.~M. Wise}, {\em Nonlinear
  simulations of solid tumor growth using a mixture model: invasion and
  branching}, Journal of Mathematical Biology, 58 (2009).

\bibitem{cristini2010multiscale}
{\sc V.~Cristini and J.~Lowengrub}, {\em Multiscale Modeling of Cancer: An
  Integrated Experimental and Mathematical Modeling Approach}, Cambridge
  University Press, 2010.

\bibitem{dangelo2008coupling}
{\sc C.~D'ANGELO and A.~Quarteroni}, {\em On the coupling of 1d and 3d
  diffusion-reaction equations: application to tissue perfusion problems},
  Mathematical Models and Methods in Applied Sciences, 18 (2008),
  pp.~1481--1504.

\bibitem{diestel1977vector}
{\sc J.~Diestel and J.~Uhl}, {\em Vector Measures}, American Mathematical
  Society, 1977.

\bibitem{engwer2017structured}
{\sc C.~Engwer, C.~Stinner, and C.~Surulescu}, {\em On a structured multiscale
  model for acid-mediated tumor invasion: The effects of adhesion and
  proliferation}, Mathematical Models and Methods in Applied Sciences, 27
  (2017), pp.~1355--1390.

\bibitem{evans2010partial}
{\sc L.~C. Evans}, {\em Partial {D}ifferential {E}quations}, American
  Mathematical Society, 2010.

\bibitem{frieboes2010three}
{\sc H.~B. Frieboes, F.~Jin, Y.-L. Chuang, S.~M. Wise, J.~S. Lowengrub, and
  V.~Cristini}, {\em Three-dimensional multispecies nonlinear tumor growth --
  {II}: {T}umor invasion and angiogenesis}, Journal of Theoretical Biology, 264
  (2010), pp.~1254--1278.

\bibitem{frigeri2015diffuse}
{\sc S.~Frigeri, M.~Grasselli, and E.~Rocca}, {\em On a diffuse interface model
  of tumour growth}, European Journal of Applied Mathematics, 26 (2015),
  pp.~215--243.

\bibitem{frigeri2018on}
{\sc S.~Frigeri, K.~F. Lam, E.~Rocca, and G.~Schimperna}, {\em On a
  multi-species cahn--hilliard--darcy tumor grwoth model with singular
  potentials}, Communications in Mathematical Sciences, 16 (2018),
  pp.~821--856.

\bibitem{fritz2020part2}
{\sc M.~Fritz, P.~K. Jha, T.~K\"oppl, J.~T. Oden, A.~Wagner, and B.~Wohlmuth},
  {\em Modeling and simulation of vascular tumors embedded in evolving
  capillary networks}, arXiv preprint arXiv:2001.10183,  (2020).

\bibitem{fritz2019local}
{\sc M.~Fritz, E.~Lima, V.~Nikolic, J.~T. Oden, and B.~Wohlmuth}, {\em Local
  and nonlocal phase-field models of tumor growth and invasion due to {ECM}
  degradation}, Mathematical Models and Methods in Applied Sciences, 29 (2019),
  pp.~2433--2468.

\bibitem{fritz2018unsteady}
{\sc M.~Fritz, E.~Lima, J.~T. Oden, and B.~Wohlmuth}, {\em On the unsteady
  {D}arcy--{F}orchheimer--{B}rinkman equation in local and nonlocal tumor
  growth models}, Mathematical Models and Methods in Applied Sciences, 29
  (2019), pp.~1691--1731.

\bibitem{garcke2016global}
{\sc H.~Garcke and K.~F. Lam}, {\em Global weak solutions and asymptotic limits
  of a {C}ahn--{H}illiard--{D}arcy system modelling tumour growth}, AIMS
  Mathematics, 1 (2016), pp.~318--360.

\bibitem{garcke2017well}
\leavevmode\vrule height 2pt depth -1.6pt width 23pt, {\em Well-posedness of a
  {C}ahn--{H}illiard system modelling tumour growth with chemotaxis and active
  transport}, European Journal of Applied Mathematics, 28 (2017), pp.~284--316.

\bibitem{garcke2018multiphase}
{\sc H.~Garcke, K.~F. Lam, R.~N{\"u}rnberg, and E.~Sitka}, {\em A multiphase
  {C}ahn--{H}illiard--{D}arcy model for tumour growth with necrosis},
  Mathematical Models and Methods in Applied Sciences, 28 (2018), pp.~525--577.

\bibitem{ginzburg1963frictional}
{\sc B.~Ginzburg and A.~Katchalsky}, {\em The frictional coefficients of the
  flows of non-electrolytes through artificial membranes}, The Journal of
  General Physiology, 47 (1963), pp.~403--418.

\bibitem{grosse2013sobolev}
{\sc N.~Gro{\ss}e and C.~Schneider}, {\em Sobolev spaces on {R}iemannian
  manifolds with bounded geometry: {G}eneral coordinates and traces},
  Mathematische Nachrichten, 286 (2013), pp.~1586--1613.

\bibitem{hanahan2011hallmarks}
{\sc D.~Hanahan and R.~A. Weinberg}, {\em Hallmarks of cancer: {T}he next
  generation}, Cell, 144 (2011), pp.~646--674.

\bibitem{hawkins2012numerical}
{\sc A.~Hawkins-Daarud, K.~G. van~der Zee, and J.~T. Oden}, {\em Numerical
  simulation of a thermodynamically consistent four-species tumor growth
  model}, International Journal for Numerical Methods in Biomedical
  Engineering, 28 (2012), pp.~3--24.

\bibitem{hillen2013convergence}
{\sc T.~Hillen, K.~J. Painter, and M.~Winkler}, {\em Convergence of a cancer
  invasion model to a logistic chemotaxis model}, Mathematical Models and
  Methods in Applied Sciences, 23 (2013), pp.~165--198.

\bibitem{jiang2015well}
{\sc J.~Jiang, H.~Wu, and S.~Zheng}, {\em Well-posedness and long-time behavior
  of a non-autonomous {C}ahn--{H}illiard--{D}arcy system with mass source
  modeling tumor growth}, Journal of Differential Equations, 259 (2015),
  pp.~3032--3077.

\bibitem{kim2009existence}
{\sc H.~Kim}, {\em Existence and regularity of very weak solutions of the
  stationary {N}avier--{S}tokes equations}, Archive for Rational Mechanics and
  Analysis, 193 (2009), pp.~117--152.

\bibitem{koppl20203d}
{\sc T.~K{\"o}ppl, E.~Vidotto, and B.~Wohlmuth}, {\em A {3D}-{1D} coupled blood
  flow and oxygen transport model to generate microvascular networks},
  International journal for numerical methods in biomedical engineering, 36
  (2020), p.~e3386.

\bibitem{koppl2018mathematical}
{\sc T.~K{\"o}ppl, E.~Vidotto, B.~Wohlmuth, and P.~Zunino}, {\em Mathematical
  modeling, analysis and numerical approximation of second-order elliptic
  problems with inclusions}, Mathematical Models and Methods in Applied
  Sciences, 28 (2018), pp.~953--978.

\bibitem{kremheller2019approach}
{\sc J.~Kremheller, A.-T. Vuong, B.~A. Schrefler, and W.~A. Wall}, {\em An
  approach for vascular tumor growth based on a hybrid embedded/homogenized
  treatment of the vasculature within a multiphase porous medium model},
  International Journal for Numerical Methods in Biomedical Engineering, 35
  (2019), p.~e3253.

\bibitem{leoni2007necessary}
{\sc G.~Leoni and M.~Morini}, {\em Necessary and sufficient conditions for the
  chain rule in ${W}_{\textup{loc}}^{1,1}(\mathbb{R}^n;\mathbb{R}^d)$ and
  ${BV}_{\textup{loc}}(\mathbb{R}^n;\mathbb{R}^d)$}, Journal of the European
  Mathematical Society, 9 (2007), pp.~219--252.

\bibitem{lima2015analysis}
{\sc E.~Lima, R.~C. Almeida, and J.~T. Oden}, {\em Analysis and numerical
  solution of stochastic phase-field models of tumor growth}, Numerical Methods
  for Partial Differential Equations, 31 (2015), pp.~552--574.

\bibitem{lima2014hybrid}
{\sc E.~Lima, J.~T. Oden, and R.~Almeida}, {\em A hybrid ten-species
  phase-field model of tumor growth}, Mathematical Models and Methods in
  Applied Sciences, 24 (2014), pp.~2569--2599.

\bibitem{lima2017selection}
{\sc E.~Lima, J.~T. Oden, B.~Wohlmuth, A.~Shahmoradi, D.~Hormuth~II,
  T.~Yankeelov, L.~Scarabosio, and T.~Horger}, {\em Selection and validation of
  predictive models of radiation effects on tumor growth based on noninvasive
  imaging data}, Computer Methods in Applied Mechanics and Engineering, 327
  (2017), pp.~277--305.

\bibitem{lions2012non}
{\sc J.~L. Lions and E.~Magenes}, {\em Non-Homogeneous Boundary Value Problems
  and Applications}, Springer Science \& Business Media, 2012.

\bibitem{lowengrub2013analysis}
{\sc J.~Lowengrub, E.~Titi, and K.~Zhao}, {\em Analysis of a mixture model of
  tumor growth}, European Journal of Applied Mathematics, 24 (2013),
  pp.~691--734.

\bibitem{murat2003chain}
{\sc F.~Murat and C.~Trombetti}, {\em A chain rule formula for the composition
  of a vector-valued function by a piecewise smooth function}, Bollettino
  dell'Unione Matematica Italiana, 6 (2003), pp.~581--595.

\bibitem{nargis2016effects}
{\sc N.~Nargis and R.~Aldredge}, {\em Effects of matrix metalloproteinase on
  tumour growth and morphology via haptotaxis}, J. Bioengineer. \& Biomedical
  Sci., 6 (2016).

\bibitem{oden2010general}
{\sc J.~T. Oden, A.~Hawkins, and S.~Prudhomme}, {\em General diffuse-interface
  theories and an approach to predictive tumor growth modeling}, Mathematical
  Models and Methods in Applied Sciences, 20 (2010), pp.~477--517.

\bibitem{oden2016toward}
{\sc J.~T. Oden, E.~Lima, R.~C. Almeida, Y.~Feng, M.~N. Rylander, D.~Fuentes,
  D.~Faghihi, M.~M. Rahman, M.~DeWitt, M.~Gadde, et~al.}, {\em Toward
  predictive multiscale modeling of vascular tumor growth}, Archives of
  Computational Methods in Engineering, 23 (2016), pp.~735--779.

\bibitem{reichold2009vascular}
{\sc J.~Reichold, M.~Stampanoni, A.~L. Keller, A.~Buck, P.~Jenny, and
  B.~Weber}, {\em Vascular graph model to simulate the cerebral blood flow in
  realistic vascular networks}, Journal of Cerebral Blood Flow \& Metabolism,
  29 (2009), pp.~1429--1443.

\bibitem{roubicek}
{\sc T.~Roub{\'\i}{\v{c}}ek}, {\em Nonlinear {P}artial {D}ifferential
  {E}quations with {A}pplications}, Springer Science \& Business Media, 2013.

\bibitem{saedpanah2014well}
{\sc F.~Saedpanah}, {\em Well-posedness of an integro-differential equation
  with positive type kernels modeling fractional order viscoelasticity},
  European Journal of Mechanics-A/Solids, 44 (2014), pp.~201--211.

\bibitem{santagiuliana2016simulation}
{\sc R.~Santagiuliana, M.~Ferrari, and B.~Schrefler}, {\em Simulation of
  angiogenesis in a multiphase tumor growth model}, Computer Methods in Applied
  Mechanics and Engineering, 304 (2016), pp.~197--216.

\bibitem{santagiuliana2019coupling}
{\sc R.~Santagiuliana, M.~Milosevic, B.~Milicevic, G.~Scium{\`e}, V.~Simic,
  A.~Ziemys, M.~Kojic, and B.~A. Schrefler}, {\em Coupling tumor growth and bio
  distribution models}, Biomedical microdevices, 21 (2019), p.~33.

\bibitem{secomb2000theoretical}
{\sc T.~Secomb, R.~Hsu, N.~Beamer, and B.~Coull}, {\em Theoretical simulation
  of oxygen transport to brain by networks of microvessels: effects of oxygen
  supply and demand on tissue hypoxia}, Microcirculation, 7 (2000),
  pp.~237--247.

\bibitem{sfakianakis2020hybrid}
{\sc N.~Sfakianakis, A.~Madzvamuse, and M.~A. Chaplain}, {\em A hybrid
  multiscale model for cancer invasion of the extracellular matrix}, Multiscale
  Modeling \& Simulation, 18 (2020), pp.~824--850.

\bibitem{simon1986compact}
{\sc J.~Simon}, {\em Compact sets in the space ${L}^p(0, {T}; {B})$}, Annali di
  Matematica pura ed applicata, 146 (1986), pp.~65--96.

\bibitem{tao2011global}
{\sc Y.~Tao}, {\em Global existence for a haptotaxis model of cancer invasion
  with tissue remodeling}, Nonlinear Analysis: Real World Applications, 12
  (2011), pp.~418--435.

\bibitem{tao2011chemotaxis}
{\sc Y.~Tao and M.~Winkler}, {\em A chemotaxis-haptotaxis model: the roles of
  nonlinear diffusion and logistic source}, SIAM Journal on Mathematical
  Analysis, 43 (2011), pp.~685--704.

\bibitem{tao2014energy}
\leavevmode\vrule height 2pt depth -1.6pt width 23pt, {\em Energy-type
  estimates and global solvability in a two-dimensional chemotaxis--haptotaxis
  model with remodeling of non-diffusible attractant}, Journal of Differential
  Equations, 257 (2014), pp.~784--815.

\bibitem{vidotto2019hybrid}
{\sc E.~Vidotto, T.~Koch, T.~K\"oppl, R.~Helmig, and B.~Wohlmuth}, {\em Hybrid
  models for simulating blood flow in microvascular networks}, Multiscale
  Modeling \& Simulation, 17 (2019), pp.~1076--1102.

\bibitem{walter1998ordinary}
{\sc W.~Walter}, {\em Ordinary Differential Equations}, Springer Science \&
  Business Media, 1998.

\bibitem{wise2008three}
{\sc S.~M. Wise, J.~S. Lowengrub, H.~B. Frieboes, and V.~Cristini}, {\em
  Three-dimensional multispecies nonlinear tumor growth -- {I}: {M}odel and
  numerical method}, Journal of Theoretical Biology, 253 (2008), pp.~524--543.

\bibitem{xu2016mathematical}
{\sc J.~Xu, G.~Vilanova, and H.~Gomez}, {\em A mathematical model coupling
  tumor growth and angiogenesis}, PloS one, 11 (2016).

\bibitem{xu2017full}
\leavevmode\vrule height 2pt depth -1.6pt width 23pt, {\em Full-scale,
  three-dimensional simulation of early-stage tumor growth: The onset of
  malignancy}, Computer Methods in Applied Mechanics and Engineering, 314
  (2017), pp.~126--146.

\bibitem{xu2020phase}
\leavevmode\vrule height 2pt depth -1.6pt width 23pt, {\em Phase-field model of
  vascular tumor growth: Three-dimensional geometry of the vascular network and
  integration with imaging data}, Computer Methods in Applied Mechanics and
  Engineering, 359 (2020), p.~112648.

\end{thebibliography}
\bibliographystyle{siam}

\end{document}